\documentclass[a4paper,twoside]{article}

\usepackage[utf8]{inputenc}
\usepackage[T1]{fontenc}
\usepackage[margin=2.5cm]{geometry}

\usepackage[dvipsnames]{xcolor}
\usepackage{hyperref}
\hypersetup{colorlinks,linkcolor={red},citecolor={blue},urlcolor={green!70!black}}
\usepackage{enumitem}

\usepackage{amsfonts,amssymb,bm}
\usepackage{siunitx}
\usepackage{float}

\usepackage{amsmath,empheq}
\numberwithin{equation}{section}

\usepackage{amsthm}
\newtheoremstyle{mytheoremstyle}{7pt}{7pt}{\normalfont}{}{\normalfont\bfseries}{:}{.5em}{}
\theoremstyle{mytheoremstyle}
\newtheorem{definition}{Definition}[section]
\newtheorem{lemma}[definition]{Lemma}
\newtheorem{proposition}[definition]{Proposition}
\newtheorem{theorem}[definition]{Theorem}
\newtheorem{corollary}[definition]{Corollary}

\newtheorem{remark}[definition]{Remark}

\usepackage[square,authoryear]{natbib}

\newcommand{\SetFont}[1]{\mathbb{#1}}

\newcommand{\R}{\SetFont{R}}
\newcommand{\LieGroupFont}[1]{\mathsf{#1}}
\newcommand{\Diff}{\LieGroupFont{Diff}}

\newcommand{\GL}{\LieGroupFont{GL}}
\newcommand{\LieAlgebraFont}[1]{\mathfrak{#1}}
\newcommand{\gl}{\LieAlgebraFont{gl}}
\newcommand{\dif}{\bm{\mathrm{d}}}
\newcommand{\lie}{\bm{\mathrm{L}}}
\newcommand{\ins}{\bm{\mathrm{i}}}
\newcommand{\ad}{\mathrm{ad}}
\newcommand{\Ad}{\mathrm{Ad}}
\newcommand{\Id}{\mathrm{Id}}
\newcommand{\Lag}{L}
\newcommand{\Lagd}{\mathcal{L}}
\newcommand{\lag}{\ell}
\newcommand{\dive}{\operatorname{div}}
\newcommand{\grad}{\operatorname{grad}}
\newcommand{\curl}{\operatorname{curl}}
\newcommand{\Tr}{\operatorname{Tr}}
\newcommand{\Den}{\operatorname{Den}}
\newcommand{\Def}{\operatorname{Def}}
\newcommand{\Orb}{\operatorname{Orb}}
\newcommand{\smooth}[1]{\mathcal{C}^\infty({#1})}

\makeatletter
\newsavebox{\@brx}
\newcommand{\llangle}[1][]{\savebox{\@brx}{\(\m@th{#1\langle}\)}%
	\mathopen{\copy\@brx\mkern2mu\kern-0.9\wd\@brx\usebox{\@brx}}}
\newcommand{\rrangle}[1][]{\savebox{\@brx}{\(\m@th{#1\rangle}\)}%
	\mathclose{\copy\@brx\mkern2mu\kern-0.9\wd\@brx\usebox{\@brx}}}
\makeatother

\usepackage{tikz}
\usepackage{rotating,hvfloat}

\newcommand*\circled[1]{\tikz[baseline=(char.base)]{\node[shape=circle,draw,inner sep=1pt] (char) {\scriptsize #1};}}

\usepackage{mdframed}
\mdfdefinestyle{box}
{
	frametitlerule=true,
	frametitlebackgroundcolor=black!5,
}

\begin{document}

\title{\textbf{Variational discretization\\ of the Navier-Stokes-Fourier system}}
\author{Benjamin \textsc{Cou\'eraud}\\ CNRS -- LMD -- IPSL\\ \'Ecole Normale Sup\'erieure de Paris -- PSL\\ 24 rue Lhomond, 75005 Paris, France\\ \href{mailto:coueraud@lmd.ens.fr}{\texttt{coueraud@lmd.ens.fr}}\and François \textsc{Gay-Balmaz}\\ CNRS -- LMD -- IPSL\\ \'Ecole Normale Sup\'erieure de Paris -- PSL\\ 24 rue Lhomond, 75005 Paris, France\\ \href{mailto:francois.gay-balmaz@lmd.ens.fr}{\texttt{francois.gay-balmaz@lmd.ens.fr}}}
\date{\today}

\maketitle

\begin{abstract}
	This paper presents the variational discretization of the compressible Navier-Stokes-Fourier system, in which the viscosity term and the heat conduction are handled within the variational approach to nonequilibrium thermodynamics developed in \cite{GBYo2017a,GBYo2017b}, as well as its discrete counterpart whose foundations have been laid in \cite{GBYo2018}, \cite{CoGB2018} and more particularly \cite{BaGB2018} for the case of the compressible Euler equation. In a first part, we review the variational framework introduced in \cite{GBYo2017b} for the Navier-Stokes-Fourier (NSF) system in the smooth setting. In a second part, we first review the discrete exterior calculus developed in \cite{PaMuToKaMaDe2011}, which is based on discrete diffeomorphisms and therefore particularly well-adapted for the discretization of Euler-Poincaré systems whose configuration space is an infinite-dimensional Lie group. Then, the variational principle underlying the NSF system is spatially discretized through the use of this discrete exterior calculus, which yields a semi-discrete nonholonomic variational principle for the NSF system, as well as semi-discrete evolution equations. In order to avoid important technical difficulties, further treatment of the phenomenological constraint is needed. In a third part we discretize in time the spatial variational principle underlying the NSF system by extending \cite{CoGB2018}, which at last yields a nonholonomic variational integrator for the NSF system, as well as (fully) discrete evolution equations.
\end{abstract}

\tableofcontents

%---------------------------------------------------------------------------------------------------

\section*{Introduction}
\addcontentsline{toc}{section}{Introduction}

In \cite{GBYo2017a,GBYo2017b} a new variational formalism was proposed for nonequilibrium thermodynamics. This formalism is an extension of the Hamilton principle that allows the inclusion of irreversible phenomena in discrete\footnote{In this context, recall that a \emph{discrete} system is a system composed of multiple \emph{simple} systems, interconnected in some sense (see \cite{GBYo2017a} for the precise definition, but such systems will not be used in what follows).} and continuum systems, by using a nonholonomic nonlinear constraint, the so-called \emph{phenomenological constraint}, and its associated \emph{variational constraint}. In this formalism, the entropy of the system is promoted to a full dynamic variable, and to each irreversible process corresponds a \emph{thermodynamic displacement}, whose rate equals the \emph{thermodynamic affinity} of the process. Thanks to the introduction of these variables, together with the phenomenological and variational constraints, this new variational formalism yields the time evolution equations of the system in accordance with the two fundamental laws of thermodynamics, see, e.g., \citet[chapter 1]{St1974}.

Equipped with such a variational formalism, it is natural to try to devise new \emph{variational integrators} based on this approach, with the aim of developing new algorithms for the simulation of multiphysics systems with the advantages that are known to variational integrators based on Lagrangian mechanics, see \cite{MaWe2001}. Indeed, variational integrators were proved to be superior to more classical algorithms thanks to the fact that they were designed to preserve as much as possible of the geometric structures underlying the mechanical system they discretize. Some important features are that the discrete symplectic structure of the discrete Lagrangian system, as well as the discrete momenta in case when there are symmetries, are all preserved, and the discrete total energy of the system \emph{remains bounded} and oscillates around its correct value during the simulation. In \cite{GBYo2018}, variational integrators for the nonequilibrium thermodynamics of simple closed systems were developed. Because of the presence of thermal effects, the flow of the smooth equations is not symplectic anymore, but rather satisfies a generalized \emph{structure-preserving property}, which reduces to the conservation of the symplectic form if thermal effects are absent. The associated variational integrators satisfy a discrete version of this structure-preserving property. In \cite{CoGB2018}, this approach was extended to \emph{a particular type} of Euler-Poincaré systems that undergo \emph{irreversible thermodynamic processes}. We recall that \emph{classical} Euler-Poincaré systems have a Lie group as configuration space, and a Lagrangian which is invariant under the action of this group, leading to a reduced variational principle as well as reduced Euler-Lagrange equations famously known as the Euler-Poincaré equations (see \cite{HoMaRa1998}). This approach can be extended to the discrete setting and leads to a variational integrator on the Lie algebra of the configuration space through the use of a \emph{group difference map} (see \cite{MaPeSh1999} and \cite{BRMa2008}).

Classical Euler-Poincaré reduction covers a broad range of physical systems, and is particularly suitable for the geometric treatment of various types of fluids (see for instance \cite{HoMaRa1998}, \cite{HoMaRa1999}, \cite{Ho2012}). Its discrete counterpart has been used in a number of works to address the \emph{time discretization} of partial differential equations that describe the evolution of \emph{reversible} fluids (meaning inviscid, non heat conducting fluids), see \cite{PaMuToKaMaDe2011}, \cite{GMPMD2011}, and \cite{DeGaGBZe2014}. Note that paramount to this approach is the use of a discrete exterior calculus (first defined in \cite{PaMuToKaMaDe2011} and then extended to the compressible case in \cite{BaGB2018}) based on the notion of \emph{discrete diffeomorphisms}, that is used for the \emph{spatial discretization} and leveraged in the derivation of the semi-discrete and fully discrete variational principles. As expected, the obtained variational integrators derived in this way all exhibit a good long-term behavior for energy as well as conservation of a Kelvin-Noether quantity. In this paper, we are concerned with the design of such a variational integrator for the Navier-Stokes-Fourier system, that we will simply refer to the NSF system. More particularly, the NSF system will refer for us to the set of partial differential equations that govern the compressible flow on a 2D or 3D domain $M$ of a Newtonian fluid that conducts heat assuming that Fourier's law is used in the mathematical model representing the fluid. Therefore, in this model the fluid simultaneously undergoes two irreversible processes that are responsible for entropy production in the system: viscosity and heat conduction. The NSF system can therefore be considered as an Euler-Poincaré system with thermodynamics, whose underlying variational principle has been given in \cite{GBYo2017b}. However, the approach developed in \cite{CoGB2018} needs to be extended in order to handle the NSF system. Indeed, in the variational formalism proposed in \cite{GBYo2017b}, the Lie group (of diffeomorphisms) acts on the entropy variable, and more importantly a new variable called \emph{thermal displacement} is introduced\footnote{See \cite{PoGu2009} for an historical account.}, which couples the variational condition to the variational constraint through a specific \emph{thermal coupling term}. Note also that the phenomenological constraint for the NSF system is the local balance of entropy, which is equivalent to the heat equation (see \cite{GBYo2017c}). The paper is organized as follows.

In a first part, after quickly reviewing Euler-Poincaré reduction, we recall the variational principles for the NSF system as developed in \cite{GBYo2017b}. Contrary to the original article, we specialize to the case of 2D or 3D euclidean space (typically meaning that the domain $M$ is considered to be an open and bounded set of $\R^2$ or $\R^3$ with smooth boundary $\partial M$), instead of a generic Riemannian manifold, to keep things simple and in line with the application that we will deal with. 

In the second part of this article, we proceed with the spatial discretization of the NSF system. Firstly, we recall and extend the discrete exterior calculus that is used in the spatial discretization of the partial differential equations arising from variational principles associated to Euler-Poincaré systems on a group of diffeomorphisms. At the heart of this approach is the notion of \emph{discrete diffeomorphism} (see \cite{PaMuToKaMaDe2011} for the incompressible case and \cite{BaGB2018} for the compressible case). This allows us to discretize the infinite dimensional Lie group of diffeomorphisms $\Diff_0(M)$ into a matrix Lie group $\mathsf{D}(\mathbb{M})$, whose dimension is directly related to the number of cells of the mesh $\mathbb{M}$ discretizing our domain $M$. The corresponding Lie algebra is the Lie algebra $\LieAlgebraFont{d}(\mathbb{M})$ of \emph{discrete vector fields}. These discrete vector fields are matrices which satisfy a certain relation, akin to a \emph{discrete divergence theorem}. Moreover, there are two important subsets (which \emph{are not} Lie subalgebras of $\LieAlgebraFont{d}(\mathbb{M})$), one being the vector subspace $\mathcal{S}$ of discrete vector fields that satisfy a certain sparsity constraint based on cell adjacency, whereas the other one, denoted by $\mathcal{V}$, is the vector subspace of \emph{physical} discrete velocity fields. The space of discrete momenta is then shown to be $\mathcal{V}^*\cong\Omega_d^1(\mathbb{M})$, where $\Omega_d^1(\mathbb{M})$ is the vector space of \emph{discrete differential 1-forms}. Then, a discrete exterior calculus is developed based on these objects. Of particular importance are the \emph{discrete Lie derivatives} of discrete differential 1-forms and discrete density-valued differential 1-forms. An explicit formula is given for simplicial grids. Secondly the original phenomenological constraint is recast in such a way that we will not have to deal with the discretization of a general covariant derivative, which is a difficult subject of its own, see \cite{DiscreteConnections1} and \cite{DiscreteConnections2} for instance. Thirdly we describe how to obtain a semi-discrete nonholonomic variational principle, that will need further treatment in the next part to obtain a fully discrete variational principle. The associated Euler-Lagrange equations yield semi-discrete evolution equations for the NSF system. We note already at this point that the nonholonomic character of the obtained variational principle is threefold: it comes from the presence of the nonholonomic phenomenological constraint, but also from the fact that we require discrete vector fields modeling the velocity to be in both $\mathcal{S}$ and $\mathcal{V}$, which are not Lie subalgebras of $\LieAlgebraFont{d}(\mathbb{M})$ as stated before.

\begin{sidewaysfigure}	
	\centering
	\begin{tikzpicture}
		\tikzset{
			big box/.style = {
				shape = rectangle,
				align = center,
				draw,
				inner sep = 0.2cm},
			small box/.style = {
				shape = rectangle,
				align = center,
				draw,
				fill = #1,
				inner sep = 0.2cm},
			comment/.style = {
				shape = rectangle,
				align = center,
				draw,
				rounded corners}
		}
	
		\matrix[column sep=1cm, row sep=1cm]{
			
			\node[big box] (smooth-variational) 
			{
				\textbf{Smooth variational formalism}\\[0.2cm] 
				Infinite-dimensional Euler-Poincaré system on $\Diff_0(M)$ with thermodynamics\\
				Material $\Diff_0(M)$ (section \ref{ssec:smooth-material}) $\xrightarrow{\text{reduction}}$ Spatial $\LieAlgebraFont{X}_0(M)$ (section \ref{ssec:smooth-spatial})
			};
		
			& &
			
			\node[small box = black!5] (smooth-NSF) {Classical Navier-Stokes-Fourier\\ system of PDE \eqref{eq:smooth-classical-Navier-Stokes-Fourier}}; \\
			
			\node[comment] (space-discretization) {Spatial discretization using Pavlov's DEC\\ (section \ref{sec:DEC} and summary \ref{tab:smooth-discrete-summary})};
			
			& 
			
			\node[comment] (geometric-mechanics) {Geometric mechanics\\ (section \ref{ssec:smooth-Euler-Poincaré})};
			
			& \\
			
			\node[big box] (semi-discrete-variational) 
			{
				\textbf{Semi-discrete variational formalism}\\[0.2cm] 
				Finite-dimensional Euler-Poincaré system on $\mathsf{D}(\mathbb{M})$ with thermodynamics\\
				Material $\mathsf{D}(\mathbb{M})$ $\xrightarrow{\text{reduction}}$ Spatial $\LieAlgebraFont{d}(\mathbb{M})$ (section \ref{ssec:semi-discrete-spatial})
			};
		
			& &
			
			\node[small box = white] (semi-discrete-NSF) {Semi-discrete Navier-Stokes-Fourier\\ system of PDE \eqref{eq:semi-discrete-classical-Navier-Stokes-Fourier}}; \\
			
			\node[comment] (time-discretization) {Time discretization according to Veselov\\ (section \ref{ssec:discrete-Euler-Poincaré})};
			
			&  & \\
			
			\node[big box] (discrete-variational) 
			{
				\textbf{Fully discrete variational formalism}\\[0.2cm]
				Discrete Euler-Poincaré system with thermodynamics\\
				Material $\mathsf{D}(\mathbb{M})$ $\xrightarrow{\text{reduction}}$ Spatial $\LieAlgebraFont{d}(\mathbb{M})$ (section \ref{ssec:discrete-spatial})
			};
		
			& & 
			
			\node[small box = white] (discrete-NSF-1) {Variational integrator for the\\ Navier-Stokes-Fourier system on $\mathsf{D}(\mathbb{M})$}; \\
			
			& &
			
			\node[comment] (Cayley-map) {Cayley map\\ (section \ref{ssec:discrete-Euler-Poincaré})}; \\
			
			& &
			
			\node[small box = black!5] (discrete-NSF-2) {Variational integrator for the\\ Navier-Stokes-Fourier system on $\LieAlgebraFont{d}(\mathbb{M})$\\ (section \ref{ssec:discrete-spatial})}; \\
		};
	
		\draw[thick] (smooth-variational) edge (space-discretization);
		\draw[thick,->,>=stealth] (space-discretization) edge (semi-discrete-variational);
		\draw[thick] (semi-discrete-variational) edge (time-discretization);
		\draw[thick,->,>=stealth] (time-discretization) edge (discrete-variational);
		
		\draw[thick,->,>=stealth] (smooth-variational) -- (smooth-NSF) node [above, pos=0.5] {Euler-Lagrange equations};
		\draw[thick,->,>=stealth] (semi-discrete-variational) -- (semi-discrete-NSF) node [above, pos=0.5] {Euler-Lagrange equations};
		\draw[thick,->,>=stealth] (discrete-variational) -- (discrete-NSF-1) node [above, pos=0.5] {Discrete Euler-Lagrange equations};
		
		\draw[thick] (discrete-NSF-1) edge (Cayley-map);
		\draw[thick,->,>=stealth] (Cayley-map) edge (discrete-NSF-2);
		
		\draw[very thick,dashed,->,>=stealth] ([xshift=1cm]smooth-NSF.east) |- (discrete-NSF-2.east);
		\draw[very thick,dashed] ([xshift=1cm]smooth-NSF.east) -- (smooth-NSF.east);
		
		\draw[thick] (smooth-NSF.south) |- (geometric-mechanics.east);
		\draw[thick,->,>=stealth] (geometric-mechanics.west) -| ([xshift=5cm]smooth-variational.south);
	\end{tikzpicture}
	\caption{Structure of the article.}
	\label{fig:article}
\end{sidewaysfigure}
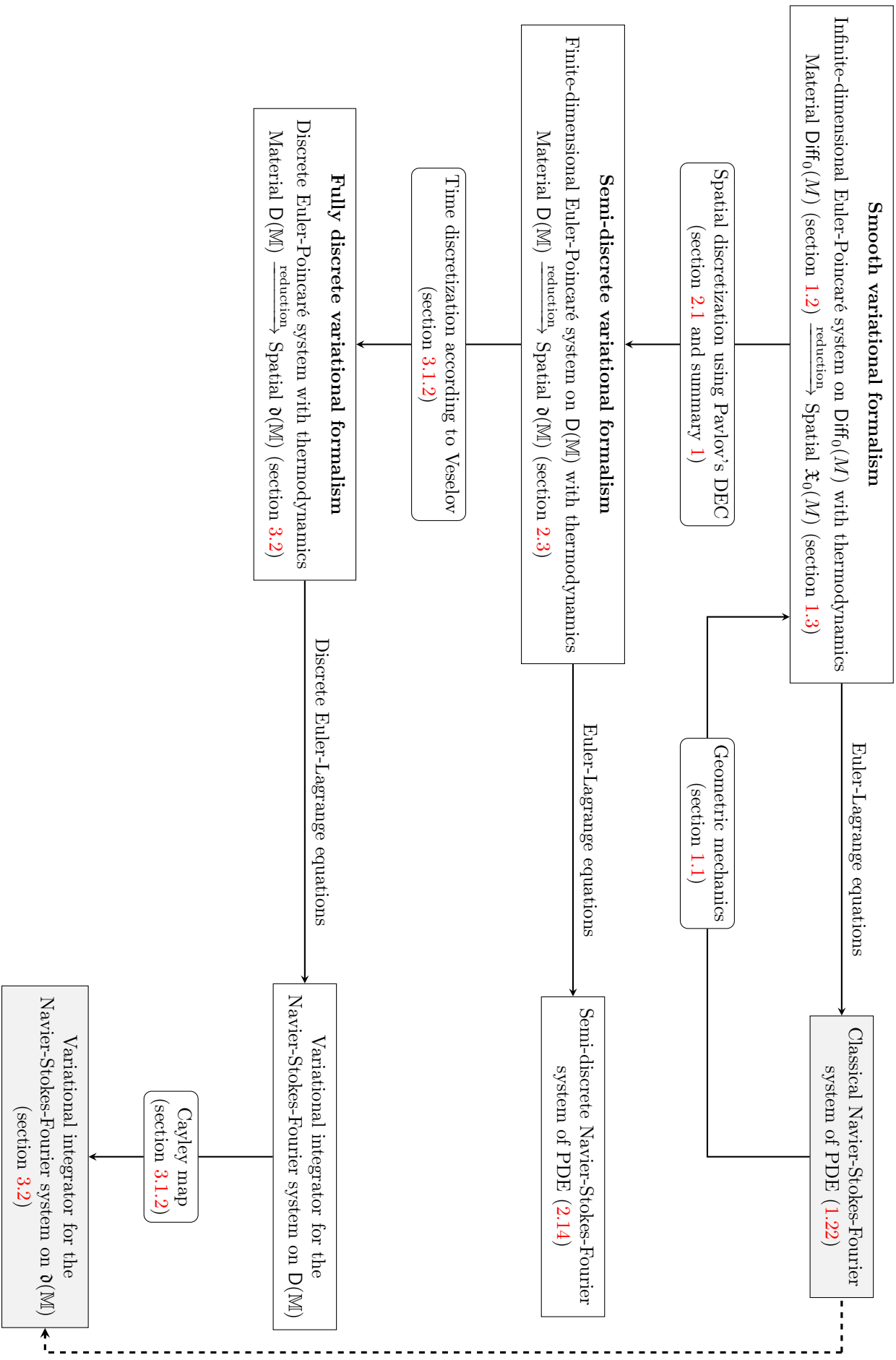

In the last part, after reviewing the variational discretization of Euler-Poincaré systems and group difference maps, we deal with the time discretization of the semi-discrete variational principle obtained in the previous part. This is done \emph{à la Veselov}, meaning that we do not discretize the semi-discrete evolution equations but rather the semi-discrete variational principle directly, in a way that extends \cite{CoGB2018}. Thus we obtain a variational integrator for the NSF system, in the sense of \cite{MaWe2001} and \cite{GBYo2018}. Given the discrete velocity field, the discrete density and the discrete entropy of the fluid at some discrete time $t_n$, the variational integrator yields these variables at the next discrete time $t_{n+1}$, allowing us to understand the evolution of the fluid state and other relevant quantities. The structure of the article is summarized in figure \ref{fig:article}.

%---------------------------------------------------------------------------------------------------

\section{The variational framework for the Navier-Stokes-Fourier system}

%---------------------------------------------------------------------------------------------------

\subsection{Review of Euler-Poincaré systems}\label{ssec:smooth-Euler-Poincaré}

In this section we consider mechanical systems whose configuration space is a finite-dimensional Lie group $G$. Given $g\in G$, denote by $R_g$ the right multiplication by $g$ in the group $G$, $\LieAlgebraFont{g}$ its Lie algebra, and $\omega_G\in\Omega^1(G,\LieAlgebraFont{g})$ its right Maurer-Cartan form; recall that $\omega_G$ is the $\LieAlgebraFont{g}  $-valued one-form on $G$ defined by $\omega_G(\dot{g})=T_g R_{g^{-1}}(\dot{g})\in T_eG\cong\LieAlgebraFont{g}$ for any $\dot{g}\in T_gG$, see \citet[chapter 3, Definition 1.3]{Sh1997}. The right action of $G$ on itself can be lifted to $TG$, and we get $TG/G\cong\LieAlgebraFont{g}$, where the diffeomorphism is given by the Maurer-Cartan form $\omega_G$. 

Given a $G$-invariant Lagrangian $\Lag:TG\to\R$ on $G$ (relatively to the lifted action), and the associated Euler-Lagrange equations, it is natural to ask how one can obtain equivalent equations directly on the Lie algebra $\LieAlgebraFont{g}$, which is the realization of the \emph{reduced} velocity phase space $TG/G$. This process is called \emph{Euler-Poincar\'e reduction} and the equations obtained in this way are called the \emph{Euler-Poincar\'e equations} on $ \LieAlgebraFont{g}$, see \citet[section 13.5]{MaRa1999} for details as well as an historical overview. Euler-Poincar\'e reduction is a particular instance of Lagrangian reduction, see \cite{MaSc1993a,MaSc1993b} and \cite{CeMaRa2001} for instance, in which one considers a Lagrangian $L:TQ \rightarrow \mathbb{R}  $ invariant under the tangent lifted action of a free and proper group action of a Lie group on the configuration space $Q$. 

Euler-Poincar\'e equations can also accommodate \emph{advected parameters}, which are very useful in applications. These parameters, initially fixed, acquire dynamics after reduction in the form of an \emph{advection equation}. In this case the Lagrangian is only invariant under the isotropy subgroup of a given reference parameter. The case of advected parameters taking values in (the dual of) a vector space has been studied in \citet[section 3]{HoMaRa1998} and the general case of advected parameters taking values in manifolds has been developed in \cite{GBTr2010}. We shall follow here this more general setting.

Denoting by $P$ the manifold in which the parameters live we consider a right action of $G$ on $P$, simply denoted by concatenation as $(g,a) \mapsto ag$. The infinitesimal generator associated to $ \xi\in\LieAlgebraFont{g}$ is the vector field on $P$ denoted by $\xi _P$. Given a reference parameter $a_\text{ref}\in M$, we use the notation $G_{a_{\rm ref}} \subset G$ and $\Orb(a_\text{ref}) \subset P$ for the isotropy subgroup and the orbit of $a_{\rm ref}$, respectively. Assuming $G_{a_{\rm ref}}$-invariance, the reduced velocity phase space is $TG/G_{a_\text{ref}}\cong\LieAlgebraFont{g}\times\Orb(a_\text{ref})\subset\LieAlgebraFont{g}\times P$, where the isomorphism is the map $[g,\dot{g}]\mapsto(\dot{g}g^{-1},a_\text{ref}g^{-1})$. In what follows, unless necessary, all actions will be typed with concatenations for the sake of simplicity. Let $\Lag_{a_\text{ref}}:TG\to\R$ be a $G_{a_\text{ref}}$-invariant Lagrangian, this means that 
\[
	\Lag_{a_{\rm ref}}(hg,\dot{h}g)=\Lag_{a_{\rm ref}}(h,\dot{h}),
\]
for all $g \in G_{a_{\rm ref}}$ and $(h,\dot{h})\in TG$. Then we can define a reduced Lagrangian $\lag:\LieAlgebraFont{g}\times\Orb(a_\text{ref})\to\R$ by setting
\[
	\lag(\dot{g}g^{-1},a_{\rm ref}g^{-1})=\Lag_{a_{\rm ref}}(g,\dot{g}),
\]
for all $(g,\dot{g})\in TG$. A similar construction can be used for any $G_{a_\text{ref}}$-invariant map.

The proof of the Euler-Poincar\'e reduction theorem is based on the following technical result proven in \citet[proposition 5.1]{BlKrMaRa1994}. We give below a simpler proof of this well-known result.

\begin{lemma}\label{lem:Euler-Poincaré-lemma}
	Let $G$ be a Lie group, $\LieAlgebraFont{g}$ its Lie algebra, $U$ an open set of $\R^2$, $g:U\to G$ a smooth map, and define two $\LieAlgebraFont{g}$-valued maps given for any $(t,\varepsilon)\in U$ by:
	\begin{equation}\label{eq:Euler-Poincaré-xi-eta}
	\xi(t,\varepsilon)=\frac{\partial g}{\partial t}(t,\varepsilon)g(t,\varepsilon)^{-1},\quad
	\eta(t,\varepsilon)=\frac{\partial g}{\partial\varepsilon}(t,\varepsilon)g(t,\varepsilon)^{-1}.
	\end{equation}
	Then on the open set $U$ we have:
	\begin{equation}\label{eq:EP-constraint}
	\frac{\partial\xi}{\partial\varepsilon}-\frac{\partial\eta}{\partial t}=-[\xi,\eta].
	\end{equation}
\end{lemma}

\begin{proof}
	We first handle the case of matrix Lie groups. In the case where $G$ is a matrix Lie group, the derivative of the right multiplication is the right multiplication itself, owing to the fact that $G$ is a subset of a vector space. Therefore $\xi$ and $\eta$ are now matrix-valued. Differentiating, using the product rule and equality of mixed partial derivatives, we obtain that
	\[
		\frac{\partial\xi}{\partial\varepsilon}-\frac{\partial\eta}{\partial t}=
		\frac{\partial^2 g}{\partial\varepsilon\partial t}g^{-1}
		-\frac{\partial g}{\partial t}g^{-1}\frac{\partial g}{\partial\varepsilon}g^{-1}
		-\frac{\partial^2 g}{\partial t\partial\varepsilon}g^{-1}
		+\frac{\partial g}{\partial\varepsilon}g^{-1}\frac{\partial g}{\partial t}g^{-1}
		=-[\xi,\eta].
	\]
	We now turn to the case of general Lie groups. Given a smooth map $g:U\to G$, denoting by $\omega_g=g^*\omega_G$ the right Darboux derivative\footnote{Also known as the right \emph{logarithmic derivative} of $g$.} of $g$, we have
	\[
		\dif\omega_g-\frac{1}{2}[\omega_g,\omega_g]=\dif(g^*\omega_G)-\frac{1}{2}[g^*\omega_G,g^*\omega_G]=g^*\left(\dif\omega_G-\frac{1}{2}[\omega_G,\omega_G]\right)=0,
	\]
	because the right Maurer-Cartan form $\omega_G$ is known to satisfy the \emph{structural equation} $\dif\omega_G-\frac{1}{2}[\omega_G,\omega_G]=0$. Now, evaluating $\dif\omega_g-\frac{1}{2}[\omega_g,\omega_g]$ on $\frac{\partial}{\partial t}$ and $\frac{\partial}{\partial\varepsilon}$ yields \eqref{eq:Euler-Poincaré-xi-eta}. Indeed, since $\left[\frac{\partial}{\partial t},\frac{\partial}{\partial\varepsilon}\right]=0$, the intrinsic formula for the exterior derivative allows us to compute that
	\begin{align*}
		\frac{\partial}{\partial t}\left[g^*\omega_G\left(\frac{\partial}{\partial\varepsilon}\right)\right]&-\frac{\partial}{\partial\varepsilon}\left[g^*\omega_G\left(\frac{\partial}{\partial t}\right)\right]-\frac{1}{2}\left[g^*\omega_G\left(\frac{\partial}{\partial t}\right),g^*\omega_G\left(\frac{\partial}{\partial\varepsilon}\right)\right]\\
		&\quad+\frac{1}{2}\left[g^*\omega_G\left(\frac{\partial}{\partial \varepsilon}\right),g^*\omega_G\left(\frac{\partial}{\partial t}\right)\right]=\frac{\partial\eta}{\partial t}-\frac{\partial\xi}{\partial\varepsilon}-\frac{1}{2}[\xi,\eta]+\frac{1}{2}[\eta,\xi]=0.
	\end{align*}
\end{proof}

We can now recall the Euler-Poincaré reduction theorem with advected parameters.

\begin{theorem}[\textbf{Euler-Poincar\'e reduction with advected parameters}]\label{thm:Euler-Poincaré-reduction}
	Let $G$ be a Lie group and $\LieAlgebraFont{g}$ its Lie algebra, and let $P$ be a manifold on which $G$ acts on the right. For a fixed parameter $a_\text{ref}\in P$ let $\Lag_{a_\text{ref}}:TG\times\R\to\R$ be a $G_{a_\text{ref}}$-invariant Lagrangian and denote the corresponding reduced Lagrangian by $\lag:\LieAlgebraFont{g}\times\Orb(a_\text{ref})\to\R$. Then the following assertions are equivalent:
	\begin{enumerate}[label=(\arabic*), ref=\thetheorem.(\arabic*)]
		
		\item The curve $g(t)\in G$ is critical for the variational principle
		\begin{equation}
			\delta\int_0^T\Lag_{a_\text{ref}}(g,\dot{g})\,\mathrm{d}t=0
			\label{eq:Euler-Poincaré-before-reduction}
		\end{equation}
		where variations $\delta g$ of $g$ vanish at $t=0, T$.
		
		\item The curve $g(t)\in G$ satisfies the Euler-Lagrange equations obtained from the variational principle \eqref{eq:Euler-Poincaré-before-reduction}.
		
		\item The curve $(\xi(t),a(t))\in\LieAlgebraFont{g} \times\operatorname{Orb}(a_{\rm ref})$, defined by $\xi(t)=\dot g(t)g(t)^{-1}$, $a(t)=a_{\rm ref}g(t)^{-1}$, is critical for the reduced variational principle
		\[
			\delta\int_0^T\lag(\xi,a)\,\mathrm{d}t=0
		\]
		subject to the variational constraints
		\[
			\delta\xi=\dot{\eta}-[\xi,\eta],\quad\delta a=-\eta_P(a),
		\]
		where $\eta$ is any curve in $\LieAlgebraFont{g}$ vanishing at $t=0,T$.	
			
		\item The curve $(\xi(t),a(t))\in\LieAlgebraFont{g}\times\operatorname{Orb}(a_{\rm ref})$ satisfies the equations
		\begin{empheq}[left=\empheqlbrace]{align}
			& \frac{\mathrm{d}}{\mathrm{d}t}\frac{\delta\lag}{\delta\xi}(\xi,a)=-\ad^*_\xi\frac{\delta\lag}{\delta\xi}(\xi,a)-\mathbf{J}\left(\frac{\delta\lag}{\delta a}(\xi,a)\right),\label{eq:Euler-Poincaré-equation}\\[2mm]
			& \dot{a}+\xi_P(a)=0\label{eq:advection-equation},
		\end{empheq}
		where $\mathbf{J}:T^*P\to\LieAlgebraFont{g}^*$ is the momentum map associated to the cotangent lift of the action of $G$ on $T^*P$; it is defined by $\big\langle\mathbf{J}(\alpha_x),\xi\big\rangle=\big\langle\alpha_x,\xi_P(x)\big\rangle$ for any $x\in P$, $\alpha_x\in T_x^*P$ and $\xi\in\LieAlgebraFont{g}$.
		
	\end{enumerate}
\end{theorem}

In these equations, \eqref{eq:advection-equation} is called the \emph{advection equation}. Equation \eqref{eq:Euler-Poincaré-equation} is called the \emph{Euler-Poincar\'e equation}. Also, the relation $\dot{g}=g\xi$ is called the \emph{reconstruction equation}. In \eqref{eq:Euler-Poincaré-equation} the variational derivatives are computed relatively to the natural duality pairing between $\LieAlgebraFont{g}^*$ and $\LieAlgebraFont{g}$. More on Euler-Poincaré reduction can be found in \cite{HoMaRa1998} or \cite{MaRa1999}. Of course, the group can also act on the left of itself rather than on the right, but in the case of fluids the group acts on the right as we will see below; $G$ can also act on the left of $P$ too.

\begin{proof} 
	The equivalence between (1) and (2) follows from a direct computation, see \citet[chapter 7]{MaRa1999}. We now show that (3) and (4) are equivalent. Taking variations from the left hand side of the variational condition in (3), we obtain using all the available constraints:
	\begin{align*}
		\int_0^T\left\langle\frac{\delta\lag}{\delta\xi},\delta\xi\right\rangle\,\mathrm{d}t
		+\int_0^T\left\langle\frac{\delta\lag}{\delta a},\delta a\right\rangle\,\mathrm{d}t
		&=\int_0^T\left\langle\frac{\delta\lag}{\delta\xi},\dot{\eta}-[\xi,\eta]\right\rangle\,\mathrm{d}t
		-\int_0^T\left\langle\frac{\delta\lag}{\delta a},\eta_P(a)\right\rangle\,\mathrm{d}t\\
		&=\int_0^T\left\langle-\frac{\mathrm{d}}{\mathrm{d}t}\frac{\delta\lag}{\delta\xi}-\ad^*_\xi\frac{\delta\lag}{\delta\xi}-\mathbf{J}\left(\frac{\delta\lag}{\delta a}\right),\eta\right\rangle\,\mathrm{d}t+\left[\left\langle\frac{\delta\lag}{\delta\xi},\eta\right\rangle\right]_0^T,
	\end{align*}
	which yields equations \eqref{eq:Euler-Poincaré-equation}, the last term in the right hand side being zero. The advection equation comes from a slightly more technical computation. Denoting by $\sigma:G\times P\to P$ the right action of $G$ on $P$ and by $ \mathbf{d}\sigma_{(g,a)}:T_g G\times T_a\operatorname{Orb}(a_{\rm ref})\rightarrow \mathbb{R}$ its derivative at $(g,a)$, we compute that
	\begin{equation*}
		\dot{a}=\dif\sigma_{(g^{-1},a_\text{ref})}(-g^{-1}\dot{g}g^{-1},0).
	\end{equation*}
	Introducing for $g\in G$ the map $\Psi_g:G\times P\to G\times P$ defined by $\Psi_g(h,a)=(g^{-1}h,ag)$ we obtain
	\[
		\dot{a}=\Psi_{g}^*\dif\sigma_{(e,a)}(-\xi,0).
	\]
	Since the pullback commutes with the differential and $\Psi_g^*\sigma=\sigma$ for any $g\in G$, we get the advection equation $\dot{a}=\mathbf{d}\sigma_{(e,a)}(-\xi,0)=-\xi_P(a)$.
	
	We are now concerned with the equivalence of (1) and (3). Suppose that (1) holds first. Setting $\eta=\delta g g^{-1}$, which is a $\LieAlgebraFont{g}$-valued curve, we obtain the desired variational formulation from left $G_{a_\text{ref}}$-invariance of the Lagrangian. It remains to understand how variations $\delta g$ of the curve $g$ in $G$ induce variations $\delta\xi$ of the curve $\xi$ in $\LieAlgebraFont{g}$ and variations $\delta a$ of the curve $a$ in $P$, and conversely. First recall how variations $\delta g$ of $g$ are constructed: consider a map $(t,\varepsilon)\in U\mapsto\bar{g}(t,\varepsilon)\in G$ such that $\bar{g}(t,0)=g(t)$ for all $t\in[0,T]$ (a curve of curves in $G$ with the same endpoints as $g$ for all the curves) and where $U$ is some open set in $\R^2$. Then $\delta g(t)=\frac{\partial\bar{g}}{\partial\varepsilon}(t,0)\in T_{g(t)}G$ for all $t\in[0,T]$ (see \citet[proposition 8.1.2]{MaRa1999} for more details). Now consider the maps $\bar{\xi}:U\to\LieAlgebraFont{g}$ and $\bar{\eta}:U\to\LieAlgebraFont{g}$ defined by
	\[
		\bar{\xi}(t,\varepsilon)=\frac{\partial\bar{g}}{\partial t}(t,\varepsilon)\bar{g}(t,\varepsilon)^{-1},\quad
		\bar{\eta}(t,\varepsilon)=\frac{\partial\bar{g}}{\partial\varepsilon}(t,\varepsilon)\bar{g}(t,\varepsilon)^{-1}.
	\]
	Notice that for all $t\in[0,T]$, $\xi(t)=\bar{\xi}(t,0)$ and $\bar{\eta}(t,0)=\eta(t)$. Now thanks to Lemma \ref{lem:Euler-Poincaré-lemma} we have for all $(t,\varepsilon)\in U$ that
	\[
		\frac{\partial\bar{\xi}}{\partial\varepsilon}(t,\varepsilon)-\frac{\partial\bar{\eta}}{\partial t}(t,\varepsilon)=-\big[\bar{\xi}(t,\varepsilon),\bar{\eta}(t,\varepsilon)\big],
	\]
	and taking $\varepsilon=0$ in this equality yields $\delta\xi=\dot{\eta}-[\xi,\eta]$. This describes entirely the variations $\delta\xi$ of $\xi$ in terms of $\delta g$ since $\eta$ directly depends on $\delta g$. Concerning variations $\delta a$ of the curve $a$ we have $\delta a=-\eta_P(a)$ using a computation similar to the one used to obtain the advection equation above.
	
	Conversely, given a curve $\xi$ in $\LieAlgebraFont{g}$ and a variation $\delta\xi$ of $\xi$, we define $g$ to be the curve in $G$ obtained from solving the reconstruction equation $\dot{g}=\xi g$. Then we set $\delta g=\eta g$. Since $\eta$ is arbitrary and zero at endpoints, $\delta g$ is arbitrary and zero at endpoints. We finish the proof using the $G_{a_\text{ref}}$-invariance to obtain the variational principle (1) from the variational principle (3).
\end{proof}

%---------------------------------------------------------------------------------------------------

\subsubsection{Compressible Euler equations}\label{ssec:compressible-Euler}

As it is fundamental regarding this article, we now detail the case of compressible Euler equations of inviscid fluid dynamics. A reference on this topic is \citet[section 6]{HoMaRa1998}. The configuration space is the infinite-dimensional Lie group $G=\Diff(M)$ of diffeomorphisms of $M$, where $M$ is an open bounded set of $\R^3$ (or $\R^2$). Its Lie algebra is the Lie algebra of vector fields $\LieAlgebraFont{X}(M)$ with Lie bracket given by minus the Jacobi-Lie bracket of vector fields $-[\mathbf{u},\mathbf{v}]=\ad_\mathbf{u}\mathbf{v}$, for all $\mathbf{u}$, $\mathbf{v}\in\LieAlgebraFont{X}(M)$, and such that $\mathbf{u}\cdot\mathbf{n}=0$ on $\partial M$, the so-called \emph{tangential} boundary condition, where $\mathbf{n}$ denotes the outward normal vector field. Recall that densities on $M$, denoted by $\Den(M)$, can be considered as the dual of smooth functions $\smooth{M}$, relatively to the pairing
\begin{equation}\label{eq:density-function-inner-product}
	\langle f,\mu\rangle=\int_M f\mu,
\end{equation}
for any $f\in\smooth{M}$ and $\mu\in\Den(M)$. Then, the dual of $\LieAlgebraFont{X}(M)$ can be identified with density-valued 1-forms $\Omega^1(M)\otimes\Den(M)$ with the help of the pairing
\begin{equation}\label{eq:1-form-densities-vector-fields-pairing}
	\langle\alpha\otimes\mu,\mathbf{u}\rangle=\int_M\langle\alpha,\mathbf{u}\rangle\mu,
\end{equation}
defined for all $\alpha\in\Omega^1(M)$, $\mu\in\Den(M)$ and $\mathbf{u}\in\LieAlgebraFont{X}(M)$. With respect to this pairing the adjoint of $\ad_{\mathbf{u}}$ is given by
\[
	\ad_\mathbf{u}^*(\alpha\otimes\mu)=\big[\lie_\mathbf{u}\alpha+(\dive_\mu\mathbf{u})\alpha\big]\otimes\mu=\lie_\mathbf{u}(\alpha\otimes\mu).
\]
In particular for $\mu=\rho\,\mathrm{d}x$ where $\mathrm{d}x$ is the euclidean volume form\footnote{We will often simply write $\rho$ for the density $\rho\,\mathrm{d}x$.}, notice that:
\[
	\big(\!\dive_{\rho\,\mathrm{d}x}\mathbf{u}\big)\rho\,\mathrm{d}x
	=\big[\rho\dive\mathbf{u}+\grad\rho\cdot\mathbf{u}\big]\mathrm{d}x
	=\dive(\rho\mathbf{u})\,\mathrm{d}x,
\]
where $\dive$ denotes the divergence with respect to $\mathrm{d}x$ and where we used the flat operator associated with the euclidean metric.

Although $\Diff(M)$ is infinite-dimensional, the framework of Euler-Poincaré reduction still makes sense (see \citet[section 6]{HoMaRa1998} for details). The action of $\Diff(M)$ on $M$ is known as the \emph{particle relabeling symmetry} in this case. Before giving explicitly the Lagrangian, we need to introduce a few standard notions. $M$ is thought as the domain where the \emph{material} particles live (their coordinates will be denoted by $X$). As the fluid undergoes a \emph{motion} which consists of a curve of deformations $\varphi_t:=\varphi(t)\in\Diff(M)$, a particle which was located at $X$ has moved to $x(X,t)=\varphi_t(X)$ after $t$ units of time. Therefore, $M$ can also be thought as the domain where \emph{spatial} particles live (their coordinate will be denoted by $x$ in this case). In the \emph{Lagrangian formalism} we are interested in the motion of a particle located at $X\in M$ whereas in the \emph{spatial formalism} we are interested in the particles that go along the flow across a location $x\in M$; the spatial formalism is also called the \emph{Euler specification} of the fluid flow and is the one which is prevalent in casual fluid mechanics. The \emph{Lagrangian velocity} is $\mathbf{U}(X,t)=\frac{\mathrm{d}\varphi_t}{\mathrm{d}t}(X)=\dot{\varphi}_t(X)\in T_{\varphi_t(X)}M$ for all $X\in M$; by definition $\dot{\varphi}_t\in T_{\varphi_t}\Diff(M)$ at any time $t$. The \emph{Eulerian velocity} is $\mathbf{u}(x,t)=\dot{\varphi}_t\big(\varphi_t^{-1}(x)\big)$ for any $x\in M$ and time $t$, which means that the Eulerian velocity actually looks at the velocity of the particle where it was $t$ units of time before. 

Let $P=\Den(M)$ be our parameter manifold and set $\rho_\text{ref}\in P$ to be the reference mass density of the fluid that will be advected along the flow. The lifted action of $\Diff(M)$ on $T\Diff(M)$ is $\dot{\varphi}\psi=\dot{\varphi}\circ\psi$ and the action of $\Diff(M)$ on $\Den(M)$ is given by $\rho\varphi=(\rho\circ\varphi)J_{\varphi}$, where $J_{\varphi}(X)=\det{T_X\varphi}$ is the Jacobian of $\varphi$. The Lagrangian $\Lag_{\rho_\text{ref}}:T\Diff(M)\to\R$ is given by:
\[
	\Lag_{\rho_\text{ref}}(\varphi,\dot{\varphi})
	=\int_M\Lagd_{\rho_\text{ref}}\big(\varphi(X),\dot{\varphi}(X),T_X\varphi\big)\,\mathrm{d}X
	=\int_M\left[\frac{1}{2}\rho_{\text{ref}}(X)\big\|\dot{\varphi}(X)\big\|^2_{\varphi(X)}-\varepsilon\left(\frac{\rho_{\text{ref}}(X)}{J_\varphi(X)}\right)J_\varphi(X)\right]\mathrm{d}X,
\]
where $\varepsilon$ denotes the internal energy of the fluid (see \citet[section 2.3]{GBYo2017b} for more details). $\Lagd_{\rho_\text{ref}}$ is called a \emph{Lagrangian density} associated to $\Lag_{\rho_\text{ref}}$. Although it is seldom used in fluid mechanics, let's derive the Euler-Lagrange equations associated to this variational formalism as it will play an important role in what follows\footnote{We would like to warn the reader that the material variational principle and the associated Euler-Lagrange equations require more advanced differential geometry that the one used for the forthcoming spatial variational principle. We will not dive into the necessary theory as it would require a quite long exposition but the reader is referred to \citet[chapter 5, section 4]{MaHu1983} for an overview of the geometry of classical field theories.}.  We begin by computing the variations of the action functional based on the Lagrangian $\Lag_{\rho_\text{ref}}$ above, where $t\mapsto\varphi(t)$ denotes a curve in $\Diff(M)$:
\begin{equation}\label{variations_L}
	\delta\int_0^T\Lag_{\rho_\text{ref}}(\varphi,\dot{\varphi})\,\mathrm{d}t=\int_0^T\int_M\left(\frac{\partial\Lagd_{\rho_\text{ref}}}{\partial\varphi}\cdot\delta\varphi+\frac{\partial\Lagd_{\rho_\text{ref}}}{\partial\dot{\varphi}}\cdot\delta\dot{\varphi}+\frac{\partial\Lagd_{\rho_\text{ref}}}{\partial T_X\varphi}:\delta T_X\varphi\right)\mathrm{d}X\,\mathrm{d}t.
\end{equation}
Then integrating by parts in time and using the generalized divergence theorem for tensors (see \cite[Problem 7.6]{MaHu1983}) we obtain:
\[
	\delta\int_0^T\Lag_{\rho_\text{ref}}(\varphi,\dot{\varphi})\,\mathrm{d}t=\int_0^T\int_M\left(-\frac{\mathrm{D}}{\mathrm{D}t}\frac{\partial\Lagd_{\rho_\text{ref}}}{\partial\dot{\varphi}}-\operatorname{DIV}\frac{\partial\Lagd_{\rho_\text{ref}}}{\partial T_X\varphi}\right)\cdot\delta\varphi\,\mathrm{d}X\,\mathrm{d}t.
\]
In this relation, $\operatorname{DIV}$ is the usual divergence operator (albeit for two-point tensors), but the uppercase is used to emphasize the fact that we use material coordinates (see \citet[chapter 1]{MaHu1983}). Therefore the \emph{material equations of motion} are given by:
\[
	\frac{\mathrm{D}}{\mathrm{D}t}\frac{\partial\Lagd_{\rho_\text{ref}}}{\partial\dot{\varphi}}=\operatorname{DIV}\mathbf{P},\quad\mathbf{P}:=-\frac{\partial\Lagd_{\rho_\text{ref}}}{\partial{T_X\varphi}}.
\]

$\mathbf{P}$ is called the \emph{first Piola-Kirchoff stress tensor}\footnote{This tensor is related to the more familiar Cauchy stress tensor through the Piola transform.}. For the Lagrangian $\Lag_{\rho_\text{ref}}$ given above, we obtain:
\[
	\rho_\text{ref}\frac{\mathrm{D}\mathbf{U}}{\mathrm{D}t}=-(\grad p)\circ(\varphi J_\varphi),
\]
where $p$ is the pressure defined by $p=\varepsilon'(\rho)\rho-\varepsilon$. This evolution equation can be understood as the compressible Euler equation in material coordinates (which is physically the \emph{material balance of momentum}). To this equation we need to add the \emph{continuity equation}, which reads $\rho=\rho_\text{ref}\varphi^{-1}$ in material coordinates, as well as the tangential boundary conditions (recall that $\varphi$ maps $\partial M$ onto $\partial M$ at each time $t$, which is also written as $\mathbf{U}\cdot\mathbf{n}=0$ on $\partial M$), and an appropriate initial condition.

For any $\varphi$ in the isotropic subgroup of $\rho_\text{ref}$, $\Lag_{\rho_\text{ref}}$ is invariant, provided that the internal energy $\varepsilon$ is $\Diff(M)$-equivariant. Then, Euler-Poincaré reduction amounts to going from the material formalism to the spatial formalism: the new variables are the flow velocity $\mathbf{u}=\dot{\varphi}\circ\varphi^{-1}\in\LieAlgebraFont{X}(M)$, the mass density $\rho=\rho_\text{ref}\varphi^{-1}\in\Orb(\rho_\text{ref})$\footnote{Notice that $\rho_\text{ref}$ is not necessarily equal to $\rho_0$ as we have the relation $\rho_0=\rho_\text{ref}\varphi_0^{-1}$.}, and the reduced Lagrangian is:
\[
	\lag(\mathbf{u},\rho)=\Lag_{\rho_\text{ref}}(\Id,\dot{\varphi}\circ\varphi^{-1})=\int_M\left[\frac{1}{2}\rho\langle\mathbf{u}^\flat,\mathbf{u}\rangle-\varepsilon(\rho)\right]\mathrm{d}x,
\]
in which we recognize the kinetic energy of the fluid.

Now, in theorem \ref{thm:Euler-Poincaré-reduction}, going from (1) to (3) is a direct consequence of the relationship between $\Lag_{\rho_\text{ref}}$ and $\lag$ that we just described above; also the variational constraint\footnote{This constraint is also known as the \emph{Lin constraint} in the Physics literature.} $\delta\xi=\dot{\eta}-\ad_\xi\eta$ now reads $\delta\mathbf{u}=\dot{\bm{\zeta}}-\ad_\mathbf{u}\bm{\zeta}=\dot{\bm{\zeta}}+[\mathbf{u},\bm{\zeta}]$, where $\bm{\zeta}\in\LieAlgebraFont{X}(M)$ vanishes at $t=0,T$, and $\delta a=-\eta_P(a)$ reads $\delta\rho=-\lie_{\bm{\zeta}}(\rho\,\mathrm{d}x)=-(\dive_{\rho\,\mathrm{d}x}\bm{\zeta})\rho\,\mathrm{d}x=-\dive(\rho\bm{\zeta})\,\mathrm{d}x$. We use now these relations in going from (3) to (4) which will yield the usual compressible Euler equation. Indeed, we have:
\[
	\delta\int_0^T\lag(\mathbf{u},\rho)\,\mathrm{d}t=\int_0^T\left[\left\langle\frac{\delta\lag}{\delta\mathbf{u}},\delta\mathbf{u}\right\rangle+\left\langle\frac{\delta\lag}{\delta\rho},\delta\rho\right\rangle\right]\mathrm{d}t,
\]
using integration by parts we compute that
\[
	\left\langle\frac{\delta\lag}{\delta\mathbf{u}},\delta\mathbf{u}\right\rangle 
	=\int_M\left\langle\rho\mathbf{u}^\flat,\dot{\bm{\zeta}}\right\rangle\mathrm{d}x+\left\langle\mathbf{u}^\flat\otimes\rho\,\mathrm{d}x,-\ad_\mathbf{u}\bm{\zeta}\right\rangle
	=\left\langle-\frac{\partial}{\partial t}(\rho\mathbf{u}^\flat)\otimes\mathrm{d}x-\ad_\mathbf{u}^*(\mathbf{u}^\flat\otimes\rho\,\mathrm{d}x),\bm{\zeta}\right\rangle,
\]
and using the divergence theorem we also compute that
\begin{align*}
	\left\langle\frac{\delta\lag}{\delta\rho},\delta\rho\right\rangle
	&=-\int_M\frac{1}{2}\langle\mathbf{u}^\flat,\mathbf{u}\rangle\dive(\rho\bm{\zeta})\,\mathrm{d}x+\int_M\frac{\partial\varepsilon}{\partial\rho}\dive(\rho\bm{\zeta})\,\mathrm{d}x\\
	&=\int_M\frac{1}{2}\rho\left\langle\dif\langle\mathbf{u}^\flat,\mathbf{u}\rangle,\bm{\zeta}\right\rangle\mathrm{d}x-\int_M\rho\left\langle\dif\varepsilon'(\rho),\bm{\zeta}\right\rangle\mathrm{d}x\\
	&=\left\langle\frac{1}{2}\rho\,\dif\langle\mathbf{u}^\flat,\mathbf{u}\rangle\otimes\mathrm{d}x-\rho\,\dif\varepsilon'(\rho)\otimes\mathrm{d}x,\bm{\zeta}\right\rangle.
\end{align*}
Therefore, the Euler-Poincaré equation \eqref{eq:Euler-Poincaré-equation} reads:
\[
	\frac{\partial}{\partial t}(\rho\mathbf{u}^\flat)+\rho\lie_\mathbf{u}\mathbf{u}^\flat+\dive(\rho\mathbf{u})\mathbf{u}^\flat-\frac{1}{2}\rho\,\dif\langle\mathbf{u}^\flat,\mathbf{u}\rangle=-\rho\,\dif\varepsilon'(\rho),
\]
and using the advection equation \eqref{eq:advection-equation} $\frac{\partial\rho}{\partial t}=-\dive(\rho\mathbf{u})$, the relation $(\nabla_\mathbf{u}\mathbf{u})^\flat=\lie_\mathbf{u}\mathbf{u}^\flat-\frac{1}{2}\dif\langle\mathbf{u}^\flat,\mathbf{u}\rangle$ (for the Levi-Civita connection associated to the euclidean metric) and the pressure $p=\varepsilon'(\rho)\rho-\varepsilon$, we finally obtain after sharpening:
\[
	\rho\frac{\partial\mathbf{u}}{\partial t}+\rho\nabla_\mathbf{u}\mathbf{u}=-\grad p,
\]
which is the \emph{Euler equation} for compressible flow, to which we need to add the advection equation for $\rho$ (called the \emph{continuity equation} in the Physics literature) and the boundary condition $\mathbf{u}\cdot\mathbf{n}=0$ on $\partial M$. From the point of view of continuum mechanics, the continuity equation corresponds to the local balance of mass and the Euler equation corresponds to the local balance of (linear) momentum. In conclusion, we have seen that the usual compressible Euler equations (the momentum equation, the continuity equation, together with the tangential boundary condition) can be considered as an infinite-dimensional Euler-Poincaré system.

\begin{remark}
	We can compute with the help of the divergence theorem that the momentum map $\mathbf{J}:T^*\Den(M)\to\Omega^1(M)\otimes\Den(M)$ is given by $\mathbf{J}(\rho\,\mathrm{d}x,f)=-\dif f\otimes\rho\,\mathrm{d}x$, for any $(\rho\,\mathrm{d}x,f)\in T^*\Den(M)$.
\end{remark}

\begin{remark}\label{rq:entropy-advection}
	At this point, more thermodynamics can be introduced as follows. Let $S_\text{ref}\in\Den(M)$ denotes the reference entropy density. This reference density, like the reference mass density $\rho_\text{ref}$, is going to acquire dynamics through the Euler-Poincaré reduction process. We define a new variable $s=S_\text{ref}\varphi^{-1}$ for the entropy density and the reduced Lagrangian becomes:
	\[
		\lag(\mathbf{u},\rho,s)=\int_M\left[\frac{1}{2}\rho\langle\mathbf{u}^\flat,\mathbf{u}\rangle-\varepsilon(\rho,s)\right]\mathrm{d}x,
	\]
	with now the internal energy of the fluid also depending on the entropy. Then we obtain the same Euler equation as above though the pressure becomes the \emph{thermodynamic pressure} $p=\frac{\partial\varepsilon}{\partial\rho}\rho+\frac{\partial\varepsilon}{\partial s}s-\varepsilon$, and the entropy density satisfies the advection equation $\frac{\partial s}{\partial t}+\dive(s\mathbf{u})=0$. One of the main achievements of \cite{GBYo2017b} is to consistently incorporate irreversible processes into one nonholonomic variational principle, where the nonholonomic constraint, called the \emph{phenomenological constraint}, will replace the preceding advection equation and reflect how entropy is created into the system. To achieve this goal, $S_\mathrm{ref}$ has to be promoted to a full dynamic variable in the variational principle as we will see in the next two sections.
\end{remark}

%---------------------------------------------------------------------------------------------------

\subsection{The material variational principle for the Navier-Stokes-Fourier system}\label{ssec:smooth-material}

In this section and the next, we recall the variational formalism for nonequilibrium thermodynamics of fluids developed in \cite{GBYo2017b}. The irreversible processes that we are interested in are internal heat conduction and viscosity, which can be thought of as internal generalized friction forces of the fluid. The phenomenological constraint is a nonlinear, nonholonomic constraint that is added to the variational condition; it is exactly formulated as the evolution equation of the entropy density as we will see shortly, and as such contains the sum of the product of thermodynamic affinities by thermodynamic fluxes. For instance, in the material version of the phenomenological constraint, the irreversible process corresponding to viscosity has for thermodynamic flux the deviatoric part $\mathbf{P}^\text{fr}(t,X)$ of the first Piola-Kirchhoff tensor $\mathbf{P}(t,X)$ and for thermodynamic affinity the gradient of the Lagrangian velocity $\nabla\dot{\varphi}$, whereas the irreversible process corresponding to heat conduction has for thermodynamic flux the entropy density flux $\mathbf{J}_S(t,X)$ and for thermodynamic affinity the temperature $\Theta(t,X)$. Therefore, in order to consider heat conduction as a generalized friction force, we need to introduce a new variable $\Gamma(t,X)$ called \emph{thermal displacement} and whose rate equals the temperature of the fluid. Now, technically speaking, a new variable $\Sigma(t,X)$ is also needed so that the associated Euler-Lagrange equations will be exactly the evolution equations for the Navier-Stokes-Fourier system (the physical meaning of $\Sigma$ will be given later); these two new variables appear in a coupling term added to the usual Lagrangian and that we call the \emph{thermal coupling term}.

The configuration space for the Navier-Stokes-Fourier system is the infinite dimensional Lie group $\Diff_0(M)$ of diffeomorphisms of $M$ that fix the boundary $\partial M$ \emph{point-wise}; this will in turn correspond to the no-slip boundary condition in the spatial formalism of the next section. In contrast to what has been done in the remark \ref{rq:entropy-advection}, a key point of the new variational formalism is to promote the entropy density variable to a full dynamic variable $t\mapsto S(t,\cdot)\in\Den(M)$ instead of being considered as a fixed parameter $S_\text{ref}\in\Den(M)$. Therefore, the material Lagrangian becomes:
\begin{align}
	\Lag_{\rho_\text{ref}}(\varphi,\dot{\varphi},S)
	&=\int_M\Lagd_{\rho_\text{ref}}(\varphi,\dot{\varphi},T_X\varphi,S)\,\mathrm{d}X\notag\\
	&=\int_M\left[\frac{1}{2}\rho_{\text{ref}}(X)\big\|\dot{\varphi}(X)\big\|^2_{\varphi(X)}-\varepsilon\left(\frac{\rho_{\text{ref}}(X)}{J_\varphi(X)},\frac{S(X)}{J_\varphi(X)}\right)J_\varphi(X)\right]\mathrm{d}X.\label{eq:material-Lagrangian}
\end{align}
Another key point is to decompose the first Piola-Kirchoff tensor as $\mathbf{P}=\mathbf{P}^\text{cons}+\mathbf{P}^\text{fr}$, where $\mathbf{P}^\text{fr}$ is coupled in the phenomenological constraint to $\nabla\dot{\varphi}$. In geometric terms, $\nabla$ is the Levi-Civita connection on the tangent vector bundle $T\Diff(M)$ associated to the Riemannian metric $\widetilde{g}$ defined for any $\mathbf{U}$, $\mathbf{V}\in T_\varphi\Diff(M)$, $\varphi\in\Diff(M)$, by:
\[
	\widetilde{g}(\mathbf{U},\mathbf{V})_\varphi=\int_M g\big(\mathbf{U}(X),\mathbf{V}(X)\big)_{\varphi(X)}\,\mathrm{d}X,	
\]
where $g$ is a Riemannian metric on $M$.

The material version of the variational formalism for nonequilibrium thermodynamics of viscous heat conducting fluids is given as follows (see \citet[section 3.1]{GBYo2017b}). 
\begin{mdframed}[style=box,frametitle={Variational principle for viscous heat conducting fluids, material version:}]
	The motion $t\mapsto \big(\varphi_t,S(t,\cdot),\Gamma(t,\cdot),\Sigma(t,\cdot)\big)\in\Diff_0(M)\times\Den(M)\times\smooth{M}\times\Den(M)$ is critical for the \emph{variational condition}:
	\begin{equation}\label{eq:material-variational-condition}
	\delta\int_0^T\bigg[\Lag_{\rho_\text{ref}}(\varphi,\dot{\varphi},S)+\underbrace{\int_M(S-\Sigma)\dot{\Gamma}\,\mathrm{d}X}_{\text{thermal coupling term}}\bigg]\mathrm{d}t=0,
	\end{equation}
	subject to the \emph{phenomenological} and \emph{variational} constraints:
	\begin{align}
	\frac{\delta\Lag_{\rho_\text{ref}}}{\delta S}\dot{\Sigma}&=-\mathbf{P}^\text{fr}:\nabla\dot\varphi+\mathbf{J}_S\cdot\dif\dot{\Gamma}-\rho_\text{ref}R,\label{eq:material-phenomenological-constraint}\\
	\frac{\delta\Lag_{\rho_\text{ref}}}{\delta S}\delta\Sigma&=-\mathbf{P}^\text{fr}:\nabla\delta\varphi+\mathbf{J}_S\cdot\dif\delta\Gamma,\label{eq:material-variational-constraint}
	\end{align}
	where $\mathbf{P}^\text{fr}(t,X)$ is the friction Piola-Kirchhoff tensor (the deviatoric part of the first Piola-Kirchhoff tensor $\mathbf{P}(t,X)$), $\mathbf{J}_S(t,X)$ is an entropy flux density and $\rho_\text{ref}(X)R(t,X)$ is the heat power supply density. The dot in $\dot{\Gamma}$ and $\dot{\Sigma}$ is a shortcut for the material derivative $\frac{\mathrm{D}}{\mathrm{D}t}$. In \eqref{eq:material-variational-condition}, the variation of the action functional has to be taken with respect to variations $\delta\varphi$, $\delta S$, $\delta\Sigma$ and $\delta\Gamma$, such that $\delta\varphi$ and $\delta\Gamma$ only vanish at the endpoints $t=0,T$.
\end{mdframed}
Notice how the variational constraint is automatically obtained from the phenomenological one: all time derivatives are replaced by variations and the external heat term is removed. Note also that although the variables $\Gamma$ and $\Sigma$ are needed in the variational formalism, they will vanish when deriving the Euler-Lagrange equations associated to the above variational principle, a task we tackle now. Taking variations in the variational condition \eqref{eq:material-variational-condition}, we obtain:

\begin{align*}
	&\delta\int_0^T\bigg[\Lag_{\rho_\text{ref}}(\varphi,\dot{\varphi},S)+\int_M(S-\Sigma)\dot{\Gamma}\,\mathrm{d}X\bigg]\mathrm{d}t\\
	&\quad=\int_0^T\int_M\bigg[\underbrace{\frac{\partial\Lagd_{\rho_\text{ref}}}{\partial\varphi}\cdot\delta\varphi+\frac{\partial\Lagd_{\rho_\text{ref}}}{\partial\dot{\varphi}}\cdot\delta\dot{\varphi}+\frac{\partial\Lagd_{\rho_\text{ref}}}{\partial T_X\varphi}:\delta T_X\varphi}_{\circled{1}}+\underbrace{\frac{\partial\Lagd_\text{ref}}{\partial S}\delta S}_{\circled{2}}+\underbrace{(S-\Sigma)\delta\dot{\Gamma}}_{\circled{4}}+\underbrace{\delta S\dot{\Gamma}}_{\circled{2}}\underbrace{-\delta\Sigma\dot{\Gamma}}_{\circled{3}}\bigg]\mathrm{d}X\,\mathrm{d}t.
\end{align*}

Now terms in \circled{1} are handled exactly in the same way as in the compressible Euler equation; they transform into:
\[
	\int_0^T\int_M\left(-\frac{\mathrm{D}}{\mathrm{D}t}\frac{\partial\Lagd_{\rho_\text{ref}}}{\partial\dot{\varphi}}-\operatorname{DIV}\frac{\partial\Lagd_{\rho_\text{ref}}}{\partial T_X\varphi}\right)\cdot\delta\varphi\,\mathrm{d}X\,\mathrm{d}t.
\]
The terms in \circled{2} combine together to yield $\frac{\partial\Lagd_{\rho_\text{ref}}}{\partial S}=-\dot{\Gamma}:=-\Theta$; the smooth map $\Theta:M\to\R$ is called the \emph{material temperature} of the fluid, which is the usual temperature scalar field written in material coordinates. Then we use the variational constraint \eqref{eq:material-variational-constraint} with the term \circled{3}, which yields using the divergence theorem:
\[
	\int_0^T\int_M\left[\operatorname{DIV}\mathbf{P}^\text{fr}\cdot\delta\varphi\underbrace{-(\operatorname{DIV}\mathbf{J}_S)\delta\Gamma}_{\circled{4}}\right]\mathrm{d}X\,\mathrm{d}t.
\]
Consequently, we obtain the evolution equations:
\[
	\frac{\mathrm{D}}{\mathrm{D}t}\frac{\partial\Lagd_{\rho_\text{ref}}}{\partial\dot{\varphi}}=\operatorname{DIV}\mathbf{P^\text{cons}}+
	\operatorname{DIV}\mathbf{P}^\text{fr}=\operatorname{DIV}\mathbf{P},\quad\mathbf{P}^\text{cons}:=-\frac{\partial\Lagd_{\rho_\text{ref}}}{\partial{T_X\varphi}}.
\]
For the Lagrangian \eqref{eq:material-Lagrangian} we obtain the equation (balance of momentum):
\[
	\rho_\text{ref}\frac{\mathrm{D}\mathbf{U}}{\mathrm{D}t}=-(\grad p)\circ(\varphi J_\varphi)+\operatorname{DIV}\mathbf{P}^\text{fr},
\]
which together with the phenomenological constraint \eqref{eq:material-phenomenological-constraint}, the continuity equation $\rho=\rho_\text{ref}\varphi^{-1}$ (balance of mass in material coordinates), and the no-slip boundary conditions ($\varphi(t,X)=X$ for any $X\in\partial M$ at any time $t$, which amounts to say that $\mathbf{U}|_{\partial M}=0$), constitute the Navier-Stokes-Fourier system of partial differential equations, in material coordinates. Although it is seldom used in fluid mechanics, the material variational principle is important for the purpose of this article since the classical Navier-Stokes-Fourier system of partial differential equations arises from the \emph{reduction} of this principle, and, as we will see in section \ref{sec:variational-discretization-Navier-Stokes-Fourier}, it is important to obtain the right variational discretization of the classical Navier-Stokes-Fourier system.

Returning to the computation above, we see that terms in \circled{4} yield the identity $\dot{S}=-\operatorname{DIV}\mathbf{J}_S+\dot{\Gamma}$. This identity allows us to find the balance of entropy in material coordinates for the fluid with the help of the phenomenological constraint \eqref{eq:material-phenomenological-constraint}:
\[
	\dot{S}+\operatorname{DIV}\mathbf{J}_S=\frac{1}{\Theta}\left(\mathbf{P}^\text{fr}:\nabla\dot{\varphi}-\mathbf{J}_S\cdot\dif\Theta+\rho_\text{ref}R\right).
\]
Therefore the entropy changes in the system are due to the internal irreversible processes of friction (viscosity) and heat conduction, but also to the external heat supply.

\begin{remark}
	The thermal coupling term can be written as the inner product $\big\langle S-\Sigma,\dot{\Gamma}\big\rangle$ given in \eqref{eq:density-function-inner-product} between the density $S-\Sigma$ and the function $\dot{\Gamma}$. This remark will be useful in the discretization of the variational principle in section \ref{sec:semi-discrete}.
\end{remark}

\begin{remark}[Geometric structure of the phenomenological constraint]
	We now relate the material variational principle composed of \eqref{eq:material-variational-condition}, \eqref{eq:material-phenomenological-constraint} and \eqref{eq:material-variational-constraint}, to variational formulations in nonholonomic mechanics. Let us consider formally the manifold ${Q}:=\Diff_0(M)\times\Den(M)\times\smooth{M}\times\Den(M)$ and denote an element in ${Q}$ by $q:=(\varphi,S,\Gamma,\Sigma)$. We shall write $(q,\dot{q},\delta q)$ for an element in the fiber over $q\in Q$ of $TQ\times_Q TQ$. The variational constraint \eqref{eq:material-variational-constraint} defines a subset $C_\text{var}\subset TQ\times_Q TQ$ as follows: $(q,\dot{s},\delta q)\in C_\text{var}$ if and only if $(q,\dot{q},\delta q)$ satisfies \eqref{eq:material-variational-constraint},
	where $(q,\dot{q},\delta q)= (\varphi,S,\Gamma,\Sigma,\dot{\varphi},\dot{S},\dot{\Gamma},\dot{\Sigma},\delta\varphi,\delta S,\delta\Gamma,\delta\Sigma)$. This variational constraint satisfies the following property: for each $(q,\dot{q})\in TQ$, the set $C_\text{var}(q,\dot{q}):=C_\text{var}\cap\big[\{(q,\dot{q})\}\times T_qQ\big]$ is a vector space. The phenomenological constraint \eqref{eq:material-phenomenological-constraint} on $(q,\dot{q})= (\varphi,S,\Gamma,\Sigma,\dot{\varphi},\dot{S},\dot{\Gamma},\dot{\Sigma})$ also defines a subset $C_\text{pheno} \subset TQ$. For the case of adiabatically closed systems (i.e., $\rho_\text{ref}R=0$), the subset $C_\text{pheno}$ is obtained from the variational constraint $C_\text{var}$ via the following general construction:
	\begin{equation}\label{eq:Cpheno-Cvar} 
		C_\text{pheno}=\big\{(q,\dot{q})\in TQ:(q,\dot{q})\in C_\text{var}(q,\dot{q})\big\}.
	\end{equation} 
	In terms of these constraint sets, the material variational formulation of the Navier-Stokes-Fourier system can thus be written as follows: a curve $t\mapsto q(t)=(\varphi(t),S(t),\Gamma(t),\Sigma(t))\in Q$ satisfies the material variational principle if and only if it satisfies the variational condition
	\begin{equation}
		\delta \int_0^T\widetilde{\Lag}_{\rho_\text{ref}}(q,\dot{q})\,\mathrm{d}t=0,\quad
		\widetilde{\Lag}_{\rho_\text{ref}}(q,\dot{q}):=\Lag_{\rho_\text{ref}}(\varphi,\dot{\varphi},S)+\int_M(S-\Sigma)\dot{\Gamma}\,\mathrm{d}X,
	\end{equation} 
	for all variations $\delta q\in C_\text{var}(q,\dot{q})$ vanishing at the endpoints and curves $t\mapsto q(t)\in C_\text{pheno}$. Therefore, the material variational principle is an infinite-dimensional, \emph{nonlinear} extension of the Lagrange-d'Alembert principle used for the treatment of nonholonomic mechanical systems with \emph{linear} constraints, see e.g., \cite{Bl2003}. Such linear constraints are given by a \emph{distribution} $\Delta\subset TQ$ and the associated variational constraint is given by $\Delta\times_Q TQ$. For the case of \emph{nonlinear} constraints a generalization of the Lagrange-d'Alembert principle has been considered in \cite{Ch1934}; see also \cite{Ap1911}, \cite{Pi1983}. In Chetaev's approach, the variational constraint is derived from the kinematic constraint. However, it has been pointed out in \cite{Ma1998} that this principle does not always lead to the correct equations of motion and in general one has to consider the kinematic and variational constraints as independent notions. A general geometric variational approach for nonholonomic systems with nonlinear and (possibly) higher order kinematic and variational constraints has been described in \cite{CeIbdLdD2004}. This setting generalizes both the Lagrange-d'Alembert and Chetaev approaches. It is important to point out that for these generalizations, including Chetaev's approach, energy may not be conserved along the solution of the equations of motion. The material variational principle and more generally the variational principle for nonequilibrium thermodynamics developed in \cite{GBYo2017a} and \cite{GBYo2017b} fall into the general setting described in \cite{CeIbdLdD2004}, extended to the infinite dimensional setting. However, for these particular constraints of thermodynamic type, total energy is conserved, consistently with the fact that we supposed the system to be isolated.
\end{remark}

%---------------------------------------------------------------------------------------------------

\subsection{The spatial variational principle for the Navier-Stokes-Fourier system}\label{ssec:smooth-spatial}

Just like we did in \ref{ssec:smooth-Euler-Poincaré} for the case of compressible fluid dynamics, in \citet[section 3.2]{GBYo2017c}, a reduction akin to Euler-Poincaré reduction is performed onto the material variational principle for viscous heat conducting fluids. Consequently a spatial version of the variational principle is obtained, whose Euler-Lagrange equations will turn out to be the classical Navier-Stokes-Fourier system of PDE. 

As before, we consider the group $\Diff_0(M)$ that acts on itself on the right, as well as on densities, functions, and other tensors. Its Lie algebra is the Lie algebra $\LieAlgebraFont{X}_0(M)$ of vector fields $\mathbf{u}$ on $M$ that satisfy the no-slip boundary condition $\mathbf{u}|_{\partial M}=0$. The reduced variables are $\mathbf{u}=\dot{\varphi}\circ\varphi^{-1}\in\LieAlgebraFont{X}_0(M)$ (the fluid velocity), $\rho=(\rho_\text{ref}\circ\varphi^{-1})J_{\varphi^{-1}}\in\Den(M)$ (the mass density of the fluid), $s=(S\circ\varphi^{-1})J_{\varphi^{-1}}\in\Den(M)$ (the entropy density of the fluid), $r=R\circ\varphi^{-1}$ (the external heat power) and the auxiliary new variables $\gamma=\Gamma\circ\varphi^{-1}\in\smooth{M}$ and $\sigma=(\Sigma\circ\varphi^{-1})J_{\varphi^{-1}}\in\Den(M)$. The reduced analogues of $\mathbf{P}^\text{fr}$ and $\mathbf{J}_S$ are defined via the inverse Piola transform\footnote{Recall that the Piola transform is a map that is used to transform tensors given in Eulerian (spatial) coordinates to tensors given in Lagrangian (material coordinates), see \cite[Defintion 7.18]{MaHu1983} and below for details.}:
\[
	\mathbf{j}_S=(\varphi_*\mathbf{J}_S)J_{\varphi^{-1}},\quad\bm{\sigma}^\text{fr}(x)(\alpha_x,\beta_x)=\mathbf{P}^\text{fr}(X)\big(\alpha_x,T_X^*\varphi(\beta_x)\big),
\]
where $x=\varphi(X)\in M$ and $\alpha_x$, $\beta_x\in T_x^*M$.

%\BC{Je ne comprends pas l'utilisation de la notation $T_X^*\varphi$. Est-ce que ce ne serait pas plutôt l'application $T_{\varphi(X)}^*\varphi^{-1}:T_{\varphi(X)}^*M\to T_X M$?}
%\FGB{Non, la convention de la notation est $T_xf:T_xM\rightarrow T_{f(x)}N$ et $T_x^*f:T_{f(x)}^*N\rightarrow T_x^*M$}

As the material Lagrangian \eqref{eq:material-Lagrangian} is $\Diff_0(M)_{\rho_\text{ref}}$-invariant, we obtained a reduced Lagrangian called the \emph{spatial Lagrangian}, which is given by:
\begin{equation}\label{eq:spatial-Lagrangian}
	\lag(\mathbf{u},\rho,s)=\int_M\left[\frac{1}{2}\rho\langle\mathbf{u}^\flat,\mathbf{u}\rangle-\varepsilon(\rho,s)\right]\mathrm{d}x.
\end{equation}
Also, in \eqref{eq:material-phenomenological-constraint} and \eqref{eq:material-variational-constraint}, the material derivatives $\frac{\mathrm{D}\Gamma}{\mathrm{D}t}(t,X)$ and $\frac{\mathrm{D}\Sigma}{\mathrm{D}t}(t,X)$ are replaced as follows:
\begin{align*}
	\frac{\mathrm{D}\Gamma}{\mathrm{D}t}(t,X)&=\frac{\mathrm{d}\gamma}{\mathrm{d}t}\big(t,\varphi_t(X)\big)=\frac{\partial\gamma}{\partial t}(t,x)+\dif\gamma_x(\mathbf{u}),\\
	\frac{\mathrm{D}\Sigma}{\mathrm{D}t}(t,X)&=\frac{\mathrm{d}}{\mathrm{d}t}\left[\sigma\big( t,\varphi_t(X)\big)J_{\varphi_t}(X)\right]=\left(\frac{\partial\sigma}{\partial t}\big(t,x\big)+\dive(\sigma\mathbf{u})(x)\right)J_{\varphi_t}(X),
\end{align*}
where $x(t,X)=\varphi_t(X)\in M$. The spatial version of the variational formalism for nonequilibrium thermodynamics of viscous heat conducting fluids is given by (see \citet[section 3.2]{GBYo2017b}) as follows.

\begin{mdframed}[style=box,frametitle={Variational principle for viscous heat conducting fluids, spatial version:}]
	The motion $t\mapsto\big(\mathbf{u}_t,\rho(t,\cdot),s(t,\cdot),\gamma(t,\cdot),\sigma(t,\cdot)\big)\in\LieAlgebraFont{X}_0(M)\times\Den(M)\times\Den(M)\times\smooth{M}\times\Den(M)$ is critical for the \emph{variational condition}:
	\begin{equation}\label{eq:spatial-variational-condition}
		\delta\int_0^T\bigg[\lag(\mathbf{u},\rho,s)+\underbrace{\int_M(s-\sigma)\big(\partial_t{\gamma}+\dif\gamma(\mathbf{u})\big)\,\mathrm{d}x}_{\text{thermal coupling term}}\bigg]\mathrm{d}t=0,
	\end{equation}
	subject to the \emph{phenomenological} and \emph{variational} constraints:
	\begin{align}
		\frac{\delta\lag}{\delta s}\big(\partial_t{\sigma}+\dive(\sigma\mathbf{u})\big)&=-\bm{\sigma}^\text{fr}:\nabla\mathbf{u}+\mathbf{j}_S\cdot\dif\big(\partial_t{\gamma}+\dif\gamma(\mathbf{u})\big)-\rho r,\label{eq:spatial-phenomenological-constraint}\\
		\frac{\delta\lag}{\delta s}\big(\delta\sigma+\dive(\sigma\bm{\zeta})\big)&=-\bm{\sigma}^\text{fr}:\nabla\bm{\zeta}+\mathbf{j}_S\cdot\dif\big(\delta\gamma+\dif\gamma(\bm{\zeta})\big),\label{eq:spatial-variational-constraint}
	\end{align}
	where $\bm{\sigma}^\text{fr}(t,x)$ is the friction Cauchy stress tensor (the deviatoric part of the Cauchy stress tensor $\bm{\sigma}(t,x)$), $\mathbf{j}_S(t,x)$ is an entropy flux density and $\rho r$ is the external heat power. In \eqref{eq:spatial-variational-condition}, the variation of the action functional has to be taken with respect to variations $$\delta\mathbf{u}=\partial_t\bm{\zeta}+[\mathbf{u},\bm{\zeta}],\quad\delta\rho=-\dive(\rho\bm{\zeta}),$$ $\delta s$, $\delta\sigma$ and $\delta\gamma$, such that $\bm{\zeta}\in\LieAlgebraFont{X}_0(M)$ and $\delta\gamma$ both vanish at the endpoints $t=0,T$, and $\delta\gamma|_{\partial M}=0$.
\end{mdframed}

Notice again how the variational constraint is automatically obtained from the phenomenological one: all time derivatives are replaced by variations, $\mathbf{u}$ is replaced by $\bm{\zeta}$ and the external heat term is removed. Note also that although the variables $\gamma$ and $\sigma$ are needed in the variational formalism, they will vanish when deriving the Euler-Lagrange equations associated to the above variational principle.

We are now going to derive the Euler-Lagrange equations associated to the spatial variational principle, and relate them directly to the Navier-Stokes-Fourier system of PDE. We begin by taking the variations of the action functional in \eqref{eq:spatial-variational-condition}:
\begin{align*}
	\delta\int_0^T\bigg[\lag(\mathbf{u},\rho,s)&+\int_M(s-\sigma)\big(\partial_t{\gamma}+\dif\gamma(\mathbf{u})\big)\,\mathrm{d}x\bigg]\mathrm{d}t\\
	=\int_0^T&\bigg[\underbrace{\left\langle\frac{\delta\lag}{\delta\mathbf{u}},\delta\mathbf{u}\right\rangle+\left\langle\frac{\delta\lag}{\delta\rho},\delta\rho\right\rangle}_{\circled{5}}+\underbrace{\left\langle\frac{\delta\lag}{\delta s},\delta s\right\rangle}_{\circled{1}}\\
	&\quad+\underbrace{\int_M\delta s\big(\partial_t\gamma+\dif\gamma(\mathbf{u})\big)\,\mathrm{d}x}_{\circled{1}}\underbrace{-\int_M\delta\sigma\big(\partial_t\gamma+\dif\gamma(\mathbf{u})\big)\,\mathrm{d}x}_{\circled{2}}\\
	&\quad+\underbrace{\int_M(s-\sigma)\partial_t(\delta\gamma)\,\mathrm{d}x}_{\circled{3}}+\underbrace{\int_M(s-\sigma)\dif(\delta\gamma)(\mathbf{u})\,\mathrm{d}x}_{\circled{3}}+\underbrace{\int_M(s-\sigma)\dif\gamma(\delta u)\,\mathrm{d}x}_{\circled{4}}\bigg],
\end{align*}
and this expression has to be equal to zero for all variations $\delta\mathbf{u}$, $\delta\rho$, $\delta s$ and $\delta\gamma$ satisfying all the variational constraints mentioned above. Then, all terms that couple with the variation $\delta s$ are in \circled{1}, therefore we obtain $\frac{\delta\lag}{\delta s}=-\big(\partial_t\gamma+\dif\gamma(\mathbf{u})\big)=:-\theta=-\Theta\circ\varphi^{-1}$ which is the \emph{spatial temperature} of the fluid. Then, in \circled{2}, we use the spatial variational constraint \eqref{eq:spatial-variational-constraint} to get:
\begin{align*}
	-\int_M\delta\sigma\big(\partial_t\gamma+\dif\gamma(\mathbf{u})\big)\,\mathrm{d}x
	&=\int_M\frac{\delta\lag}{\delta s}\delta\sigma\,\mathrm{d}x\\
	&=-\int_M\frac{\delta\lag}{\delta s}\dive(\sigma\bm{\zeta})\mathrm{d}x
	-\int_M\bm{\sigma}^\text{fr}:\nabla\bm{\zeta}\,\mathrm{d}x+\int_M\mathbf{j}_S\cdot\dif\big(\delta\gamma+\dif\gamma(\bm{\zeta})\big)\,\mathrm{d}x.
\end{align*}
Now using the well-known formula $\dive(f\mathbf{u})=f\dive\mathbf{u}+\dif f(\mathbf{u})$ and the traditional divergence theorem, the generalized divergence theorem for tensors (see \cite[Problem 7.6]{MaHu1983}), as well as the boundary conditions $\bm{\zeta}|_{\partial M}=0$ as well as $\delta\gamma|_{\partial M}=0$\footnote{In the case where the constraint $\delta\gamma|_{\partial M}=0$ is not imposed, then the variational formalism above yields the condition $\mathbf{j}_S\cdot\mathbf{n}=0$ on $\partial M$, and if in addition the external heat power supply is zero, then the fluid is adiabatically closed.}, we obtain:
\[
	-\int_M\delta\sigma\big(\partial_t\gamma+\dif\gamma(\mathbf{u})\big)\,\mathrm{d}x
	=\underbrace{\int_M\sigma\dif\left(\frac{\delta\lag}{\delta s}\right)(\bm{\zeta})\,\mathrm{d}x}_{\circled{4}}+\underbrace{\int_M(\dive\bm{\sigma}^\text{fr})(\bm{\zeta})\,\mathrm{d}x}_{\circled{5}}\underbrace{-\int_M(\dive\mathbf{j}_S)\delta\gamma\,\mathrm{d}\mathbf{x}}_{\circled{3}}\underbrace{-\int_M(\dive{\mathbf{j}_S})\dif\gamma(\bm{\zeta})\,\mathrm{d}x}_{\circled{4}}.
\]
After using an integration by parts relatively to time as well as the divergence theorem, we get from the terms in \circled{3} that couple with the variation $\delta\gamma$ the following relation:
\begin{equation}\label{eq:eq1}
	\partial_t(s-\sigma)+\dive\big((s-\sigma)\mathbf{u}\big)=-\dive\mathbf{j}_S.
\end{equation}
It remains to analyze the terms in \circled{4}. Again, using integration by parts, the divergence theorem and boundary conditions, we obtain:
\begin{align*}
	\int_0^T\bigg[\int_M\sigma&\dif\left(\frac{\delta\lag}{\delta s}\right)(\bm{\zeta})\,\mathrm{d}x+\int_M(s-\sigma)\dif\gamma(\delta u)\,\mathrm{d}x-\int_M(\dive{\mathbf{j}_S})\dif\gamma(\bm{\zeta})\,\mathrm{d}x\bigg]\mathrm{d}t\\
	&=\int_0^T\!\!\int_M\sigma\dif\left(\frac{\delta\lag}{\delta s}\right)(\bm{\zeta})\,\mathrm{d}x\,\mathrm{d}t+\int_0^T\!\!\int_M\big[\underbrace{-\partial_t(s-\sigma)-\dive\big((s-\sigma)u\big)-\dive\mathbf{j}_S}_{=0\text{ thanks to \eqref{eq:eq1}}}\big]\dif\gamma(\bm{\zeta})\,\mathrm{d}t\\
	&\quad\quad\quad-\int_0^T\!\!\int_M(s-\sigma)\dif\theta(\bm{\zeta})\,\mathrm{d}x\,\mathrm{d}t\\
	&=\underbrace{\int_0^T\!\!\int_M s\frac{\delta\lag}{\delta s}(\bm{\zeta})\,\mathrm{d}x\,\mathrm{d}t}_{\circled{5}}.
\end{align*}
Finally, after collecting all terms in \circled{5} that couple with $\bm{\zeta}$ (some of them have already been computed in section \ref{ssec:smooth-Euler-Poincaré} explicitly) we obtain the following equation:
\[
	(\partial_t+\ad_\mathbf{u}^*)\frac{\delta\lag}{\delta\mathbf{u}}=\rho\,\dif\frac{\delta\lag}{\delta\rho}+s\,\dif\frac{\delta\lag}{\delta s}+\dive\bm{\sigma}^\mathrm{fr}.
\]
This is the Euler-Lagrange equation associated to the spatial (reduced) variational principle, which physically corresponds to the balance of momentum. Together with the advection equation for the mass density (the continuity equation) which is obtained by reduction of $\rho=(\rho_\text{ref}\circ\varphi^{-1})J_{\varphi^{-1}}$, the phenomenological constraint \eqref{eq:spatial-phenomenological-constraint} (the local balance of entropy), and using the relation \eqref{eq:eq1} to make $\sigma$ disappear, we obtain the following \emph{abstract} Navier-Stokes-Fourier system of equations:
\begin{equation}\label{eq:smooth-abstract-Navier-Stokes-Fourier}
	\begin{cases}
		\displaystyle(\partial_t+\ad_\mathbf{u}^*)\frac{\delta\lag}{\delta\mathbf{u}}=\rho\,\dif\frac{\delta\lag}{\delta\rho}+s\,\dif\frac{\delta\lag}{\delta s}+\dive\bm{\sigma}^\mathrm{fr}\\[7pt]
		\partial_t\rho+\dive(\rho\mathbf{u})=0\\[7pt]
		\displaystyle\frac{\delta\lag}{\delta s}\big(\partial_t{s}+\dive(s\mathbf{u})+\dive\mathbf{j}_S\big)=-\bm{\sigma}^\text{fr}:\nabla\mathbf{u}+\mathbf{j}_S\cdot\dif\theta-\rho r
	\end{cases}
\end{equation}
together with the non-slip boundary condition $\mathbf{u}|_{\partial M}=0$ and also appropriate initial conditions. Notice how the auxiliary variables $\gamma$ and $\sigma$ have disappeared from this system of partial differential equations, which describes the evolution of the state $(\mathbf{u},\rho,s)$ of the fluid with respect to time in accordance with the principles of thermodynamics. Note also that the three partial differential equations are fully coupled.

In order to recover the \emph{classical} Navier-Stokes-Equations, we need phenomenological laws that give us expressions for the symmetric tensor $\bm{\sigma}^\mathrm{fr}$ and the vector field $\mathbf{j}_S$. We choose to focus on Newtonian fluids that conduct heat according to Fourier's law, therefore:
\begin{equation}\label{eq:Newton-Fourier}
	\bm{\sigma}:=-p\Id+\bm{\sigma}^\mathrm{fr},\quad
	\bm{\sigma}^\mathrm{fr}:=2\mu\Def\mathbf{u}+\left(\zeta-\frac{2}{3}\mu\right)(\dive\mathbf{u})\Id,\quad
	\mathbf{j}_Q^\flat:=\theta\mathbf{j}_S^\flat=-\lambda\dif\theta,
\end{equation}
with $\mu\geq 0$ the shear viscosity, $\zeta\geq 0$ the bulk viscosity, $\lambda\geq 0$ the thermal conductivity and $p$ the thermodynamic pressure (given below); and where $\Def\mathbf{u}=\frac{1}{2}(\nabla\mathbf{u}+\nabla\mathbf{u}^\mathsf{T})$ is the linear strain tensor and $\mathbf{j}_Q$ the heat flux density. Using these expressions as well as computing the variational derivatives in \eqref{eq:smooth-abstract-Navier-Stokes-Fourier} of the spatial Lagrangian \eqref{eq:spatial-Lagrangian} and sharpening, we finally recover the classical Navier-Stokes-Fourier system of partial differential equations:
\begin{equation}\label{eq:smooth-classical-Navier-Stokes-Fourier}
	\begin{cases}
		\displaystyle\rho\left(\frac{\partial\mathbf{u}}{\partial t}+\nabla_{\mathbf{u}}\mathbf{u}\right)=-\grad p+\mu\Delta\mathbf{u}+\left(\zeta-\frac{2}{3}\mu\right)\grad\dive\mathbf{u}\\[7pt]
		\displaystyle\frac{\partial\rho}{\partial t}+\dive(\rho\mathbf{u})=0\\[7pt]
		\displaystyle\frac{\partial s}{\partial t}+\dive(s\mathbf{u})=\frac{1}{\theta}\left(\bm{\sigma}^\textrm{fr}:\nabla\mathbf{u}+\lambda\Delta\theta+\rho r\right)
	\end{cases}
\end{equation}
where the thermodynamic pressure $p$ is given by $p=\frac{\partial\varepsilon}{\partial\rho}\rho+\frac{\partial\varepsilon}{\partial s}s-\varepsilon$, $\Delta\mathbf{u}$ is the vector Laplacian of the velocity vector field $\mathbf{u}$, and $\Delta\theta$ is the scalar Laplacian of the temperature scalar field $\theta$\footnote{Note that in the more general case of a Riemannian manifold $(M,g)$, the vector Laplacian is replaced by the Hodge-De Rham Laplacian $\Delta=\bm{\delta}\circ\dif+\dif\circ\bm{\delta}$, and the scalar Laplacian is replaced by the Laplace-Beltrami operator. This can be observed in \cite{GBYo2017b}, therefore answering the question of the appropriate generalization of the Navier-Stokes equations to an arbitrary Riemannian manifold (see \cite{ChCzDi2017}).}.

\begin{remark}[Other fluid models]
	More sophisticated fluid models can be used. Starting from the abstract Navier-Stokes-Fourier equations \eqref{eq:smooth-abstract-Navier-Stokes-Fourier} and specifying different a friction stress $\bm{\sigma}^\mathrm{fr}$, entropy flux $\mathbf{j}_S$ and internal energy $\varepsilon$ will yield other partial differential equations for these fluid models.
\end{remark}

\begin{remark}[External forces]
	Additional external forces (that are not friction-type forces and non conservative) can be added to the variational principle in the same way as external forces are added to Hamilton's principle to yield the Lagrange-d'Alembert principle.
\end{remark}

\begin{remark}[Free energy Lagrangian]
	We have seen that at the end of the derivation of the constrained Euler-Lagrange equations associated to either the material or the variational principle, the phenomenological constraint \eqref{eq:spatial-phenomenological-constraint} becomes the local balance of entropy in the system of equations \eqref{eq:smooth-abstract-Navier-Stokes-Fourier}. Using the \emph{free energy Lagrangian} instead of the Lagrangian \eqref{eq:material-Lagrangian}, one can obtain at the end the \emph{heat equation} for the fluid instead of the local balance of entropy, which should be more practical in applications (see \cite{GBYo2017c}).
\end{remark}

\begin{corollary}[Local and global entropy balance, meaning of $\sigma$]\label{cor:smooth-spatial-entropy-balance}
	Thanks to the relation \eqref{eq:eq1} and the phenomenological constraint \eqref{eq:spatial-phenomenological-constraint}, we have
	\begin{align*}
		\partial_t s+\dive(s\mathbf{u})+\dive(\mathbf{j}_S)&=\partial_t\sigma+\dive(\sigma\mathbf{u})=\underbrace{\frac{1}{\theta}(\bm{\sigma}^\textrm{fr}:\nabla\mathbf{u}-\mathbf{j_S}\cdot\dif\theta)}_{\circled{1}}+\underbrace{\frac{\rho r}{\theta}}_{\circled{2}},
	\end{align*}
%	\FGB{Je ne comprends pas cette \'egalit\'e. Moi j'ai
%	\begin{align*}
%	\partial_t s+\dive(s\mathbf{u})+\dive(\mathbf{j}_S)&=\partial_t\sigma+\dive(\sigma\mathbf{u})= \frac{1}{\theta} \left(\bm{\sigma}^\textrm{fr}:\nabla\mathbf{u}-\mathbf{j}_S\cdot\mathbf{d}\theta \right)\\
%	&=\frac{1}{\theta} \left(\bm{\sigma}^\textrm{fr}:\nabla\mathbf{u} + \frac{\lambda}{\theta}\mathbf{d}\theta \cdot\mathbf{d}\theta \right),
%	\end{align*}}
	where \circled{1} is the \emph{internal} entropy production (which can be further separated into viscosity-related and heat-related contributions) and \circled{2} is the \emph{external} entropy production. Therefore, we can interpret $\dot{\sigma}$ as the rate at which (total) entropy is produced in the fluid. The global balance of entropy is the classical one for nonequilibrium thermodynamics of continuum systems: define $\mathsf{S}(t)=\int_M s(t,x)\,\mathrm{d}x$ then applying the Reynold's transport theorem to $\mathsf{S}$ and using the spatial phenomenological constraint \eqref{eq:spatial-phenomenological-constraint} we obtain:
	\[
		\frac{\mathrm{d}\mathsf{S}}{\mathrm{d}t}+\int_{\partial M}\mathbf{j}_S\cdot\mathrm{d}\mathbf{S}=\underbrace{\int_M\frac{1}{\theta}(\bm{\sigma}^\mathrm{fr}:\nabla\mathbf{u}-\mathbf{j}_S\cdot\dif\theta)\,\mathrm{d}x}_{\text{internal source}}+\underbrace{\int_M\frac{\rho r}{\theta}\,\mathrm{d}x}_{\text{external source}}.
	\]
	Note that this balance can be obtained by integration over $M$ of the local one.
\end{corollary}

\begin{corollary}[Local and global energy balance]\label{cor:smooth-spatial-energy-balance}
	We can compute the global balance of energy as it is usually done for Lagrangian systems and in particular Euler-Poincaré systems (see \cite[section 1.3]{CoGB2018}). The spatial energy density of the fluid is defined as the Legendre transform of the spatial Lagrangian:  $\mathsf{E}(\mathbf{u},\rho,s)=\left\langle\frac{\delta\lag}{\delta\mathbf{u}}(\mathbf{u},\rho,s),\mathbf{u}\right\rangle-\lag(\mathbf{u},\rho,s)$. For the spatial Lagrangian \eqref{eq:spatial-Lagrangian} we obtain $\mathsf{E}(\mathbf{u},\rho,s)=\int_M\left[\frac{1}{2}\rho\langle\mathbf{u}^\flat,\mathbf{u}\rangle+\varepsilon(\rho,s)\right]\mathrm{d}x$ which is precisely the \emph{total} energy of the fluid: kinetic energy plus internal energy. Then we compute the time derivative of $\mathsf{E}$ along a solution of \eqref{eq:smooth-abstract-Navier-Stokes-Fourier} as follows:
	\begin{align*}
		\frac{\mathrm{d}\mathsf{E}}{\mathrm{d}t}
		&=\left\langle\frac{\partial}{\partial t}\frac{\delta\lag}{\delta\mathbf{u}},\mathbf{u}\right\rangle-\left\langle\frac{\delta\lag}{\delta\rho},\dot{\rho}\right\rangle-\left\langle\frac{\delta\lag}{\delta s},\dot{s}\right\rangle \\
		&=\left\langle\frac{\partial}{\partial t}\frac{\delta\lag}{\delta\mathbf{u}}-\rho\dif\frac{\delta\lag}{\delta\rho},\mathbf{u}\right\rangle+\langle\theta,\dot{s}\rangle\\
		&=\left\langle-\ad_\mathbf{u}^*\frac{\delta\lag}{\delta\mathbf{u}}+\dive\bm{\sigma}^\text{fr},\mathbf{u}\right\rangle+\big\langle\theta,\partial_t{s}+\dive(s\mathbf{u})\big\rangle\\
		&=\langle\dive\bm{\sigma}^\mathrm{fr},\mathbf{u}\rangle+\langle-\dive\bm{\sigma}^\mathrm{fr},\mathbf{u}\rangle-\int_M\dive\mathbf{j}_Q\,\mathrm{d}x+\int_M\rho r\,\mathrm{d}x\\
		&=-\int_{\partial M}\mathbf{j}_Q\cdot\mathrm{d}\mathbf{S}+\int_M\rho r\,\mathrm{d}x,
	\end{align*}
	which means that the rate at which the total energy of the fluid changes is determined by the exchanges of heat with its environment. We will see that this global total energy balance will still hold at the discrete level. We could also obtain this global balance from a \emph{local} one. A quite long but straightforward computation using the relations $(\nabla_\mathbf{u}\mathbf{u})^\flat=\lie_\mathbf{u}\mathbf{u}^\flat-\frac{1}{2}\dif\langle\mathbf{u}^\flat,\mathbf{u}\rangle$ and $\dive(\bm{\sigma}\cdot\mathbf{u})=\dive(\bm{\sigma})\cdot\mathbf{u}+\bm{\sigma}:\nabla\mathbf{u}$ yields
	\[
		\partial_t e+\dive(e\mathbf{u})=\dive(\widetilde{\bm{\sigma}}\cdot\mathbf{u})-\dive\mathbf{j}_Q+\rho r,
	\]
	where $\widetilde{\bm{\sigma}}=-\widetilde{p}\,\Id+\bm{\sigma}^\mathrm{fr}$ with $\widetilde{p}:=p-\frac{1}{2}\rho\langle\mathbf{u}^\flat,\mathbf{u}\rangle$ (the so-called \emph{static} pressure). Then applying Reynold's transport theorem to $\mathsf{E}(t)=\int_M e(t,x)\,\mathrm{d}x$ as well as the boundary condition $\mathbf{u}|_{\partial M}=0$, we recover the global balance of energy.
\end{corollary}

\begin{corollary}[Viscous Kelvin's theorem]\label{cor:viscous-Kelvin-theorem}
	For any Euler-Poincaré system, a Kelvin-Noether theorem is available, which is akin to a generalization of the famous Kelvin's theorem in fluid mechanics (see \citet[section 4]{HoMaRa1998}). We now derive a Kelvin-Noether theorem for the abstract system \eqref{eq:smooth-abstract-Navier-Stokes-Fourier} and then in a particular case that will yield the \emph{viscous Kelvin's theorem} of fluid mechanics. Let $\mathcal{C}$ be a manifold; we will suppose that $\Diff_0(M)$ acts on the left of $\mathcal{C}$. Let $\mathcal{K}:\mathcal{C}\times\Orb(\rho_\text{ref})\times\Den(M)\to\R$ be an equivariant map where $\Diff_0(M)$ acts on $\LieAlgebraFont{X}(M)^{**}$ by the dual of the coadjoint action of $\Diff_0(M)$ on $\LieAlgebraFont{X}(M)^*$. Finally let the \emph{Kelvin-Noether quantity} (associated to $\mathcal{C}$ and $\mathcal{K}$) be the map $I:\mathcal{C}\times\LieAlgebraFont{X}(M)\times\Orb(\rho_\text{ref})\times\Den(M)\to\R$ defined by $I(\gamma,\mathbf{u},\rho,s)=\left\langle\mathcal{K}(\gamma,\rho),\frac{\delta\lag}{\delta\mathbf{u}}(\mathbf{u},\rho,s)\right\rangle$. Then along a solution $(\varphi,\mathbf{u},\rho,s)$ of \eqref{eq:smooth-abstract-Navier-Stokes-Fourier}, where the motion $t\mapsto\varphi_t$ satisfies the \emph{reconstruction equation} $\mathbf{u}\circ\varphi_t=\dot{\varphi}_t$, we obtain:
	\[
		\frac{\mathrm{d}I}{\mathrm{d}t}=\big\langle\mathcal{K}(\gamma,\rho),\dive\widetilde{\bm{\sigma}}\big\rangle,
	\]
	with $\gamma:=\gamma_\text{ref}\varphi$, for some $\gamma_\text{ref}\in\mathcal{C}$. Indeed, using the equivariance property of $\mathcal{K}$, as well as \citet[formula (9.3.7)]{MaRa1999}, and the local balance of momentum, we have successively:
	\begin{align*}
		\frac{\mathrm{d}I}{\mathrm{d}t}&=\left\langle\mathcal{K}(\gamma_\text{ref},\rho_\text{ref}),\frac{\mathrm{d}}{\mathrm{d}t}\Ad_{\varphi_t^{-1}}^*\frac{\delta\lag}{\delta\mathbf{u}}(\mathbf{u},\rho,s)\right\rangle\\
		&=\left\langle\mathcal{K}(\gamma_\text{ref},\rho_\text{ref}),\Ad_{\varphi_t^{-1}}^*\left[\ad_\mathbf{u}^*\frac{\delta\lag}{\delta\mathbf{u}}(\mathbf{u},\rho,s)+\frac{\partial}{\partial t}\frac{\delta\lag}{\delta\mathbf{u}}(\mathbf{u},\rho,s)\right]\right\rangle\\
		&=\left\langle\mathcal{K}(\gamma,\rho),\rho\,\dif\frac{\delta\lag}{\delta\rho}+s\,\dif\frac{\delta\lag}{\delta s}+\dive\bm{\sigma}^\mathrm{fr}\right\rangle\\
		&=\big\langle\mathcal{K}(\gamma,\rho),\dive\widetilde{\bm{\sigma}}\big\rangle,
	\end{align*}
	where $\widetilde{\bm{\sigma}}=-\widetilde{p}\,\Id+\bm{\sigma}^\mathrm{fr}$ with $\widetilde{p}:=p-\frac{1}{2}\rho\langle\mathbf{u}^\flat,\mathbf{u}\rangle$. The viscous Kelvin's theorem is recovered by choosing $\mathcal{C}$ to be the space of closed loops $\mathcal{C}^0(S^1)$ and $\mathcal{K}$ to be the circulation map defined by $\big\langle\mathcal{K}(\gamma,\rho),\alpha\big\rangle=\oint_\gamma\frac{\alpha}{\rho}$, where $\frac{\alpha}{\rho}\in\Omega^1(M)$ is defined by $\frac{\alpha}{\rho}=\frac{f}{\rho}\beta$ for $\alpha=\beta\otimes f\,\mathrm{d}x$. Then along a solution of \eqref{eq:smooth-classical-Navier-Stokes-Fourier} we have:
	\[
		\frac{\mathrm{d}I}{\mathrm{d}t}=\oint_\gamma\frac{\mu}{\rho}(\Delta\mathbf{u})^\flat,
	\]
	from which we recover the classical Kelvin's theorem for inviscid flows.
\end{corollary}

%---------------------------------------------------------------------------------------------------

\section{Spatial discretization of the Navier-Stokes-Fourier system}\label{sec:semi-discrete}

%---------------------------------------------------------------------------------------------------

\subsection{Discrete exterior calculus based on discrete diffeomorphisms}\label{sec:DEC}

In this section we review the objects that will allow us to conveniently discretize both the material (to some extent) and spatial variational principles presented in section \ref{ssec:smooth-material} and \ref{ssec:smooth-spatial}, respectively. This amounts to discretize the infinite-dimensional Lie group $\Diff_0(M)$ into a matrix Lie group. In the case of incompressible flows, this has been done in \cite{PaMuToKaMaDe2011} whereas in the case of compressible flows, this has been done in \cite{BaGB2018}; we shall present this case below and relate it to the former. 

%---------------------------------------------------------------------------------------------------

\subsubsection{Discrete diffeomorphisms}

\begin{figure}
	\centering
	\includegraphics[scale=2.5]{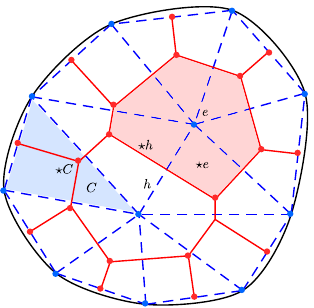}
	\caption{A simplicial mesh $\mathbb{M}$ (in blue) approximating the fluid domain $M$ (in black), as well as its dual mesh ${\star\mathbb{M}}$ (in red). Also is represented a primal cell $C$ (filled in blue), and its dual node ${\star C}$, a primal edge $h$ (in dashed blue) and its dual edge ${\star h}$ (in plain red), and a primal node $e$ and its dual cell ${\star e}$ (filled in red).}
	\label{fig:mesh}
\end{figure}

First of all, the 2D or 3D fluid domain $M$ has to be discretized into a 2D or 3D mesh $\mathbb{M}$ with $N$ cells $C_i$, $i\in\{1,\dots,N\}$, that tile the domain $M$. We do not explicitly describe what is a \emph{nice} mesh in our framework but for our purpose, the cells need to be disjoint except on shared edges (or faces in 3D). It is understood that, smaller is the characteristic mesh size $h$, better is the approximation $\mathbb{M}$ of the fluid domain $M$\footnote{Such a limit, when $h\to 0$, can be understood as taking \emph{nice} successive mesh refinements of $\mathbb{M}$.}. However, it is important to note that we do not require the mesh $\mathbb{M}$ to be structured nor regular. We will mostly illustrate the forthcoming concepts with the help of a simplicial mesh $\mathbb{M}$ of $M$ (see figure \ref{fig:mesh}). Note also that the \emph{circumcentric dual mesh}\footnote{In the case where $\mathbb{M}$ is a simplicial Delaunay triangulation, ${\star\mathbb{M}}$ is the associated Voronoi diagram.} ${\star\mathbb{M}}$ of $\mathbb{M}$, which is obtained by connecting the circumcenters of adjacent primal cells, will play a central role as we will see. We will denote by $\Omega$ the diagonal matrix with diagonal elements $\Omega_{ii}=|C_i|$ the unsigned volume (or area in 2D) of cell $C_i$.

\begin{definition}[Group of discrete diffeomorphisms]\label{def:discrete-diffeomorphisms}
	Given a mesh $\mathbb{M}$ approximating the fluid domain $M$, the matrix Lie group
	\[
		\mathsf{D}(\mathbb{M})=\big\{q\in\GL(N)^+:q\mathbf{1}=\mathbf{1}\big\}
	\]
	where $\mathbf{1}=(1,\dots,1)^\mathsf{T}\in\R^N$, is called the \emph{group of discrete diffeomorphisms} associated to $\mathbb{M}$.
\end{definition}

Thanks to this group, we will make sense of the $\Diff_0(M)$ relabeling symmetry at the discrete level and the spatial variational principle will be properly discretized in section \ref{ssec:semi-discrete-spatial}. Note that, at this point, the no-slip boundary condition has not been introduced yet, so $\mathsf{D}(\mathbb{M})$ is a discretization of $\Diff(M)$ rather than $\Diff_0(M)$.

\begin{definition}[Discrete function]\label{def:discrete-function}
	Given a mesh $\mathbb{M}$ approximating the fluid domain $M$, a \emph{discrete function} is a vector $F\in\R^N$\footnote{We will see shortly that each component $F_i$ is considered as an approximation of the cell average of an associated function $f\in L^2(M)$ on cell $C_i$, $i\in\{1,\dots,N\}$.}. The vector space of discrete functions is denoted by $\Omega^0_d(\mathbb{M})\cong\R^N$. The $L^2$ inner product of functions on $M$ is discretized into an inner product $\langle\cdot,\cdot\rangle_0$ on discrete functions by setting:
	\[
		\langle F,G\rangle_0=F^\mathsf{T}\Omega G=\sum_{i=1}^N F_i\Omega_{ii} G_i,
	\]
	for any $F$, $G\in\Omega^0_d(\mathbb{M})$.
\end{definition}

The main idea behind the definition \ref{def:discrete-diffeomorphisms} is as follows. The action of a diffeomorphism $\varphi\in\Diff(M)$ on a function $f\in L^2(M)$ is given by $f\mapsto f\circ\varphi^{-1}$. This linear action preserves constant functions. With the given definition of discrete functions and the associated inner product, we obtain that the preservation of constants at the discrete level is equivalent to the condition $q\mathbf{1}=\mathbf{1}$, for the action $F\mapsto q^{-1}F$ on the right of $\Omega^0_d(\mathbb{M})$.

\begin{remark}
	In the case of incompressible flows, we need to impose some kind of incompressibility condition. As the action of $\Diff(M)$ on $L^2(M)$ also preserves the $L^2$ inner product according to Koopman's lemma, the condition of incompressibility at the discrete level is $q^\mathsf{T}\Omega q=\Omega$, and the associated group of discrete diffeomorphisms becomes $\mathsf{D}_\text{vol}(\mathbb{M})=\big\{q\in\mathsf{D}(\mathbb{M}):q^\mathsf{T}\Omega q=\Omega\big\}$.
\end{remark}

Discrete functions can be integrated on the whole mesh $\mathbb{M}$ discretizing $M$ in the following sense: given $F\in\Omega^0_d(\mathbb{M})$, its integral over $\mathbb{M}$ is the real number $$\Tr(\Omega F^\mathsf{T})=\sum_{i=1}^N\Omega_{ii} F_i,$$ where $\Omega$ is seen as a vector of $\R^N$.

We shall now explain in which sense a discrete function $F\in\Omega^0_d(\mathbb{M})$ approximates a function $f\in \mathcal{C}^0(M)$, and in which sense a discrete diffeomorphism $q\in\mathsf{D}(\mathbb{M})$ approximates a smooth diffeomorphism $\varphi\in\Diff(M)$.

\begin{figure}
	\centering
	\includegraphics[scale=0.8]{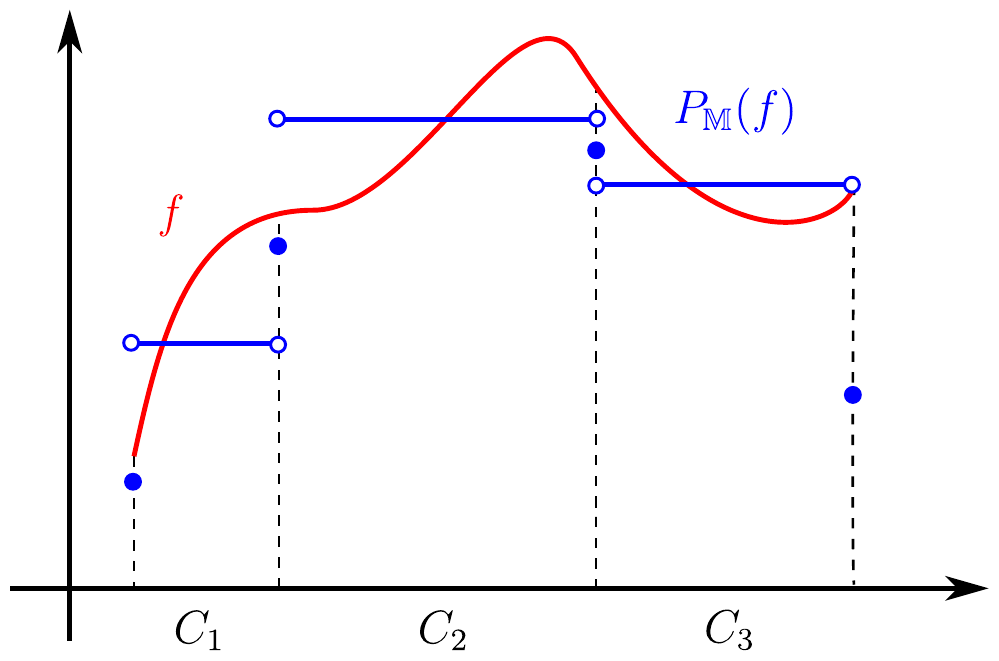}
	\caption{A continuous function $f$ (in red) and its discretization $P_\mathbb{M}(f)$ (in blue), with $\mathbb{M}$ being the 1D mesh $C_1\cup C_2\cup C_3$ (or its closure to be more precise).}
	\label{fig:reconstruction}
\end{figure}

\begin{definition}\label{def:reconstruction}
	Given a family $\{\mathbb{M}_h\}_{h>0}$ of meshes $\mathbb{M}_h$ of characteristic size $h$ approximating the fluid domain $M$, we say that a family $\{F_h\}_{h>0}$ of discrete functions $F_h\in\Omega^0_d(\mathbb{M}_h)$ \emph{approximates} the function $f\in\mathcal{C}^0(M)$ if
	\[
		\big\|S_{\mathbb{M}_h}(F_h)-f\big\|_\infty\xrightarrow[h\to 0]{}0,
	\]
	where $S_{\mathbb{M}}:\Omega^d_0(\mathbb{M})\to L^2(M)$ is the \emph{reconstruction operator} defined by $$S_\mathbb{M}(F)=\sum_{i=1}^N\left(F_i\chi_{\mathring{C}_i}+\frac{F_i}{2}\chi_{\partial C_i}\right)$$ for any mesh $\mathbb{M}$ and $F\in\Omega^0_d(\mathbb{M})$. Notice that $S_\mathbb{M}(F)$ is a piecewise continuous function (see figure \ref{fig:reconstruction}). In particular, we have the following:
	\begin{itemize}
		\item For any $x\in\mathring{C}_i$, we have
		\[
			|F_i-f(x)|=\big|S_{\mathbb{M}_h}(F)(x)-f(x)\big|\leq\sup_{x\in M}\big|S_{\mathbb{M}_h}(F)(x)-f(x)\big|=\big\|S_{\mathbb{M}_h}(F)(x)-f(x)\big\|_\infty\xrightarrow[h\to 0]{}0,
		\]
		which justifies the approximation $f(x)\approx F_i$ on $\mathring{C}_i$;
		
		\item Similarly for any $x\in C_i\cap C_j$, we have
		\[
			\left|\frac{F_i+F_j}{2}-f(x)\right|\leq\big\|S_{\mathbb{M}_h}(F)(x)-f(x)\big\|_\infty\xrightarrow[h\to 0]{}0,
		\]
		which justifies the approximation $f(x)\approx\frac{F_i+F_j}{2}$ on $C_i\cap C_j=h_{ij}$.
	\end{itemize}
\end{definition}

Given a function $f\in\mathcal{C}^0(M)$ there is a simple way to discretize it into a discrete function $F\in\Omega^0_d(\mathbb{M})$ by considering averages on cells.

\begin{definition}
	Given a mesh $\mathbb{M}$ approximating the fluid domain $M$, the \emph{discretization operator} is the map $P_\mathbb{M}:\mathcal{C}^0(M)\to\Omega^d_0(\mathbb{M})$ defined by
	\[
		P_{\mathbb{M}}(f)_i=\frac{1}{\Omega_{ii}}\int_{C_i}f\,\mathrm{d}x,
	\]
	for any $f\in\mathcal{C}^0(M)$.
\end{definition}

It can be shown that $P_{\mathbb{M}}(f)$ approximates $f$ in the sense of the previous definition. All functions that we consider hereafter will be discretized in this way. Therefore, the coefficients $F_i$ of a discrete function $F\in\Omega^0_d(\mathbb{M})$ will always be thought of as approximations of cell averages $\frac{1}{\Omega_{ii}}\int_{C_i}f\,\mathrm{d}x$, $i\in\{1,\dots,N\}$.

\begin{definition}\label{def:diffeomorphism-approximation}
	Given a family $\{\mathbb{M}_h\}_{h>0}$ of meshes $\mathbb{M}_h$ of characteristic size $h$ approximating the fluid domain $M$, we say that a family $\{q_h\}_{h>0}$ of matrices $q_h\in\mathsf{D}(\mathbb{M}_h)$ \emph{approximates} a diffeomorphism $\varphi\in\Diff(M)$ if for any $f\in\mathcal{C}^0(M)$,
	\[
		\big\|S_{\mathbb{M}_h}\big(q_h P_{\mathbb{M}_h}(f)\big)-f\circ\varphi^{-1}\big\|_\infty\xrightarrow[h\to 0]{}0.
	\]
\end{definition}

Notice that a discrete diffeomorphism $q$ naturally acts on a vector $F$ by matrix-vector multiplication $qF$, so we have a right action of $\mathsf{D}(\mathbb{M})$ on $\Omega^0_d(M)$ defined by $F\cdot q:=q^{-1}F$; whereas $\Diff(M)$ naturally acts on the right of $L^2(M)$ as $f\cdot\varphi:=f\circ\varphi$. Note that this action is unitary for $\langle\cdot,\cdot\rangle_0$ if and only if in addition $q\in\mathsf{D}_\text{vol}(\mathbb{M})$.

%---------------------------------------------------------------------------------------------------

\subsubsection{Discrete velocity fields}

We already know that in the reduction process we will need the Lie algebra of $\mathsf{D}(\mathbb{M})$.

\begin{definition}[Discrete vector field]
	Given a mesh $\mathbb{M}$ approximating the fluid domain $M$, the Lie algebra of $\mathsf{D}(\mathbb{M})$ is
	\[
		\LieAlgebraFont{d}(\mathbb{M})=\big\{A\in\gl(N):A\mathbf{1}=0\big\}
	\]
	equipped with the matrix commutator $[A,B]=AB-BA$; it is called the \emph{Lie algebra of discrete vector fields}.
\end{definition}

Actually, the physical discrete velocity fields belong to a certain \emph{linear subspace} $\mathcal{V}$ of $\LieAlgebraFont{d}(\mathbb{M})$, not a Lie subalgebra of $\LieAlgebraFont{d}(\mathbb{M})$, which, as we will see, will have important consequences on the spatial semi-discrete (discretized in space) variational principle. Once we understand what the elements of a discrete vector field $A$ approximate, we will have a clear understanding of the condition $A\mathbf{1}=0$ as a discrete divergence theorem for $A$. Note also that up to that moment we did not mention how to discretize the no-slip boundary conditions. It will be made clear at the end of this section, but right now we shall say that $\LieAlgebraFont{d}(\mathbb{M})$ is the Lie algebra of discrete vector fields that are \emph{tangent to the boundary} (that is, they satisfy a discrete version of the tangential boundary condition $\mathbf{u}\cdot\mathbf{n}=0$).

\begin{remark}\label{rq:discrete-incompressibility}
	In the case of incompressible flows, the Lie algebra of $\mathsf{D}(\mathbb{M})$ should reflect the fact that the discrete vector fields are divergence-free, for a discrete notion of divergence. Computing the Lie algebra $\LieAlgebraFont{d}_\text{vol}(\mathbb{M})$ of $\mathsf{D}_\text{vol}(\mathbb{M})$, we find $\LieAlgebraFont{d}_\text{vol}(\mathbb{M})=\big\{A\in\LieAlgebraFont{d}(\mathbb{M}):A^\mathsf{T}{\Omega}+\Omega A=0\big\}$, so the divergence-free condition corresponds to the discrete condition $A^\mathsf{T}{\Omega}+\Omega A=0$ and we will see later how it is related to a notion of discrete divergence.
\end{remark}

Let $\{q_h(t)\}_{h>0}$, $q_h(t)\in\mathsf{D}(\mathbb{M})$, be a family of discrete motions approximating a smooth motion $\varphi_t\in\Diff(M)$ in the sense of definition \ref{def:diffeomorphism-approximation}. Then, for $h>0$ fixed and for $f\in\mathcal{C}^0(M)$, set $F=P_{\mathbb{M}_h}(f)$. The definition means that at all time, $q_h(t)F:=F_t$ converges to the function $f\circ\varphi_t^{-1}:=f_t$. Define $A_h(t)=\dot{q}_h(t)q_h^{-1}(t)$, which is the discrete analogue of the velocity vector field $\mathbf{u}_t=\dot{\varphi}_t\circ\varphi_t^{-1}$. Therefore, in the discrete case we have $\dot{F}_t=A_h(t)F_t$, whereas in the continuous case we have $\dot{f}_t=-\lie_{\mathbf{u}_t}f_t$, which suggests the following definition.

\begin{definition}\label{def:vector-field-approximation}
	Given a family $\{\mathbb{M}_h\}_{h>0}$ of meshes $\mathbb{M}_h$ of characteristic size $h$ approximating the fluid domain $M$, we say that a family $\{A_h\}_{h>0}$ of discrete vector fields $A_h\in\LieAlgebraFont{d}(\mathbb{M}_h)$ approximates the smooth vector field $\mathbf{u}\in\LieAlgebraFont{X}(M)$ if for any $f\in\mathcal{C}^0(M)$,
	\[
		\big\|S_{\mathbb{M}_h}\big(A_h P_{\mathbb{M}_h}(f)\big)-(-\lie_\mathbf{u}f)\big\|_\infty\xrightarrow[h\to 0]{}0.
	\]
\end{definition}

We have the following result, see \citet[lemma 2]{PaMuToKaMaDe2011} for the proof.

\begin{proposition}
	Given a family $\{\mathbb{M}_h\}_{h>0}$ of meshes $\mathbb{M}_h$ of characteristic size $h$ approximating the fluid domain $M$, let $\{A_h(t)\}_{h>0}$ be a family of time-dependent discrete vector fields $A_h(t)\in\LieAlgebraFont{d}(\mathbb{M}_h)$ that approximates the time-dependent smooth vector field $\mathbf{u}_t\in\LieAlgebraFont{X}(M)$. Then there exists a family $\{q_h(t)\}_{h>0}$ of discrete motions $q_h(t)\in\mathsf{D}(\mathbb{M}_h)$ such that $A_h(t)=\dot{q}_h(t)q^{-1}_h(t)$ and $\{q_h(t)\}_{h>0}$ approximates the motion $\varphi_t\in\Diff(M)$ associated to $\mathbf{u}_t$, in the sense of definition \ref{def:diffeomorphism-approximation}.
\end{proposition}

\begin{figure}
	\centering
	\includegraphics[scale=2.5]{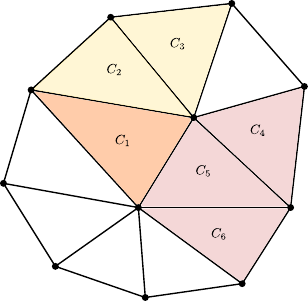}
	\caption{Here the cells $C_2$ and $C_5$ are two cells away from each other. However, given $A$ and $B\in\mathcal{S}$, we have that $[A,B]_{25}\neq 0$ as $[A,B]_{25}=\sum_{i\in N(2)\cap N(5)}(A_{2i}B_{i5}-B_{2i}A_{i5})=A_{21}B_{15}-B_{21}A_{15}\neq 0$ since $C_1$ is the only common cell to $C_2$ and $C_5$. Consequently $[A,B]\notin\mathcal{S}$.}
	\label{fig:sparsity}
\end{figure}

Given a discrete motion $t\mapsto q(t)\in\mathsf{D}(\mathbb{M})$, the coefficient $A_{ij}(t)$ of the associated discrete vector field $A(t)=\dot{q}(t)q(t)^{-1}\in\LieAlgebraFont{d}(\mathbb{M})$ can be thought of how much of the fluid is transported between cells $C_i$ and $C_j$. Consequently it is physically sound to have $A_{ij}(t)$ non-zero only if cells $C_i$ and $C_j$ share an edge in 2D, or a face in 3D; so we introduce
\[
	\mathcal{S}=\big\{A\in\LieAlgebraFont{d}(\mathbb{M}):A_{ij}=0,\quad\forall j\notin N(i)\big\},
\]
where $N(i)$ denotes the set of all indices, \emph{including} $i$, of cells sharing an edge in 2D or a face in 3D with the cell $C_i$. Note that $\mathcal{S}$ is only a linear subspace of $\LieAlgebraFont{d}(\mathbb{M})$, not a Lie subalgebra, because of the matrix multiplications appearing in the matrix commutator\footnote{The coefficient $[A,B]_{ij}$ may be non zero if cells $C_i$ and $C_j$ are two cells away from each other, see figure \ref{fig:sparsity}.}. As a constraint on discrete vector fields, $\mathcal{S}$ is therefore \emph{nonholonomic} (we have $[\mathcal{S},\mathcal{S}]\subset\LieAlgebraFont{d}(\mathbb{M})$ and $\mathcal{S}\cap[\mathcal{S},\mathcal{S}]=\{0\}$), we will call it the \emph{sparsity constraint}. This requirement on discrete vector fields let us describe the coefficients $A_{ij}$, $j\in N(i)$, of a discrete vector field $A\in\mathcal{S}$.

\begin{proposition}\label{pro:vector-field-approximation-coefficients}
	Let $A\in\mathcal{S}$ approximate a smooth vector field $\mathbf{u}\in\LieAlgebraFont{X}(M)$, tangent to the boundary $\partial M$, in the sense of definition \ref{def:vector-field-approximation}. Then:
	\begin{align*}
		A_{ij}&\approx-\frac{1}{2\Omega_{ii}}\int_{h_{ij}}\mathbf{u}\cdot\mathrm{d}\mathbf{S},\quad\forall j\in N(i),\ j\neq i,\\
		A_{ii}&\approx\frac{1}{2\Omega_{ii}}\int_{C_i}\dive\mathbf{u}\,\mathrm{d}x,
	\end{align*}
	where $h_{ij}$ denotes the shared edge in 2D, or face in 3D, between the cells $C_i$ and $C_j$, and oriented by the exterior normal to $C_i$.
\end{proposition}

\begin{proof}
	We only give a sketch of the proof; a more detailed one is currently reviewed, in a more analytical context. Given any discrete function $F\in\Omega_0(\mathbb{M})$ we have for any cell $C_i$ that
	\[
		-\Omega_{ii}(AF)_i\approx\int_{C_i}\dif f(\mathbf{u})\,\mathrm{d}x=\underbrace{\sum_{j\in N(i)}\int_{h_{ij}}f\mathbf{u}\cdot\mathbf{n}\,\mathrm{d}S}_{\circled{1}}-\underbrace{\int_{C_i}f\dive\mathbf{u}\,\mathrm{d}x}_{\circled{2}}.
	\]
	In \circled{1} we use that $f\approx\frac{F_i+F_j}{2}$ on $h_{ij}$ as explained before. :
	Therefore \circled{1} can be approximated by
	\[
		\frac{F_i}{2}\int_{C_i}\dive\mathbf{u}\,\mathrm{d}x+\sum_{j\in N(i)}\frac{F_j}{2}\int_{h_{ij}}\mathbf{u}\cdot\mathbf{n}\,\mathrm{d}S,
	\]
	whereas \circled{2} can be approximated by
	\[
		F_i\int_{C_i}\dive\mathbf{u}\,\mathrm{d}x,
	\]
	since $f\approx F_i$ on $\mathring{C}_i$. Consequently, $-\Omega_{ii}(AF)_i$ can be seen as an approximation of
	\[
		-\frac{F_i}{2}\int_{C_i}\dive\mathbf{u}\,\mathrm{d}x+\sum_{j\in N(i)}\frac{F_j}{2}\int_{h_{ij}}\mathbf{u}\cdot\mathbf{n}\,\mathrm{d}S.
	\]
	Since $-\Omega_{ii}(AF)_i=-\Omega_{ii}A_{ii}F_i-\sum_{j\in N(i)}\Omega_{ii}A_{ij}F_j$ and $F$ is arbitrary, we obtain the aforementioned approximations.
\end{proof}

\begin{figure}
	\centering
	\includegraphics[scale=2.5]{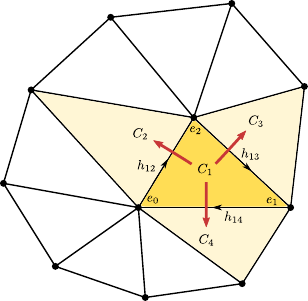}
	\caption{Fluid is coming from neighboring cells $C_2$, $C_3$ and $C_4$ out or in the cell $C_1=[e_0,e_1,e_2]$ through its boundary $\partial C_1=-[e_0,e_2]+[e_1,e_2]+[e_0,e_1]=h_{12}-h_{13}-h_{14}$ (assuming that the primal edges $h_{ij}$ are positively oriented with respect to the exterior normal), which has for effect to compress or dilate the fluid in $C_1$.}
	\label{fig:discrete-divergence-theorem}
\end{figure}

Note that, thanks to this proposition, we can understand the defining condition $A\mathbf{1}=0$ for $A\in\mathcal{S}$ as a discrete version of the classical divergence theorem, see figure \ref{fig:discrete-divergence-theorem}.

From the previous proposition, we also observe that a discrete vector field $A\in\LieAlgebraFont{d}(\mathbb{M})$ shall satisfies the relation $\Omega_{ii}A_{ij}=-\Omega_{jj}A_{ji}$, for $j\in N(i)$ with $j\neq i$, because of the orientation of the shared edge or face $h_{ij}$. This entails that the matrix $A^\mathsf{T}\Omega+\Omega A$ is a diagonal one. Therefore, besides the sparsity constraint $\mathcal{S}$, we define the \emph{compressibility constraint} to be
\begin{align*}
	\mathcal{V}&=\big\{A\in\LieAlgebraFont{d}(\mathbb{M}):A^\mathsf{T}\Omega+\Omega A\text{ is diagonal}\big\}\\
	&=\big\{A\in\LieAlgebraFont{d}(\mathbb{M}):A=A^\text{vol}+A^\text{diag},\ A^\text{vol}\in\LieAlgebraFont{d}_\text{vol}(\mathbb{M})\text{ and }A^\text{diag}\text{ diagonal}\big\},
\end{align*}
with $A^\text{diag}$ satisfying $A^\mathsf{T}\Omega+\Omega A=2\Omega A^\text{diag}$. Alternatively in terms of coefficients, a discrete vector field $A\in\LieAlgebraFont{d}(\mathbb{M})$ satisfies the compressibility constraint if and only if 
\begin{equation}\label{eq:compressibility-constraint}
	\Omega_{ii}A_{ij}+A_{ji}\Omega_{jj}=2\Omega_{ii}A_{ii}\delta_{ij}
\end{equation}
for all $i$, $j\in\{1,\dots,N\}$. Note that, again, $\mathcal{V}$ is nonholonomic. $\mathcal{V}$ is called the compressibility constraint as the \emph{physical} discrete velocity field is required to satisfy this constraint.

\begin{definition}[Discrete divergence]
	Given a mesh $\mathbb{M}$ approximating the fluid domain $M$, the \emph{discrete divergence} of a discrete vector field $A\in\LieAlgebraFont{d}(\mathbb{M})$ is the discrete function $\dive A\in\Omega^0_d(\mathbb{M})$ defined by $(\dive A)_i=2 A_{ii}$, for all $i\in\{1,\dots,N\}$.
\end{definition}

Therefore the defining condition for $\mathcal{V}$ reads $\Omega_{ii}A_{ij}+A_{ji}\Omega_{jj}=\Omega_{ii}\delta_{ij}(\dive A)_i$ for for all $i$, $j\in\{1,\dots,N\}$, and we now see the relationship with the incompressibility condition $A^\mathsf{T}\Omega+\Omega A=0$ in remark \ref{rq:discrete-incompressibility}: $A\in\LieAlgebraFont{d}(\mathbb{M})$ is incompressible if and only if $\dive A=0$. Note also that this definition is coherent with proposition \ref{pro:vector-field-approximation-coefficients} in the case where in addition $A\in\mathcal{S}$, because $(\dive A)_i$ should approximate $\frac{1}{\Omega_{ii}}\int_{C_i}\dive\mathbf{u}\,\mathrm{d}x$, and we know from the proposition that $2 A_{ii}$ already approximates this cell average.

\begin{lemma}\label{lem:integral-discrete-divergence}
	Given a mesh $\mathbb{M}$ approximating the fluid domain $M$ and $A\in\mathcal{V}$, the integral of $\dive A$ over $\mathbb{M}$ is zero, that is, we have: $$\sum_{i=1}^N\Omega_{ii}(\dive A)_i=0.$$
\end{lemma}

Upon adding the sparsity constraint into the statement, one sees that this lemma is precisely the discrete counterpart of the fact that the integral of $\dive\mathbf{u}$ over the whole domain $M$ is zero for $\mathbf{u}$ any vector field tangent to the boundary $\partial M$.

We now introduce the action $F\cdot A$ of a discrete vector field $A\in\LieAlgebraFont{d}(\mathbb{M})$ on a discrete function $F\in\Omega^0_d(\mathbb{M})$. This is obtained as the infinitesimal action of $\mathsf{D}(\mathbb{M})$ on $\Omega^0_d(\mathbb{M})$: $F\cdot A=-AF$. It is the discrete analogue of the action $\lie_\mathbf{u}f=\dif f(\mathbf{u})$ of a smooth vector field $\mathbf{u}\in\LieAlgebraFont{X}(M)$ on a smooth function $f\in\smooth{M}$. When $A\in\mathcal{S}$ we compute for any $i\in\{1,\dots,N\}$ that:
\[
	(F\cdot A)_i=-\sum_{\substack{j\in N(i)}}A_{ij}F_j.
\]
Note that this is a Lie algebra action, meaning here that the map $\LieAlgebraFont{d}(\mathbb{M})\to\operatorname{End}\Omega_d^0(\mathbb{M})$, $A\mapsto\rho_A$, with $\rho_A(F):=F\cdot A=-AF$, is a Lie algebra anti-homomorphism.

We are going now to extend what has been done until now to include discrete vector fields that are not tangential to the boundary. This is an essential feature of the framework developed in this paper, as the smooth vector field $\mathbf{j}_S$ is potentially not tangent to the boundary $\partial M$ in the case the fluid is not adiabatic (see section \ref{ssec:smooth-spatial}). We will also make clear how boundary conditions are handled within this framework. First of all, we define $\overline{\mathbb{M}}=\mathbb{M}\cup C_\partial$, where $C_\partial$ is a special cell modeling the \emph{environment} outside of the fluid domain $M$ which has been discretized into a mesh $\mathbb{M}$. We will also denote by $\mathbb{M}^\circ$ the mesh composed of all cells $C_i\subset\mathbb{M}$ which do not have an edge in $\partial\mathbb{M}$. Then the Lie algebra of (general) discrete vector fields is defined by:
\[
	\LieAlgebraFont{d}(\overline{\mathbb{M}})=\big\{A\in\gl(N+1):A\cdot\mathbf{1}=0\big\},
\]
where $\mathbf{1}=(1,\dots,1)^\mathsf{T}\in\R^{N+1}$. The spatial discretization $J_S$ of the smooth vector field $\mathbf{j}_S$ will belong to this Lie algebra. Both sparsity and compressibility constraints are defined as before, except now that cells in $\mathbb{M}$ can have a new neighbor: $N(\partial)=\big\{i\in\{1,\dots,N\}:C_i\in\mathbb{M}{\setminus}\mathbb{M}^\circ\big\}\cup\{\partial\}$; we will denote these linear subspaces of $\LieAlgebraFont{d}(\overline{\mathbb{M}})$ by $\overline{\mathcal{S}}$ and $\overline{\mathcal{V}}$, respectively. By convention, the index $\partial$ of the cell $C_\partial$ is the $(N+1)^\text{th}$ cell of $\overline{\mathbb{M}}$. The proposition \ref{pro:vector-field-approximation-coefficients} is extended as follows:
\[
	A_{i\partial}\approx-\frac{1}{2\Omega_{ii}}\int_{h_{i\partial}}\mathbf{u}\cdot\mathrm{d}\mathbf{S},
	\quad A_{\partial i}\approx-\frac{1}{2\Omega_{\partial\partial}}\int_{h_{\partial i}}\mathbf{u}\cdot\mathrm{d}\mathbf{S},
	\quad A_{\partial\partial}\approx\frac{1}{2\Omega_{\partial\partial}}\int_{C_\partial}\dive\mathbf{u}\,\mathrm{d}x;
\]
where $\mathbf{u}$ is a smooth vector field not necessarily tangent to the boundary $\partial M$, $h_{i\partial}$ is the edge common to $C_i$ and the exterior environment, oriented from $C_i$ to $C_\partial$, and provided that we can make sense of the volume $\Omega_{\partial\partial}$ of $C_\partial$. To that end, notice that:
\[
	A_{\partial\partial}=-\sum_{\substack{i\in N(\partial)\\ i\neq\partial}}A_{\partial i}\approx\sum_{\substack{i\in N(\partial)\\ i\neq\partial}}\frac{1}{2\Omega_{\partial\partial}}\int_{h_{\partial i}}\mathbf{u}\cdot\mathrm{d}\mathbf{S}=\frac{1}{2\Omega_{\partial\partial}}\int_{\partial\mathbb{M}}\mathbf{u}\cdot\mathrm{d}\mathbf{S}=\frac{1}{2\Omega_{\partial\partial}}\int_{\mathbb{M}}\dive\mathbf{u}\,\mathrm{d}x,
\]
where $\mathbb{M}$ is interpreted as a piecewise smooth manifold with boundary. Therefore, in accordance with proposition \ref{pro:vector-field-approximation-coefficients}, we set $C_\partial=\mathbb{M}$ and $\Omega_{\partial\partial}=\sum_{i=1}^N\Omega_{ii}$. Note that the left hand side of this approximation means that $A_{\partial\partial}$ is the discretization of the net flux through the boundary of the fluid domain of the vector field approximated by $A$. Moreover, notice that the right hand side of this approximation is given by
\[
	\sum_{i=1}^N\sum_{\substack{j\in N(i)\\ j\neq i}}\frac{1}{2|\mathbb{M}|}\int_{h_{ij}}\mathbf{u}\cdot\mathrm{d}\mathbf{S},
\]
therefore we can define $A_{\partial\partial}$ by setting
\[
	A_{\partial\partial}=\sum_{i=1}^N\sum_{\substack{j\in N(i)\\ j\neq i}}\frac{\Omega_{ii}}{|\mathbb{M}|}A_{ij}.
\]
In this context, tangential boundary conditions are implemented into the following Lie subalgebra:
\[
	\LieAlgebraFont{d}(\mathbb{M})=\big\{A\in\LieAlgebraFont{d}(\overline{\mathbb{M}}):A_{i\partial}=0\ \forall C_i\in\overline{\mathbb{M}}\big\}
	\simeq\big\{A\in\gl(N):A\cdot\mathbf{1}=0\big\}.
\]
In that case, we recover the usual sparsity and compressibility constraints. Finally, we also have the Lie subalgebra of discrete vector fields satisfying the no-slip boundary condition, defined as:
\[
	\mathfrak{d}_0(\mathbb{M})=\big\{A\in\LieAlgebraFont{d}(\overline{\mathbb{M}}):A_{ij}=0\ \forall C_i\in\overline{\mathbb{M}},\ \forall C_j\in \overline{\mathbb{M}}{\setminus}\mathbb{M}^\circ\}\simeq\big\{A\in \mathfrak{gl}(N-N^\partial):A\cdot\mathbf{1}=0\big\},
\]
where $N^\partial$ is the number of cells in $\mathbb{M}$ that have an edge in $\partial\mathbb{M}$. Practically this means that a cell that has an edge in $\partial\mathbb{M}$ does not transfer anything to the exterior environment nor to its neighboring cells, and we take this definition as the discretization of the no-slip boundary condition $\mathbf{u}\cdot\mathbf{n}=0$ where $\mathbf{u}\in\LieAlgebraFont{X}(M)$. Note that we can also define the Lie groups $\mathsf{D}(\overline{\mathbb{M}})$ and $\mathsf{D}_0(\mathbb{M})$ of the Lie algebras $\LieAlgebraFont{d}(\overline{\mathbb{M}})$ and $\mathfrak{d}_0(\mathbb{M})$ respectively, but we will not use them in what follows.

We define as well the space of extended discrete functions $\Omega_d^0(\overline{\mathbb{M}})=\Omega_0^d(\mathbb{M})\oplus\R$. That is, an extended discrete function $F\in\Omega_d^0(\overline{\mathbb{M}})$ is a usual discrete function with an additional component $F_\partial$, which is used to store the value of the discrete function on the exterior environment of the fluid domain (this value can also be included in the reconstruction operator $S_\mathbb{M}$ defined in \ref{def:reconstruction}). The inner product $\langle\cdot,\cdot\rangle_0$, as well as discrete integration, remain however the same on these extended discrete functions, as the exterior environment does not play any role, only the mesh $\mathbb{M}$; and thus summations are performed from $1$ to $N$. Given $A\in\LieAlgebraFont{d}(\mathbb{M})$, its divergence $\dive A\in\Omega_d^0(\overline{\mathbb{M}})$ is defined as usual by $(\dive A)_i=2A_{ii}$, for $i\in\{1,\dots,N+1\}$, and the action of $A$ on $F\in\Omega_d^0(\overline{\mathbb{M}})$ is also defined as usual by $F\cdot A=-AF\in\Omega_d^0(\overline{\mathbb{M}})$. We finish this section by adapting lemma \ref{lem:integral-discrete-divergence} to the case where the discrete vector field is not tangent to the boundary.

\begin{lemma}\label{lem:integral-discrete-divergence-general}
	Given a mesh $\mathbb{M}$ approximating the fluid domain $M$ and $A\in\overline{\mathcal{V}}$ we have:
	\[
		\sum_{i=1}^N\Omega_{ii}(\dive A)_i=\sum_{i=1}^N\Omega_{ii}A_i^\partial,
	\]
	where $A^\partial\in\Omega_d^0(\overline{\mathbb{M}})$ is defined by $A^\partial_i=-2A_{i\partial}$, $i\in\{1,\dots,N+1\}$.
\end{lemma}

\begin{proof}
	The proof goes as in the one of lemma \ref{lem:integral-discrete-divergence}:
	\[
		\sum_{i=1}^N\Omega_{ii}(\dive A)_i=2\sum_{1\leq i,j\leq N}\Omega_{ii}A_{ij}=-2\sum_{i=1}^N\Omega_{ii}A_{i\partial}.
	\]
\end{proof}

Notice that this identity is exactly the discrete analogue of $\int_M\dive\mathbf{u}\,\mathrm{d}x=\int_{\partial M}\mathbf{u}\cdot\mathrm{d}\mathbf{S}$ for $\mathbf{u}\in\LieAlgebraFont{X}(M)$.

\begin{remark}
	Suppose that the fluid domain $M$ is not simply connected, for instance suppose that there is a hole in $M$. Then any mesh $\mathbb{M}$ discretizing $M$ will also have a hole, and can be considered as a part of the exterior environment of the fluid. Then for each connected component of the exterior environment can be modeled by adding a special cell to the mesh.
\end{remark}

We finish this section by saying a few words on the \emph{initialization} and \emph{reconstruction} processes that we will use later on, that is, how we obtain a discrete vector field $A$ from a smooth vector field $\mathbf{u}$, and how we obtain a smooth vector field $\mathbf{u}$ from a discrete vector field $A$, for the case of a \emph{simplicial} mesh $\mathbb{M}$. Given a cell $C_i$ and a smooth vector field $\mathbf{u}$ on the fluid domain $M$, we denote by $\mathbf{u}_{ij}$ the normal component of $\mathbf{u}$ with respect to the exterior normal $\mathbf{n}_{ij}$ from $C_i$ to an adjacent cell $C_j$, that is, $\mathbf{u}_{ij}=\mathbf{u}(\overline{x}_{ij})\cdot\mathbf{n}_{ij}$, with $\overline{x}_{ij}$ the midpoint of the shared edge $h_{ij}$. These values used to approximate the fluxes given by proposition \ref{pro:vector-field-approximation-coefficients}:
\[
	A_{ij}=-\frac{|h_{ij}|\mathbf{u}_{ij}}{2\Omega_{ii}},\quad j\in N(i),\ j\neq i,
\]
and then we set $A_{ii}=-\sum_{j\in N(i),\ j\neq i}A_{ij}$. Thus we obtain discrete vector field $A$ from a smooth vector field $\mathbf{u}$. Of course, other flux discretizations are possible for this initialization procedure. 

Concerning the reconstruction process, we use Raviart-Thomas finite elements of the lowest degree. Given a pair of adjacent cells $C_i$ and $C_j$ we define a smooth vector field $\bm{\Psi}_{ij}$ by
\[
	\bm{\Psi}_{ij}(x)=\frac{|h_{ij}|(x-x_{ij})}{2\Omega_{ii}},
\]
with $x_{ij}$ being the (coordinates of the) primal node opposite to the primal edge $h_{ij}$. Then we obtain a smooth vector field ${\mathbf{u}}_i$ for the whole cell $C_i$ be setting ${\mathbf{u}}_i(x)=\sum_{j=1}^{3}\mathbf{u}_{ij}\bm{\Psi}_{ij}(x)$. Finally the reconstructed vector field for the whole mesh is $\mathbf{u}=\sum_{i=1}^N\mathbf{u}_i$. Other reconstruction methods are possibly feasible too.

%---------------------------------------------------------------------------------------------------

\subsubsection{Discrete differential forms}

The subject of discrete differential forms has a long history and appear in several frameworks, notably in the discrete exterior calculus (DEC) of \cite{DeHiLeMa2015}. Here we follow the approach developed in \cite{PaMuToKaMaDe2011}, which is rich in relationships with \cite{DeHiLeMa2015} but is well adapted to Euler-Poincaré reduction thanks to its interplay with discrete vector fields as defined in the previous section\footnote{It is interesting to compare the approach of \cite{PaMuToKaMaDe2011} with the DEC of \cite{DeHiLeMa2015} because recently a systematic study of the convergence of DEC to smooth exterior calculus has begun, see \cite{ScTs2018}.}.

In terms of DEC, our discrete differential forms are a bit more general than cochains on the dual mesh. Indeed, in the traditional chain-cochain approach, on either the primal or dual mesh, chains are defined relatively to a \emph{simplicial complex}, meaning that each face of a simplex belongs to the complex itself; this requirement is relaxed here. Let $\mathbb{M}$ be a mesh approximating the fluid domain $M$, whose cells are denoted by $C_i$ and their circumcenter by ${\star C_i}$, $i\in\{1,\dots,N\}$. We denote by $\mathsf{C}_k({\star\mathbb{M}})$ the free vector space generated by the oriented $k$-simplices of ${\star\mathbb{M}}$. For instance, a $0$-simplex is just a circumcenter, a $1$-simplex is the segment joining two circumcenters (not necessarily circumcenters of adjacent primal cells), and a $2$-simplex is a triangle joining three circumcenters; of which we can take linear combinations that live in $\mathsf{C}_k({\star\mathbb{M}})$.

\begin{definition}
	Given a mesh $\mathbb{M}$ approximating the fluid domain $M$ with $N$ cells, a \emph{discrete differential form} of degree $k$ (or \emph{discrete $k$-form}) is an element of the linear dual $\mathsf{C}_k(\star\mathbb{M})^*$, or, alternatively\footnote{The linear isomorphism is given by $\Psi:\Lambda^{k+1}(\R^N)^*\to\mathsf{C}_k({\star\mathbb{M}})^*$ defined by $\big\langle\Psi(\omega),[\star C_{i_1},\dots,\star C_{i_k}]\big\rangle=\omega_{i_1\dots i_k}$, where $[C_{i_0},\dots,C_{i_k}]$ denotes the $k$-simplex based on the circumcenters $\star C_{i_0},\dots,\star C_{i_k}$. In practice we will not make any distinction.}, an alternating $(k+1)$-linear map on $\R^N$, that is, an element of $\Lambda^{k+1}(\R^N)^*$. The vector space of discrete $k$-forms is denoted by $\Omega^k_d(\mathbb{M})$.
\end{definition}

Concretely, a discrete $0$-form is a vector $(F_{i})_{1\leq i\leq N}$, a discrete $1$-form is a skew-symmetric matrix $(F_{ij})_{1\leq i,j\leq N}$, and a discrete $2$-form is a completely skew-symmetric 3-tensor $(F_{ijk})_{1\leq i,j,k\leq N}$. Note that in the case of a discrete 2-form or 3-form, the indices in the coefficients do not necessarily correspond to adjacent cells; if it the case, they can be interpreted as discrete differential forms in the sense of DEC. For a discrete function $F$, $F_i$ is an approximation of a function $f$ at the circumcenter ${\star C_i}$ of the cell $C_i$. For a discrete 1-form $F$, $F_{ij}$ is an approximation of the integral of a differential 1-form on the $1$-simplex $[{\star C_i},{\star C_j}]$ (akin to a circulation). For a discrete 2-form $F$, $F_{ijk}$ is an approximation of the integral of a differential 2-form on the 2-simplex $[{\star C_i},{\star C_j},{\star C_k}]$ (akin to a flux).

We introduce now operators that are the discrete counterparts of operators in the smooth setting. 

\begin{definition}[Discrete exterior derivative]
	Let $\mathbb{M}$ be a mesh approximating the fluid domain $M$. The \emph{discrete exterior operator} is the map $\dif:\Omega_d^k(\mathbb{M})\to\Omega_d^{k+1}(\mathbb{M})$ defined for any $F\in\Omega_d^k(\mathbb{M})$ by
	\[
		(\dif F)_{i_0\dots i_k}=\sum_{0\leq j\leq k}(-1)^j F_{i_0\dots\widehat{i_j}\dots i_k},
	\]
	where the hat denotes an omission. In particular, for $F\in\Omega^0_d(\mathbb{M})$, $(\dif F)_{ij}=F_j-F_i$, and for $F\in\Omega^1_d(\mathbb{M})$, $(\dif F)_{ijk}=F_{ij}+F_{jk}+F_{ki}$. The differential of an extended discrete function $F\in\Omega_0^d(\overline{\mathbb{M}})$ is defined in a similar manner, except it yields a $(N+1)\times(N+1)$ matrix instead of a $N\times N$ matrix.
\end{definition}

We now summarize the construction of the most important operator of Pavlov's framework: the flat operator $\flat:\mathcal{S}\to\Omega^1_d(\mathbb{M})$; see \cite[Section 3.5]{PaMuToKaMaDe2011} for details. We would like to construct a flat operator $\flat:\mathcal{S}\to\Omega^1_d(\mathbb{M})$ and a pairing $\llangle\cdot,\cdot\rrangle:\Omega^1_d(\mathbb{M})\times\LieAlgebraFont{d}(\mathbb{M})\to\R$ which approximates (as the characteristic size of the mesh decreases) the smooth $L^2$ inner product of vector fields $\int_M\mathbf{u}\cdot\mathbf{v}\,\mathrm{d}x$, with $\mathbf{u}$, $\mathbf{v}\in\LieAlgebraFont{X}(M)$, and in which the Euclidean metric is used. Such a pairing will be used in forthcoming variational principles, and $A^\flat$ will be paired not only with discrete vector fields $B\in\mathcal{S}$ but also with discrete vector fields of the form $[B,C]\in[\mathcal{S},\mathcal{S}]$\footnote{We will have variations of the form $\delta A=\dot{B}-[A,B]$ for the discrete velocity field $A$.}. Since the coefficient $[B,C]_{ij}$ may not be zero when $C_i$ and $C_j$ are two cells away of each other, the coefficient $A^\flat_{ij}$ also may not be zero for such cells, which justifies the generalization which was raised at the beginning of this section. The pairing itself is the $\Omega$-weighted Frobenius inner product inherited from $\gl(N)$: $$\llangle L,B\rrangle:=\Tr(L^\mathsf{T}\Omega B)=\sum_{1\leq i,j\leq N}\Omega_{ii} L_{ij}B_{ij}.$$

\begin{definition}[Discrete flat operator]
	Let $\{\mathbb{M}_h\}_{h>0}$ be a family of meshes $\mathbb{M}_h$ approximating the fluid domain $M$. An operator $\flat_h:\mathcal{S}\to\Omega^1_d(\mathbb{M}_h)$ is called a \emph{discrete flat operator} if the following conditions are satisfied:
	\begin{enumerate}[label=(\arabic*), ref=\thetheorem.(\arabic*)]
		
		\item For any $A_h$, $B_h\in\mathcal{S}$ which approximate $\mathbf{u}$, $\mathbf{v}\in\LieAlgebraFont{X}(M)$ in the sense of definition \ref{def:vector-field-approximation},
		\[
			\llangle[\big] A_h^{\flat_h},B_h\rrangle[\big]\xrightarrow[h\to 0]{}\llangle\mathbf{u}^\flat,\mathbf{v}\rrangle=\int_M\mathbf{u}\cdot\mathbf{v}\,\mathrm{d}x,
		\]
		
		\item For any $A_h$, $B_h\in\mathcal{S}$ and $C_h$ which approximate $\mathbf{u}$, $\mathbf{v}$ and $\mathbf{w}\in\LieAlgebraFont{X}(M)$ in the sense of definition \ref{def:vector-field-approximation},
		\[
			\llangle[\big] A_h^{\flat_h},[B_h,C_h]\rrangle[\big]\xrightarrow[h\to 0]{}\llangle[\big]\mathbf{u}^\flat,[\mathbf{v},\mathbf{w}]\rrangle[\big]=\int_M\mathbf{u}\cdot[\mathbf{v},\mathbf{w}]\,\mathrm{d}x.
		\]
		
	\end{enumerate}
\end{definition}

This last condition can be seen to be equivalent to the convergence of the \emph{partial discrete vorticity} $\dif A_h^{\flat_h}$ to the smooth vorticity $\dif\mathbf{u}^\flat$ when $h\to 0$, in a certain sense, see \citet[lemma 4]{PaMuToKaMaDe2011}.

Let's fix a mesh $\mathbb{M}$ approximating $M$, and let $A\in\mathcal{S}$. For $j\in N(i)$, we define:
\begin{equation}\label{eq:flat-operator}
	A_{ij}^\flat:=2\Omega_{ii}A_{ij}\frac{|{\star h}_{ij}|}{|h_{ij}|},
\end{equation}
where $h_{ij}$ is the shared edge (in 2D, or face in 3D) between cells $C_i$ and $C_j$, and ${\star h}_{ij}$ is the dual edge between these cells, that is, the edge in the dual mesh ${\star\mathbb{M}}$ that connects the circumcenters of $C_i$ and $C_j$. Indeed, recall from \ref{pro:vector-field-approximation-coefficients} that $2\Omega_{ii}A_{ij}$ is an approximation of the flux $-\int_{h_{ij}}\mathbf{u}\cdot\mathrm{d}\mathbf{S}$, for instance suppose that $2\Omega_{ii}A_{ij}=-\mathbf{u}(x_{ij})\cdot\mathbf{n}\,|h_{ij}|$ with $x_{ij}$ the barycenter of $h_{ij}$. Then, with this approximation, $A_{ij}^\flat=-\mathbf{u}(x_{ij})\cdot\mathbf{n}\,|{\star h}_{ij}|$, which is an approximation of the circulation $-\int_{{\star h}_{ij}}\mathbf{u}\cdot\mathrm{d}\mathbf{l}$, and therefore $A^\flat$ is a discrete 1-form as defined before. Then, it can be shown that this definition for $j\in N(i)$ make $\flat$ satisfy the first condition, see \citet[section 3.5]{PaMuToKaMaDe2011}. We also need to define $A_{jk}^\flat$ in the case where the cells $C_j$ and $C_k$ are two cells away from each other; let's call $C_i$ the cell they share. Then $A_{jk}^\flat$ is defined by the relation
\begin{equation}\label{eq:two-vorticities} 
	(\dif A^\flat)_{ijk}=A_{ij}^\flat+A_{jk}^\flat+A_{ki}^\flat=K_{ijk}\,\omega_A(e_{ijk}),
\end{equation}
where the $K_{ijk}$ are constants to be determined and $\omega_A(e_{ijk})$ is the \emph{total discrete vorticity} around the primal node (in 2D, or edge in 3D)\footnote{This way the discrete vorticity could be seen as a dual discrete function.} shared by the cells $C_i$, $C_j$ and $C_k$. This notion of total discrete vorticity is defined at any primal node $e\in\mathbb{M}$ by:
\begin{equation}\label{eq:total-vorticity}
	\omega_A(e)=\sum_{{\star h}_{ij}\in\partial({\star e})}s_{ij}A_{ij}^\flat,
\end{equation}
where ${\star e}$ is the dual face to $e$ in ${\star\mathbb{M}}$, and $s_{ij}=+1$ if ${\star h}_{ij}$ is positively oriented around $e$, $s_{ij}=-1$ otherwise. Therefore, $\omega_A(e)$ can be considered as an approximation of the circulation of $\mathbf{u}$ around the discrete loop $\partial({\star e})$ around $e$, or, according to Stokes theorem, an approximation of the flux of vorticity $\dif\mathbf{u}^\flat$ through the dual cell ${\star e}$.

\begin{figure}
	\centering
	\includegraphics[scale=0.6]{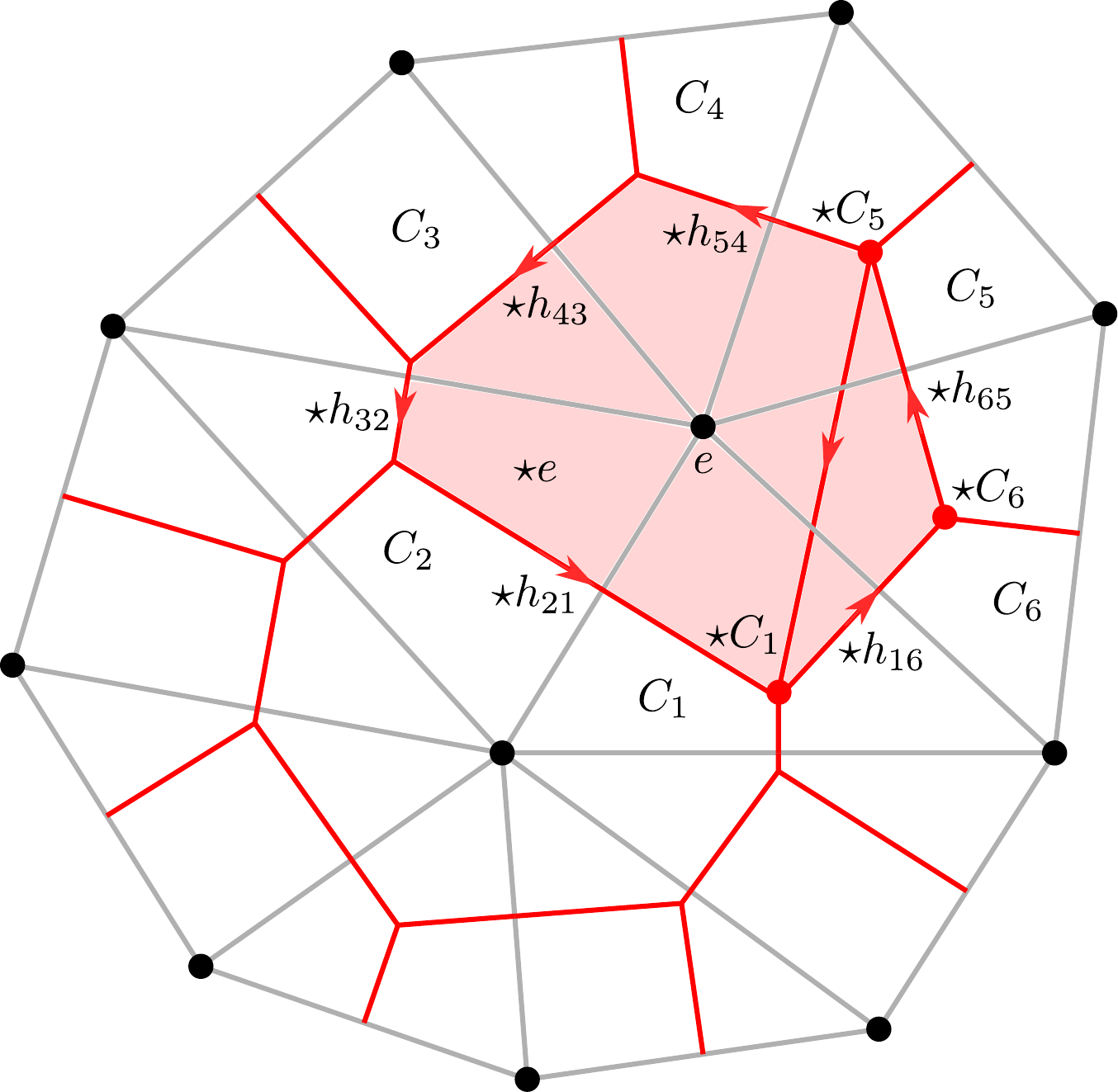}
	\caption{The \emph{total} vorticity is $\omega_A(e)=A_{21}^\flat+A_{16}^\flat+A_{65}^\flat+A_{54}^\flat+A_{43}^\flat+A_{32}^\flat$ which is an approximation of $-\int_{\partial(\star e)}\mathbf{u}^\flat$, whereas the \emph{partial} vorticity $(\dif A^\flat)_{651}$ is an approximation of $-\int_{[\star C_6,\star C_5,\star C_1]}\dif\mathbf{u}^\flat$, which is a part of the total vorticity. The signs come from those in the proposition \ref{pro:vector-field-approximation-coefficients}.}
	\label{fig:vorticity}
\end{figure}

\begin{remark}
	So far we have introduced two notions of discrete vorticity associated to a discrete vector field $A\in\mathcal{S}$: the \emph{partial} discrete vorticity $\dif A^\flat$ and the \emph{total} discrete vorticity $\omega_A$. As just explained, $\omega_A$ can be thought as the vorticity flux through a dual cell. In contrast to $\omega_A(e_{ijk})$, the relation \eqref{eq:two-vorticities} tells us that $(\dif A^\flat)_{ijk}$ only represents a fraction of this flux. More precisely, it can be thought as an approximation of the flux of vorticity through the simplex based on the triplet of adjacent cells $C_i$, $C_j$ and $C_k$, with $C_j$ and $C_k$ two cells away of each other; hence the name. See figure \ref{fig:vorticity}.
\end{remark}

The constants $K_{ijk}$ are computed such that the second condition for the flat operator is satisfied: 
\[
	K_{ijk}=s_{ijk}\frac{|{\star e}_{ijk}\cap C_i|}{|{\star e}_{ijk}|}\text{ in 2D,}\quad
	K_{ijk}=\frac{8}{3}s_{ijk}\frac{|{\star e}_{ijk}\cap C_i|}{|{\star e}_{ijk}|}\text{ in 3D},
\]
where $s_{ijk}=+1$ if the triplet of cells $C_i$, $C_j$ and $C_k$ is positively oriented around $e_{ijk}$, and $s_{ijk}=-1$ otherwise; see \citet[section 2.2.5]{Pa2009} and \citet[section 3.5]{PaMuToKaMaDe2011}.

Now that we recalled the definition of the flat operator, we can define its inverse the sharp operator $\sharp:\Omega^1_d(\mathbb{M})\to\mathcal{S}$.

\begin{definition}\label{def:sharp}
	Let $\mathbb{M}$ be a mesh approximating the fluid domain $M$. The \emph{sharp operator} $\sharp:\Omega_d^1(\mathbb{M})\to\mathcal{S}$ is defined for any $Z\in\Omega_d^1(\mathbb{M})$ by:
	\[
		\forall j\in N(i),\ j\neq i,\quad Z^\sharp_{ij}:=\frac{|h_{ij}|}{2\Omega_{ii}|{\star h}_{ij}|}Z_{ij},
	\]
	where as usual $h_{ij}$ is the shared edge (in 2D, or face in 3D) between cells $C_i$ and $C_j$, and ${\star h}_{ij}$ is the dual edge between these cells; and $Z^\sharp_{ii}=-\sum_{j\in N(i),\ j\neq i}Z^\sharp_{ij}$ thus enforcing that $Z^\sharp\in\mathcal{S}$. This definition is easily extended to the case of the differential of an extended function.
\end{definition}

Thanks to this sharp operator we can define a discrete Laplace-Beltrami operator acting on discrete functions. Interestingly, this operator corresponds to the Laplace-Beltrami operator found in \citet[section 9]{DeHiLeMa2015}.

\begin{definition}[Discrete Laplace-Beltrami operator]\label{def:discrete-Laplace-Beltrami-Laplacian}
	Let $\mathbb{M}$ be a mesh approximating the fluid domain $M$ and $F\in\Omega^0_d({\mathbb{M}})$. The \emph{Laplacian} of $F$ is the discrete function $\Delta F=\dive(\dif F)^\sharp\in\Omega^0_d({\mathbb{M}})$; that is, for any $i\in\{1,\dots,N+1\}$ by:
	\[
		(\Delta F)_i=\frac{1}{\Omega_{ii}}\sum_{j\in N(i)}(F_i-F_j)\frac{|h_{ij}|}{|{\star h}_{ij}|},
	\]
	where $h_{ij}$ is the shared edge between the cells $C_i$ and $C_j$ whereas ${\star h}_{ij}$ denotes the dual edge.
\end{definition}

We continue by defining a discrete analogue of the Hodge-De Rham Laplacian of a differential 1-form.

\begin{figure}
	\centering
	\includegraphics[scale=2.5]{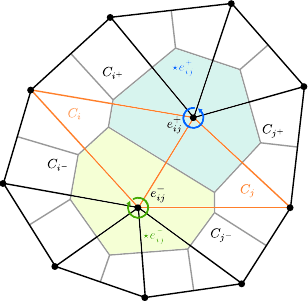}
	\caption{The primal nodes $e_{ij}^\pm$ appearing in the definition of $\Lambda(A^\flat)$, and the cells $C_{i^\pm}$ and $C_{j^\pm}$ appearing in the discrete Lie derivative $\mathbf{P}\big(\Omega^{-1}[A^\mathsf{T},\Omega L]\big)_{ij}$ in section \ref{ssec:discrete-Lie-derivatives}.}
	\label{fig:discrete-vector-Laplacian}
\end{figure}

\begin{definition}[Discrete Hodge-De Rham operator for discrete 1-forms]\label{def:discrete-Hodge-De-Rham-Laplacian}
	Let $\mathbb{M}$ be a mesh approximating the fluid domain $M$ and $A\in\mathcal{S}$. First, we define $\Lambda(A^\flat)\in\Omega^1_d(\mathbb{M})$ for any $j\in N(i)$ by
	\[
		\Lambda(A^\flat)_{ij}=\frac{1}{2}\Big[\omega_A(e_{ij}^+)\left|{\star{e}_{ij}^+}\right|-\omega_A(e_{ij}^-)\left|{\star{e}_{ij}^-}\right|\Big]\frac{|{\star h}_{ij}|}{|h_{ij}|},
	\]
	where $e_{ij}^+$ is the extremity of $h_{ij}$ such that the pair $C_i$, $C_j$ is positively oriented around $e_{ij}^+$, and negatively oriented around $e_{ij}^-$ (see figure \ref{fig:discrete-vector-Laplacian}). Then we define the \emph{Laplacian} $\Delta A^\flat\in\Omega^1_d(\mathbb{M})$ of the discrete 1-form $A^\flat$ by: $$\Delta(A^\flat)=\dif(\dive A)-\Lambda(A^\flat).$$
\end{definition}

Note that this definition is similar to the one given for the vector Laplacian in calculus. This Laplacian operator will appear as we derive the Euler-Lagrange equations associated to the variational principle for the spatially discretized Navier-Stokes-Fourier system; it will be associated to viscosity like in the smooth case.

We finish this section by a lemma which is the discrete analogue of the smooth identity $\langle f,\dive\mathbf{u}\rangle=\llangle-\dif f,\mathbf{u}\rrangle$, for any $\mathbf{u}\in\LieAlgebraFont{X}(M)$ tangent to the boundary, with on the left the $L^2$ inner product of smooth functions on $M$ and on the right the $L^2$ inner product of smooth vector fields on $M$. In the case where $\mathbf{u}$ is not tangent to the boundary, the lemma gives the discrete analogue of the familiar calculus identity
\[
	\int_Mf\dive\mathbf{u}\,\mathrm{d}x=\int_{\partial M}f\mathbf{u}\cdot\mathrm{d}\mathbf{S}-\int_M\dif f(\mathbf{u})\,\mathrm{d}x.
\]

\begin{lemma}\label{lem:div-adjoint}
	Let $\mathbb{M}$ be a mesh approximating the fluid domain $M$ and $A\in\mathcal{V}$. Then for any discrete function $F\in\Omega_d^0(\mathbb{M})$ we have:
	\[
		\langle F,\dive A\rangle_0=\llangle-\dif F,A\rrangle=-\sum_{i=1}^N\Omega_{ii}(F\cdot A)_i.
	\]
	In the case where $A\in\overline{\mathcal{V}}$, we have for any extended discrete function $F\in\Omega_d^0(\overline{\mathbb{M}})$:
	\[
		\langle F,\dive A\rangle_0=-\sum_{i=1}^N\Omega_{ii}(F\cdot A)_i+\sum_{i=1}^N\Omega_{ii}\frac{F_i+F_\partial}{2}A_i^\partial.
	\]
\end{lemma}

\begin{proof}
	Thanks to the compressibility constraint \eqref{eq:compressibility-constraint} and the discrete divergence theorem, the left hand side equals
	\begin{align*}
		\sum_{i=1}^N 2\Omega_{ii} A_{ii}F_i=\sum_{1\leq i,j\leq N}(\Omega_{ii} A_{ij}+\Omega_{jj} A_{ji})F_i=\sum_{1\leq i,j\leq N}\Omega_{jj}A_{ji}F_i,
	\end{align*}
	whereas the term in the middle equals
	\[
		\Tr\left((-\dif F)^\mathsf{T}\Omega A\right)=\sum_{1\leq i,j\leq N}(F_i-F_j)\Omega_{jj} A_{ji}=\sum_{1\leq i,j\leq N}F_i\Omega_{jj} A_{ji},
	\]
	and the right hand side reads
	\[
		\sum_{i=1}^N\Omega_{ii}(AF)_i=\sum_{1\leq i,j\leq N}\Omega_{ii} A_{ij}F_j.
	\]
	In the case $A\in\overline{\mathcal{V}}$ we have
	\[
		\langle F,\dive A\rangle_0=\sum_{1\leq i,j\leq N}(\Omega_{ii} A_{ij}+\Omega_{jj} A_{ji})F_i
		=\sum_{1\leq i,j\leq N}\Omega_{jj}A_{ji}F_i-\sum_{i=1}^N\Omega_{ii}F_iA_i^\partial,
	\]
	while
	\[
		-\sum_{i=1}^N\Omega_{ii}(F\cdot A)_i=\sum_{1\leq i,j\leq N}\Omega_{jj}A_{ji}F_i+\sum_{i=1}^N\Omega_{ii}A_{i\partial}F_\partial,
	\]
	and these expressions combine to give the result.
\end{proof}

%---------------------------------------------------------------------------------------------------

\subsubsection{Discrete momenta}\label{ssec:discrete-momenta}

In the smooth case, the physical momenta are given by $\LieAlgebraFont{X}(M)^*\cong\Omega^1(M)\otimes\Den(M)$, or by $\Omega^1(M)$ if a volume form is fixed. In the discrete case, the physical momenta are not given by $\LieAlgebraFont{d}(\mathbb{M})^*$, where $\mathbb{M}$ is a mesh approximating $M$, essentially because the physical velocity field has to satisfy in addition the compressibility constraint $\mathcal{V}$ as we have seen before. Therefore, in the discrete case, the discrete momenta are given by $\mathcal{V}^*$, which we compute now.

\begin{proposition}[Discrete momenta, {\citet[section 2]{BaGB2018}}]
	Let $\mathbb{M}$ be a mesh approximating the fluid domain $M$. Then we have:
	\begin{enumerate}[label=(\arabic*), ref=\thetheorem.(\arabic*)]
		\item $\LieAlgebraFont{d}(\mathbb{M})^*\cong\big\{L\in\gl(N):\operatorname{diag}L=0\}$,
		\item $\mathcal{V}^*\cong\Omega^1_d(\mathbb{M})\subset\gl(N)$.
	\end{enumerate}
\end{proposition}

\begin{proof}
	For the first point, the question can be reduced to $\gl(N)$ equipped with its unweighted Frobenius product by considering $\Omega\LieAlgebraFont{d}(\mathbb{M})$ instead of $\LieAlgebraFont{d}(\mathbb{M})$. We need to compute $\big[\Omega\LieAlgebraFont{d}(\mathbb{M})\big]^\bot$. Observe that all matrices $L\in\gl(N)$ such that $L_{ij}=\hat{L}_{ij}:=L_{ii}$ for all $i$, $j\in\{1,\dots,N\}$ are in $\big[\Omega\LieAlgebraFont{d}(\mathbb{M})\big]^\bot$. Moreover, the dimension of this linear subspace is $N$ whereas the codimension of $\Omega\LieAlgebraFont{d}(\mathbb{M})$ in $\gl(N)$ is $N$ according to its definition by linear equations, therefore we obtain $\LieAlgebraFont{d}(\mathbb{M})^\bot\cong\big\{L\in\gl(N):L=\hat{L}\big\}$. Consequently the Frobenius product $\llangle\cdot,\cdot\rrangle:\gl(N)/\LieAlgebraFont{d}(\mathbb{M})^\bot\times\LieAlgebraFont{d}(\mathbb{M})\to\R$ is non-degenerate and $\LieAlgebraFont{d}(\mathbb{M})^*\cong\gl(N)/\LieAlgebraFont{d}(\mathbb{M})^\bot$. This quotient is isomorphic to $\{L\in\gl(N):\operatorname{diag}L=0\big\}$, the isomorphism being given by the map $\mathbf{Q}:\gl(N)\to\{L\in\gl(N):\operatorname{diag}L=0\big\}$, $L\mapsto L-\hat L$, whose kernel is exactly $\ker\mathbf{Q}=\LieAlgebraFont{d}(\mathbb{M})^\bot$.
	
	For the second point, we consider $\LieAlgebraFont{d}(\mathbb{M})^*$ as a linear subspace of $\gl(N)$ in accordance with the first point. We need to compute $\mathcal{V}^\bot$. Take $L\in\LieAlgebraFont{d}(\mathbb{M})^*$ symmetric and $B\in\mathcal{V}$. Then we have for some diagonal matrix $D$:
	\[
		\llangle L,B\rrangle=\Tr(L\Omega B)=-\Tr\big(L(\Omega B)^\mathsf{T}\big)+\Tr(LD)=-\Tr(B^\mathsf{T}\Omega L)=-\llangle B, L\rrangle.
	\]
	Therefore, $\big\{L\in\LieAlgebraFont{d}(\mathbb{M})^*:L^\mathsf{T}=L\big\}\subset\mathcal{V}^*$. Also the dimension of this linear subspace is equal to $\dim\LieAlgebraFont{d}(\mathbb{M})^*-\dim\mathcal{V}$, so $\mathcal{V}^\bot$ is the space of symmetric $N\times N$ matrices with zero diagonal. Then $\mathcal{V}^*\cong\gl(N)^*/\mathcal{V}^\bot\cong\Omega_d^1(\mathbb{M})$ , where the isomorphism is given by the map $\mathbf{P}:\LieAlgebraFont{d}(\mathbb{M})^*\to\Omega_d^1(\mathbb{M})$, $L\mapsto\operatorname{skew}(L-\hat{L})$,  whose kernel is seen to be $\ker\mathbf{P}=\mathcal{V}^\bot$.
\end{proof}

The following result immediately follows from the previous proposition and its proof. Its importance will emerge as we derive the material and spatial semi-discrete variational principles for the Navier-Stokes-Fourier system, and beyond.

\begin{corollary}\label{cor:discrete-momenta-projection}
	Let $\mathbb{M}$ be a mesh approximating the fluid domain $M$. Then:
	\begin{enumerate}[label=(\arabic*), ref=\thetheorem.(\arabic*)]
		\item $\llangle L,B\rrangle=0$ for all $B\in\LieAlgebraFont{d}(\mathbb{M})$ is equivalent to $\mathbf{Q}(L)=0$, where $\mathbf{Q}:\gl(N)\to\LieAlgebraFont{d}(\mathbb{M})^*$ is the map $L\mapsto L-\hat{L}$;
		\item $\llangle L,B\rrangle=0$ for all $B\in\mathcal{V}$ is equivalent to $\mathbf{P}(L)=0$, where $\mathbf{P}:\gl(N)\to\mathcal{V}^*$ is the map $L\mapsto\operatorname{skew}(L-\hat{L})$.
	\end{enumerate}
\end{corollary}

\begin{remark}\label{rq:projector-incompressible-Euler}
	In the incompressible case we have $\LieAlgebraFont{d}_\text{vol}(\mathbb{M})^\bot\cong\dif\Omega^0_d(\mathbb{M})$ (the space of discrete gradients) and $\LieAlgebraFont{d}_\text{vol}(\mathbb{M})^*\cong\Omega_d^1(\mathbb{M})/\dif\Omega_d^0(\mathbb{M})$. Consequently, in this case, $\llangle L,B\rrangle=0$ for all $B\in\LieAlgebraFont{d}_\text{vol}(\mathbb{M})$ is equivalent to $\mathbf{R}(L)=0$ with $\mathbf{R}:\gl(N)\to\LieAlgebraFont{d}_\text{vol}(\mathbb{M})^*$, $L\mapsto[\operatorname{skew}L]$, that is, $L=\dif P$ for some $P\in\Omega_d^0(\mathbb{M})$; see \citet[section 1]{GMPMD2011}.
\end{remark}

\begin{lemma}\label{lem:projection}
	Let $\mathbb{M}$ be a mesh approximating the fluid domain $M$. Then for any $L\in\gl(N)$ and $A\in\mathcal{V}$ we have $\llangle L,A\rrangle=\llangle[\big]\mathbf{P}(L),A\rrangle[\big]$.
\end{lemma}

\begin{proof}
	Using the compressibility constraint \eqref{eq:compressibility-constraint} and the discrete divergence theorem $A\cdot\mathbf{1}=0$ we compute:
	\begin{align*}
	\llangle[\big]\mathbf{P}(L),A\rrangle[\big]&=
	\frac{1}{2}\sum_{1\leq i,j\leq N}\Omega_{ii}L_{ij}A_{ij}-\frac{1}{2}\sum_{1\leq i,j\leq N}\Omega_{ii}L_{ji}A_{ij}-\frac{1}{2}\sum_{1\leq i,j\leq N}\Omega_{ii}L_{jj}A_{ij}\\
	&=\sum_{1\leq i,j\leq N}\Omega_{ii}L_{ij}A_{ij}-\sum_{1\leq i,j\leq N}\Omega_{ii}L_{ji}A_{ii}\delta_{ij}+\frac{1}{2}\sum_{1\leq i,j\leq N}\Omega_{ii}L_{jj}A_{ij}\\
	&=\llangle L,A\rrangle-\sum_{i=1}^N\Omega_{ii}L_{ii}A_{ii}-\frac{1}{2}\sum_{1\leq i,j\leq N}\Omega_{jj}L_{jj}A_{ji}+\sum_{1\leq i,j\leq N}\Omega_{ii}L_{jj}A_{ii}\delta_{ij}\\
	&=\llangle L,A\rrangle.
	\end{align*}
\end{proof}

%---------------------------------------------------------------------------------------------------

\subsubsection{Densities}

In this section we define discrete densities, which will be used as a discrete counterpart for the mass density $\rho$ and entropy density $s$ of the fluid that were introduced in section \ref{ssec:smooth-spatial}.

\begin{definition}[Discrete density]
	Let $\mathbb{M}$ be a mesh approximating the fluid domain $M$. A \emph{discrete density} on $\mathbb{M}$ is an element of $\R^N$, in a similar way to discrete functions. We will denote the vector space of discrete densities by $\Den_d(\mathbb{M})$, which is in duality with $\Omega_d^0(\mathbb{M})$ by the inner product $\langle\cdot,\cdot\rangle_0$ defined in definition \ref{def:discrete-function}.
\end{definition}

Remember that discrete diffeomorphisms $\mathsf{D}(\mathbb{M})$ act linearly on the right of discrete functions $\Omega_d^0(\mathbb{M})$ by $F\cdot q=q^{-1}F$. Therefore there is a dual linear action of $\mathsf{D}(\mathbb{M})$ on $\Den_d(\mathbb{M})$ relatively to $\langle\cdot,\cdot\rangle_0$, and defined by $\langle D\bullet q,F\rangle_0=\langle D,F\cdot q^{-1}\rangle_0$, for any $D\in\Den_d(\mathbb{M})$ and $F\in\Omega^0_d(\mathbb{M})$. We can compute explicitly that:
\[
	D\bullet q=\Omega^{-1}q^\mathsf{T}\Omega D,
\]
and we can also compute the infinitesimal action of $\LieAlgebraFont{d}(\mathbb{M})$ on the right of $\Den_d(\mathbb{M})$ as well, which is given for any $A\in\LieAlgebraFont{d}(\mathbb{M})$ and $D\in\Den_d(\mathbb{M})$ by:
\[
	D\bullet A=\Omega^{-1}A^\mathsf{T}\Omega D.
\]
This action satisfies the relation $\langle D\bullet A,F\rangle_0=-\langle D,F\cdot A\rangle_0$ for all $A\in\LieAlgebraFont{d}(\mathbb{M})$, $D\in\Den_d(\mathbb{M})$ and $F\in\Omega_d^0(\mathbb{M})$.

The following lemma can be understood as the discrete analogue in $\Den_d(\mathbb{M})$ of the relation $\dive(\rho\mathbf{u})\,\mathrm{d}x=\rho\dive\mathbf{u}\,\mathrm{d}x+\dif\rho(\mathbf{u})\,\mathrm{d}x$ in $\Den(M)$.

\begin{lemma}\label{lem:two-actions}
	Let $\mathbb{M}$ be a mesh approximating the fluid domain $M$. Then for any $A\in\mathcal{V}$ and $D\in\Den_d(\mathbb{M})$ we have for any $i\in\{1,\dots,N\}$:
	\[
		(D\bullet A)_i=(\dive A)_i D_i+(D\cdot A)_i=2A_{ii}D_i-\sum_{j=1}^N A_{ij}D_j,
	\]
	where in $D\cdot A$, $D$ is viewed as an element of $\Omega_d^0(\mathbb{M})$\footnote{As vector spaces $\Den_d(\mathbb{M})$ and $\Omega_d^0(\mathbb{M})$ are identical, but not as $\mathsf{D}(\mathbb{M})$-spaces. The same can be said in the smooth case: once a measure has been fixed, one can identify densities as functions, but $\Diff(M)$ acts differently on each of them.}.
\end{lemma}

\begin{proof}
	Because $A$ satisfies the compressibility constraint \eqref{eq:compressibility-constraint}, we compute that:
	\[
		(D\bullet A)_i=\sum_{j=1}^N\Omega_{ii}^{-1}A_{ji}\Omega_{jj}D_j
		=\sum_{j=1}^N(2A_{ii}\delta_{ij}-A_{ij})D_j
		=2A_{ii}D_i-\sum_{j=1}^N A_{ij}D_j.
	\]
\end{proof}

The following lemma is the discrete analogue of $\int_M\dive(\rho\mathbf{u})\,\mathrm{d}x=0$, for any $\mathbf{u}\in\LieAlgebraFont{X}(M)$ which is tangent to the boundary.

\begin{lemma}\label{lem:integral-action-density}
	Let $\mathbb{M}$ be a mesh approximating the fluid domain $M$. Then for any $A\in\mathcal{V}$ and $D\in\Den_d(\mathbb{M})$ we have:
	\[
		\sum_{i=1}^N\Omega_{ii}(D\bullet A)_i=0.
	\]
\end{lemma}

\begin{proof}
	Using the previous lemma, the compressibility constraint \eqref{eq:compressibility-constraint}, and the discrete divergence theorem, we compute:
	\begin{align*}
		\sum_{i=1}^N\Omega_{ii}(D\bullet A)_i
		&=\sum_{1\leq i,j\leq N}2\Omega_{ii}A_{ii}D_i\delta_{ij}-\sum_{1\leq i,j\leq N}\Omega_{ii}A_{ij}D_j\\
		&=\sum_{1\leq i,j\leq N}\Omega_{ii}A_{ij}D_i+\sum_{1\leq i,j\leq N}\Omega_{jj}A_{ji}D_i-\sum_{1\leq i,j\leq N}\Omega_{ii}A_{ij}D_j\\
		&=0.
	\end{align*}
\end{proof}

We conclude this section with a lemma that will be useful when we will derive the Euler-Lagrange equations from the forthcoming variational principle.

\begin{lemma}\label{lem:pairing-change}
	Let $\mathbb{M}$ be a mesh approximating the fluid domain $M$. Then for any $F\in\Omega_0^d(\mathbb{M})$, $D\in\Den(\mathbb{M})$ and $A\in\LieAlgebraFont{d}(\mathbb{M})$ we have:
	\[
		\langle F,D\bullet A\rangle_0=\llangle[\big] DF^\mathsf{T},A\rrangle[\big].
	\]
\end{lemma}

\begin{proof}
	The left hand side is equal to:
	\[
		\langle D,AF\rangle_0=D^\mathsf{T}\Omega AF=\sum_{1\leq i,j\leq N}D_i\Omega_{ii}A_{ij}F_j,
	\]
	whereas the right hand side is equal to:
	\[
		\Tr(FD^\mathsf{T}\Omega A)=\sum_{1\leq i,j\leq N}(FD^\mathsf{T})_{ij}(\Omega A)_{ji}=\sum_{1\leq i,j\leq N}F_iD_j\Omega_{jj}A_{ji}.
	\]
\end{proof}

\begin{remark}
	We have a similar identity when considering the action of $\LieAlgebraFont{d}(\mathbb{M})$ on discrete functions. For any $F\in\Omega_0^d(\mathbb{M})$, $D\in\Den(\mathbb{M})$ and $A\in\LieAlgebraFont{d}(\mathbb{M})$, we have $\langle D,F\cdot A\rangle_0=-\llangle[\big] DF^\mathsf{T},A\rrangle[\big]$.
\end{remark}

\begin{remark}
	The pairing $\llangle L,B\rrangle=\Tr(L^\mathsf{T}\Omega B)$ introduced before can also be given a wider interpretation using discrete densities. Indeed, the diagonal matrix $\Omega$ can be seen as a discrete density and $\Omega B$ as the discrete analogue of $\mathbf{v}\otimes\mathrm{d}x$. More generally, the discrete analogue of a 1-form density $\mathbf{v}\otimes\rho\,\mathrm{d}x$ is $DB$, where $D$ is a discrete density regarded as a diagonal matrix and $B$ a discrete vector field. The pairing $\llangle,\rrangle$ is therefore the discrete analogue of the pairing \eqref{eq:1-form-densities-vector-fields-pairing}.
\end{remark}

%---------------------------------------------------------------------------------------------------

\subsubsection{Lie derivatives}\label{ssec:discrete-Lie-derivatives}

As we have already seen in section \ref{ssec:compressible-Euler}, the adjoint representation of $\LieAlgebraFont{X}(M)$ is given by $\ad_{\mathbf{u}}\mathbf{v}=-[\mathbf{u},\mathbf{v}]=-\lie_\mathbf{u}\mathbf{v}$, and its adjoint $\ad_\mathbf{u}^*:\LieAlgebraFont{X}(M)^*\to\LieAlgebraFont{X}(M)^*$, with respect to the pairing between vector fields and 1-form densities, is given by $\ad_\mathbf{u}^*(\alpha\otimes\mu)=\lie_\mathbf{u}(\alpha\otimes\mu)$. We can do the same at the discrete level to obtain \emph{discrete Lie derivatives} of discrete momenta.

On $\LieAlgebraFont{d}(\mathbb{M})$, the adjoint action is given by $\ad_A B=[A,B]$ and we compute easily that its adjoint $\ad^*_A:\LieAlgebraFont{d}(\mathbb{M})^*\to\LieAlgebraFont{d}(\mathbb{M})^*$ with respect to the pairing $\llangle\cdot,\cdot\rrangle$ is given for any $L\in\LieAlgebraFont{d}(\mathbb{M})^*$ by:
\[
	\ad_A^*L=\mathbf{Q}\big(\Omega^{-1}[A^\mathsf{T},\Omega L]\big).
\]

\begin{remark}
	In the incompressible case, we have in a similar way that the operator $\ad_A^*:\LieAlgebraFont{d}_\text{vol}(\mathbb{M})^*\to\LieAlgebraFont{d}_\text{vol}(\mathbb{M})^*$ is given by $\ad_A^*L=\mathbf{R}\big(\Omega^{-1}[A^\mathsf{T},\Omega L]\big)=\Omega^{-1}[A^\mathsf{T},\Omega L]$.
\end{remark}

Recall that the operator $\ad_A^*$ is defined by $\llangle\ad_A^*L,B\rrangle=\llangle L,\ad_A B\rrangle$, for any $L\in\LieAlgebraFont{d}(\mathbb{M})^*$ and $B\in\LieAlgebraFont{d}(\mathbb{M})$\footnote{Thus $\ad_A^*$ is the dual map associated to $\ad_A^*$, but is \emph{not} the dual (coadjoint) representation associated to $\ad_A$.}. However, because of the nonholonomic compressibility constraint, in the variational principle to come we will only have that $B\in\mathcal{V}$. Therefore, the operator that we need to understand is $\mathbf{P}\big(\Omega^{-1}[A^\mathsf{T},\Omega L]\big)=\operatorname{skew}\ad_A^*L$, which is the subject of the following proposition.

\begin{proposition}[{\citet[lemma 3.1]{BaGB2018}}]\label{pro:Lie-derivative-1-form-density}
	Let $\mathbb{M}$ be a mesh approximating the fluid domain $M$ and let $A\in\LieAlgebraFont{d}(\mathbb{M})$, $D\in\Den_d(\mathbb{M})$ and $B\in\mathcal{S}$. Define $L\in\LieAlgebraFont{d}(\mathbb{M})^*$ by $L_{ij}=D_i B_{ij}^\flat$ for all $i$, $j\in\{1,\dots,N\}$. Then we have:
	\begin{align*}
		\mathbf{P}\big(\Omega^{-1}[A^\mathsf{T},\Omega L]\big)_{ij}=
		\omega_B(e_{ij}^+)&\left(K_i^+\overline{D}_{ji^+}A_{ii^+}+K_j^+\overline{D}_{ij^+}A_{jj^+}\right)
		-\omega_B(e_{ij}^-)\left(K_i^-\overline{D}_{ji^-}A_{ii^-}+K_j^-\overline{D}_{ij^-}A_{jj^-}\right)\\
		&\quad+\overline{D}_{ij}\left(\sum_{k\in N(i)}A_{ik}B_{ik}^\flat-\sum_{k\in N(j)}A_{jk}B_{jk}^\flat\right)-\overline{D\bullet A}_{ij}B_{ij}^\flat,
	\end{align*}
	where $\overline{D}_{ij}=\frac{1}{2}(D_i+D_j)$, and $K_i^\pm=\frac{|{\star e}_{ij}^\pm\cap C_i|}{|{\star e}_{ij}^\pm|}$, and the other notations are explained in figure \ref{fig:discrete-vector-Laplacian}.
\end{proposition}

Using Cartan's formula $\lie_\mathbf{u}=\ins_\mathbf{u}\dif+\dif\ins_\mathbf{u}$, the formula of the previous proposition can be considered as the discrete analogue of the following smooth formula in $\Den(M)$, which motivates the forthcoming definition:
\[
	\lie_\mathbf{u}(\mathbf{v}^\flat\otimes\rho\,\mathrm{d}x)=\big(\ins_\mathbf{u}\dif\mathbf{v}^\flat+\dif\langle\mathbf{v}^\flat,\mathbf{u}\rangle\big)\otimes\rho\,\mathrm{d}x+\mathbf{v}^\flat\otimes\dive(\rho\mathbf{u})\,\mathrm{d}x.
\]

\begin{definition}[Lie derivative of a 1-form density]
	Let $\mathbb{M}$ be a mesh approximating the fluid domain $M$. Let $A\in\LieAlgebraFont{d}(\mathbb{M})$, $D\in\Den_d(\mathbb{M})$ and $B\in\mathcal{S}$. We will denote by $\lie_A(DB^b)_{ij}=\mathbf{P}\big(\Omega^{-1}[A^\mathsf{T},\Omega L]\big)_{ij}$ where $L\in\LieAlgebraFont{d}(\mathbb{M})^*$ is defined by $L_{ij}=D_i B_{ij}^\flat$ for all $i$, $j\in\{1,\dots,N\}$ (or $D$ is viewed as a diagonal matrix). Then $\lie_A(DB^b)\in\LieAlgebraFont{d}(\mathbb{M})^*$ is called the \emph{Lie derivative} of the 1-form density $L$.
\end{definition}

\begin{remark}
	These Lie derivatives admit an other algebraic description akin to the Cartan formula in differential geometry. Following \cite{PaMuToKaMaDe2011}, define for a discrete 1-form $F\in\Omega^1_d(\mathbb{M})$ the \emph{interior product} $\ins_A F\in\Omega_d^0(\mathbb{M})$ by a discrete vector field $A\in\LieAlgebraFont{d}(\mathbb{M})$ as:
	\[
		(\ins_AF)_i:=(AF^\mathsf{T})_{ii}=(FA^\mathsf{T})_{ii},
	\]
	as well as the \emph{interior product} $\ins_A F\in\Omega_d^1(\mathbb{M})$ of a discrete 2-form $F\in\Omega_d^2(\mathbb{M})$ by a discrete vector field $A\in\LieAlgebraFont{d}(\mathbb{M})$ as:
	\[
		(\ins_AF)_{ij}=\sum_{k}\left(F_{ikj}A_{ik}-F_{jki}A_{jk}\right).
	\]
	Then the Lie derivative $\lie_A F\in\Omega^1_d(\mathbb{M})$ of a discrete 1-form $F\in\Omega^1_d(\mathbb{M})$ is defined by $\lie_A=-(\ins_A\dif+\dif\ins_A)$. Consequently we obtain the following expression:
	\begin{align*}
		(\lie_AF)_{ij}&=-(\ins_A\dif F)_{ij}+(\dif\ins_AF)_{ij}\\
		&=-\sum_k(F_{ik}+F_{kj}+F_{ji})A_{ik}+\sum_k(F_{jk}+F_{ki}+F_{ij})A_{jk}-(FA^\mathsf{T})_{jj}+(FA^\mathsf{T})_{ii}\\
		&=-\big[AF-(AF)^\mathsf{T}\big]_{ij}\\
		&=-\big(AF+FA^\mathsf{T}\big)_{ij}.
	\end{align*}
	In the incompressible case, $A\in\LieAlgebraFont{d}_\text{vol}(\mathbb{M})$ and from the relation $A^\mathsf{T}\Omega+\Omega A=0$ we obtain that $\ad_A^*F=\Omega^{-1}[A^\mathsf{T},\Omega L]=AF+FA^\mathsf{T}=\lie_A F$. Now in the compressible case, $A\in\mathcal{V}$, and we have another expression of $\lie_A F$. From the identity $A^\mathsf{T}\Omega+\Omega A=\Omega\dive(A)$ where $\dive(A)$ is viewed as a diagonal matrix, we obtain:
	\begin{align*}
		\Omega^{-1}\big[A^\mathsf{T},\Omega F\big]=-(AF+FA^\mathsf{T})+\dive(A)F,
	\end{align*}
	hence:
	\[
		\ad_A^*L=\mathbf{P}\big(\Omega^{-1}[A^\mathsf{T},\Omega L]\big)=-(AF+FA^\mathsf{T})+\overline{\dive A}\circ F=\lie_AF+\overline{\dive A}\circ F,
	\]
	with $\big(\overline{\dive A}\circ F\big)_{ij}=\dive(A)_{ij}F_{ij}$.
\end{remark}

%---------------------------------------------------------------------------------------------------

\subsubsection{Summary comparing smooth objects and their discrete counterparts}

We provide a summary in table \ref{tab:smooth-discrete-summary} which lists smooth objects and their discrete counterparts in the discrete exterior calculus developed so far. It will become handy to draw comparisons when we tackle the spatial discretization of the Navier-Stokes-Fourier system.

\begin{table}
	\centering
	\renewcommand*{\arraystretch}{1.4}
	\begin{tabular}{|c||c|c|}
		\hline
		& smooth object on $M$ & discrete counterpart on $\mathbb{M}$ \\[3pt]
		\hline\hline
		function & $f\in\smooth{M}$ & $F\in\Omega_d^0(\mathbb{M})$ \\[3pt]
		\hline 
		global pairing of functions & $\langle f,g\rangle=\int_M fg\,\mathrm{d}x$ & $\langle F,G\rangle=\Tr(F^T\Omega G)$ \\
		\hline 
		density & $\rho\,\mathrm{d}x\in\Den(M)$ & $D\in\Den_d(M)$ \\
		\hline 
		diffeomorphism & $\varphi\in\Diff(M)$ & $q\in\mathsf{D}(\mathbb{M})$ \\
		\hline 
		vector field & $\mathbf{u}=\dot{\varphi}\circ\varphi^{-1}\in\LieAlgebraFont{X}(M)$ & $A=\dot{q}q^{-1}\in\LieAlgebraFont{d}(\mathbb{M})$ \\
		\hline
		divergence theorem & $\int_{C_i}\dive\mathbf{u}\,\mathrm{d}x=\int_{\partial C_i}\mathbf{u}\cdot\mathrm{d}\mathbf{S}$ & $A\cdot\mathbf{1}=0$ \\
		\hline
		global pairing of vector fields & $\llangle\mathbf{u}^\flat,\mathbf{v}\rrangle=\int_M\mathbf{u}\cdot\mathbf{v}\,\mathrm{d}x$ & $\llangle A^\flat,B\rrangle=\Tr\big((A^\flat)^\mathsf{T}\Omega B\big)$ \\
		\hline
		group action on functions & $f\mapsto f\circ\varphi$ & $F\mapsto F\cdot q=q^{-1}F$ \\
		\hline
		group action on densities & $\rho\,\mathrm{d}x\mapsto(\rho\circ\varphi)J_\varphi\,\mathrm{d}x$ & $D\mapsto D\bullet q=\Omega^{-1}q^\mathsf{T}\Omega D$ \\
		\hline
		Lie algebra action on functions & $f\mapsto\dif f(\mathbf{u})$ & $F\mapsto F\cdot A=-AF$ \\
		\hline 
		Lie algebra action on densities & $\rho\,\mathrm{d}x\mapsto\dive(\rho\mathbf{u})\,\mathrm{d}x$ & $D\mapsto D\bullet A=\Omega^{-1}A^\mathsf{T}\Omega D$ \\
		\hline
		Lie derivative of a 1-form & $\alpha\mapsto\lie_\mathbf{u}\alpha$ & $L\mapsto\mathbf{P}(\Omega^{-1}[A^\mathsf{T},\Omega L])$ \\
		\hline
		Lie derivative of a 1-form density & $\alpha\otimes\rho\,\mathrm{d}x\mapsto\lie_\mathbf{u}(\alpha\otimes\rho\,\mathrm{d}x)$ & $L_{ij}=D_iB^\flat_{ij}\mapsto\mathbf{P}(\Omega^{-1}[A^\mathsf{T},\Omega L])$ \\
		\hline
	\end{tabular}
	\caption{Comparison of smooth objects and their discrete counterparts.}
	\label{tab:smooth-discrete-summary}
\end{table}

%---------------------------------------------------------------------------------------------------

\subsection{Discretization of the phenomenological and variational constraints}\label{ssec:constraints-discretization}

The presence of a gradient (or more generally of a covariant derivative) in both the phenomenological and variational constraints seems to prevent a discretization according to the discrete framework presented in section \ref{sec:DEC}. Nevertheless, the phenomenological constraint, and hence the {associated} variational constraint, can be recast into constraints that are amenable to discretization in terms of the discrete exterior calculus developed so far, which is the subject of this section.
. The difficult point is the discretization of the quantity $\bm{\sigma}^\text{fr}:\nabla\mathbf{v}$, and more precisely of $\Def\mathbf{u}:\Def\mathbf{v}$, for any vector fields $\mathbf{u}$, $\mathbf{v}\in\LieAlgebraFont{X}(M)$, where $\Def\mathbf{u}=\frac{1}{2}(\nabla\mathbf{u}+\nabla\mathbf{u}^\mathsf{T})$, the so-called \emph{linear strain tensor} in continuum mechanics.

\begin{lemma}\label{lem:variational-constraint-recasting}
	For any vector fields $\mathbf{u}$, $\mathbf{v}\in\LieAlgebraFont{X}(M)$ we have:
	\[
		\Def\mathbf{u}:\Def\mathbf{v}=
		\dive\mathbf{u}\dive\mathbf{v} + \frac{1}{2}\curl\mathbf{u}\cdot\curl\mathbf{v} + \dive\big[\nabla_\mathbf{v}\mathbf{u} - (\dive\mathbf{u})\mathbf{v}\big].
	\]
\end{lemma}

\begin{proof}
	Recall that for any 2-tensor $\mathbf{T}$ and vector field $\mathbf{v}$ we have the identity $\dive(\ins_\mathbf{v}\mathbf{T})=\ins_\mathbf{v}\dive\mathbf{T}+\mathbf{T}:\nabla\mathbf{v}$. In particular, for $\mathbf{T}=\Def\mathbf{u}$ which is symmetric, we have:
	\[
		\Def\mathbf{u}:\Def\mathbf{v}=\dive\big(\ins_\mathbf{v}\Def\mathbf{u}\big)-\ins_\mathbf{v}\dive\Def\mathbf{u}.
	\]
	We compute each term separately. Firstly we get:
	\[
		\ins_\mathbf{v}\Def\mathbf{u}=\nabla_\mathbf{v}\mathbf{u}+\frac{1}{2}(\nabla\mathbf{u}^\mathsf{T}-\nabla\mathbf{u})\cdot\mathbf{v}=\nabla_\mathbf{v}\mathbf{u}-\frac{1}{2}(\curl\mathbf{u})\times\mathbf{v}.
	\]
	Secondly we have:
	\[
		\dive\Def\mathbf{u}=\frac{1}{2}\Delta\mathbf{u}+\frac{1}{2}\nabla\dive\mathbf{u}=\nabla\dive\mathbf{u}-\frac{1}{2}\curl\curl\mathbf{u},
	\]
	and owing to the fact that $\curl\curl\mathbf{u}\cdot\mathbf{v}=\curl\mathbf{u}\cdot\curl\mathbf{v}+\dive\big[(\curl\mathbf{u})\times\mathbf{v}\big]$, we obtain that
	\[
		\ins_\mathbf{v}\dive\Def\mathbf{u}=-\dive\mathbf{u}\dive\mathbf{v}-\frac{1}{2}\curl\mathbf{u}\cdot\curl\mathbf{v}+\dive\left[(\dive\mathbf{u})\mathbf{v}-\frac{1}{2}(\curl\mathbf{u})\times\mathbf{v}\right].
	\]
	Combining the terms together, we finally obtain the desired result.
\end{proof}

\begin{corollary}\label{cor:phenomenological-constraint-recasting}
	For any vector field $\mathbf{u}\in\LieAlgebraFont{X}(M)$ we have:
	\[
		\Def\mathbf{u}:\Def\mathbf{u}=
		(\dive\mathbf{u})^2 + \frac{1}{2}\|\curl\mathbf{u}\|^2 + \dive\left[\left(\lie_\mathbf{u}\mathbf{u}^\flat-\frac{1}{2}\dif\|\mathbf{u}\|^2\right)^\sharp - (\dive\mathbf{u} )\mathbf{u}\right].
	\]
\end{corollary}

\begin{proof}
	This directly follows from the previous lemma and the fact that $(\nabla_\mathbf{u}\mathbf{u})^\flat=\lie_\mathbf{u}\mathbf{u}^\flat-\frac{1}{2}\dif\|\mathbf{u}\|^2$.
\end{proof}

Before giving the full discretization of the variational principles \ref{ssec:smooth-material} and \ref{ssec:smooth-spatial}, we need to figure out a proper discretization of the phenomenological constraint and its associated variational constraint. We will begin with the variational constraint \eqref{eq:spatial-phenomenological-constraint} in the spatial formalism, and its associated variational constraint \eqref{eq:spatial-variational-constraint}. As we have just seen, an important part of this constraint reads:
\[
\Def\mathbf{u}:\Def\mathbf{v}=
\dive\mathbf{u}\dive\mathbf{v} + \frac{1}{2}\curl\mathbf{u}\cdot\curl\mathbf{v} + \dive\big[\nabla_\mathbf{v}\mathbf{u} - (\dive\mathbf{u})\mathbf{v}\big],
\]
for any vector field $\mathbf{u}$  and $\mathbf{v}\in\LieAlgebraFont{X}(M)$. Each term on the right-hand side can be discretized within the framework developed in \ref{sec:DEC} at a cell $C_i$ of a mesh $\mathbb{M}$ approximating the fluid domain $M$, except for the term $\curl\mathbf{u}\cdot\curl\mathbf{v}$ that we are going to address now. Indeed, as the covariant derivative $\nabla_\mathbf{v}\mathbf{u}$ sits inside a divergence term, \emph{any} discretization of it can be used since anyway it will vanish when deriving later the discrete Euler-Lagrange equations, thanks to lemma \ref{lem:integral-discrete-divergence} (in that case $\mathbf{v}$ needs to be zero on the boundary such that by locality of the covariant derivative we have $\nabla_\mathbf{v}\mathbf{u}=0$; we will assume that a discretization having a similar property exists). 

\begin{figure}
	\centering
	\includegraphics[scale=2.5]{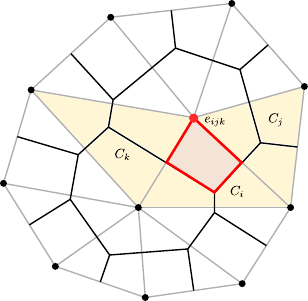}
	\caption{The geometric shape, drawn in red, determined by the intersection of $\star e_{ijk}$ and $C_i$ is called the kite associated to $e=e_{ijk}$ and $C_i$. Alternatively, in terms of primal cells, such a kite is determined by a triplet $(i,j,k)$ of adjacent cells $C_i$, $C_j$ and $C_k$. Note however that the triplet $(i,j,k)$ and $(i,k,j)$ determine the same geometric kite, albeit changing its orientation. Observe also how kites partition $\mathbb{M}$.}
	\label{fig:kite}
\end{figure}

For that matter, we first need to understand how to discretize $\int_M\curl\mathbf{u}\cdot\curl\mathbf{v}$. Note that in $2D$ this integral equals $\int\omega_\mathbf{u}\,\omega_\mathbf{v}\,\mathrm{d}{x}$ where $\omega_{\mathbf{u}}$ is the function defined by $\dif\mathbf{u}^\flat=\omega_\mathbf{u}\,\dif x\wedge\dif y$, whereas in 3D this integral equals $\int_M\dif\mathbf{u}^\flat\wedge\star\dif\mathbf{v}^\flat$. We first discretize $\omega_\mathbf{u}$ into $\omega_A$, where $A$ is a discrete vector field approximating $\mathbf{u}$ (see \ref{eq:total-vorticity}). Next, we observe that the mesh $\mathbb{M}$ can be partitioned into \emph{kites}, which are the geometric shapes determined by the intersection between a dual cell $\star e$ associated to a primal node $e$, and a primal cell $C_i$. Alternatively, such a kite is determined by a triplet of adjacent cells $C_i$, $C_j$ and $C_k$; see figure \ref{fig:kite}.

Using this partition of $\mathbb{M}$ and noting that the value of the discrete vorticity on a kite corresponds to its value on the dual cell it belongs to, we can approximate our integral as follows:
\begin{align*}
	\int_M\omega_\mathbf{u}\,\omega_\mathbf{v}\,\mathrm{d}\mathbf{x}
	&\approx\sum_{\text{kite }\kappa}\omega_A(e_\kappa)\omega_B(e_\kappa)\,|{\star e_\kappa}\cap C_\kappa|\\
	&=\sum_{\substack{\text{adjacent cells}\\ C_i,\ C_j,\ C_k}}\omega_A(e_{ijk})\,\omega_B(e_{ijk})\,|{\star e_{ijk}}\cap C_i|\\
	&=\sum_{i=1}^N\Omega_{ii}\sum_{j,k\in N(i)}W_{ijk}(\dif A^\flat)_{ijk}(\dif B^\flat)_{ijk},
\end{align*}
where we defined
\[
W_{ijk}=\frac{|{\star e_{ijk}}|^2}{2\Omega_{ii}|{\star e_{ijk}}\cap C_i|}.
\]
Consequently we can approximate the smooth function $\omega_\mathbf{u}\,\omega_\mathbf{v}$ on a cell $C_i$ by the formula
\[
\sum_{j,k\in N(i)}W_{ijk}(\dif A^\flat)_{ijk}(\dif B^\flat)_{ijk}.
\]
Notice that, although we started in the 2D setting, this formula also makes sense for the 3D setting. For our purpose we need to rearrange the terms in the previous approximation. Firstly, summing over the kites of $\mathbb{M}$ accounts to summing over the primal nodes of $\mathbb{M}$, therefore:
\begin{align*}
	\sum_{\substack{\text{adjacent cells}\\ C_i,\ C_j,\ C_k}}(\mathbf{d}A^\flat)_{ijk}(\mathbf{d}B^\flat)_{ijk}\frac{|{\star e_{ijk}}|^2}{|{\star e_{ijk}}\cap C_i|}
	& = \sum_{\text{nodes }e}\sum_{\substack{\text{adjacent cells}\\ C_i\cap C_j\cap C_k=e}}\omega_A(e_{ijk})\,\omega_B(e_{ijk})K_{ijk}^2\frac{|{\star e_{ijk}}|^2}{|{\star e_{ijk}}\cap C_i|}\\
	& = \sum_{\text{nodes }e}\sum_{\substack{\text{adjacent cells}\\ C_i\cap C_j\cap C_k=e}}\omega_A(e_{ijk})\,\omega_B(e_{ijk})\,|{\star e_{ijk}}\cap C_i|,
\end{align*}
but the primal node $e_{ijk}$ doesn't depend on $(i,j,k)$ (any triplet with the same intersection yields the same node) so that we can factor out $\omega_A(e_{ijk})\,\omega_B(e_{ijk})$ and thus obtain:
\begin{equation*}
	\sum_{\text{nodes }e}\omega_A(e)\,\omega_B(e)\sum_{\substack{\text{adjacent cells}\\ C_i\cap C_j\cap C_k=e}}|{\star e_{ijk}}\cap C_i|=\frac{1}{2}\sum_{\text{nodes }e}\omega_A(e)\,\omega_B(e)\,|{\star e}|,
\end{equation*}
where the $\frac{1}{2}$ factor comes from the fact that enumerating triplets of adjacent cells $(i,j,k)$ is equivalent to enumerating cells $i$ and considering adjacent cells $k$ and $j$, but then $e_{ijk}=e_{ikj}$ so we have twice the desired contribution. Secondly, the sum on the primal nodes can be rewritten in another way using the definition \ref{eq:total-vorticity}\footnote{The round arrow is an intuitive shortcut replacing the more sophisticated notation in the definition \ref{eq:total-vorticity}.}:
\begin{align*}
	\sum_{\text{nodes }e}\omega_A(e)\,\omega_B(e)\,|{\star e}|
	&=\sum_{\text{nodes }e}\omega_A(e)\left(\sum_{i,j\circlearrowleft e}s_{ij}B_{ij}^\flat\right)|{\star e}|\\
	&=\sum_{i=1}^N\sum_{j\in N(i)}\sum_{\substack{\text{nodes }e\\ e=C_i\cap C_j}}\omega_A(e)s_{ij}B_{ij}^\flat\,|{\star e}|,
\end{align*}
but for each pair $(i,j)$ with $j\in N(i)$, there are only two primal nodes in the intersection $C_i\cap C_j$, so we get the following formula:
\begin{equation*}
	\sum_{\text{nodes }e}\omega_A(e)\,\omega_B(e)\,|{\star e}|
	=\sum_{i=1}^N\sum_{j\in N(i)}\Big[\omega_A(e_{ij}^+)\left|{\star e_{ij}^+}\right|-\omega_A(e_{ij}^-)\left|{\star e_{ij}^-}\right|\Big]B_{ij}^\flat,
\end{equation*}
where, by convention, $e_{ij}^+$ is the primal node such that we go from $C_i$ to $C_j$ in the counterclockwise direction, and $e_{ij}^-$ is the primal node such that we go from $C_i$ to $C_j$ in the clockwise direction (see figure \ref{fig:discrete-vector-Laplacian}). Finally, using the definition \ref{eq:flat-operator} of the flat operator for adjacent cells, we obtain:
\[
\sum_{i=1}^N\Omega_{ii}\sum_{j,k\in N(i)}W_{ijk}(\dif A^\flat)_{ijk}(\dif B^\flat)_{ijk}=\llangle[\big]\Lambda(A^\flat),B\rrangle[\big],
\]
where $\Lambda(A^\flat)$ has been defined in \ref{def:discrete-Hodge-De-Rham-Laplacian}.

\begin{remark}
	From the above formula, we see that the discrete operator $\Lambda:\Omega^1_d(\mathbb{M})\to\Omega^1_d(\mathbb{M})$ is the discrete analogue of the smooth operator $\curl^*\circ\curl$, relatively to the $L^2$ pairing of vector fields tangent to the boundary. Therefore, in the forthcoming variational principle, we will exactly reproduce at the discrete level the following computation:
	\begin{align*}
		\int_M \operatorname{Def}\mathbf{u}:\operatorname{Def}\mathbf{v}\,\mathrm{d}{x}
		&=\int_{M}\left[\dive^*(\dive\mathbf{u})\cdot\mathbf{v}+\frac{1}{2}\curl^*(\curl\mathbf{u})\cdot\mathbf{v}\right]\mathrm{d}{x}\\
		&=\int_{M}\left(-\nabla\dive\mathbf{u}+\frac{1}{2}\curl\curl\mathbf{u}\right)\cdot\mathbf{v}\,\mathrm{d}{x}\\
		&=\int_{M}\left[-\frac{1}{2}\Delta\mathbf{u}\cdot\mathbf{v}-\frac{1}{2}(\nabla\dive\mathbf{u})\cdot\mathbf{v}\right]\mathrm{d}{x}.
	\end{align*}
\end{remark}

We summarize our result in the following proposition.

\begin{proposition}\label{pro:discrete-curl-curl}
	Let $\mathbb{M}$ be a mesh approximating the fluid domain $M$. Given two discrete vector fields $A$, $B\in\LieAlgebraFont{d}(\mathbb{M})$ approximating two continuous vector fields $\mathbf{u}$ and $\mathbf{v}\in\LieAlgebraFont{X}(M)$ (these are supposed to be tangent to the boundary $\partial M$), either the function $\omega_\mathbf{u}\,\omega_\mathbf{v}$ or $\dif\mathbf{u}^\flat\wedge\star\dif\mathbf{v}^\flat$ can be approximated on a cell $C_i$ by:
	\[
	\big(\dif A^\flat\wedge\star\dif B^\flat\big)_i:=\sum_{j,k\in N(i)}W_{ijk}(\dif A^\flat)_{ijk}(\dif B^\flat)_{ijk}.
	\]
	Moreover we have the following identities:
	\begin{align*}
		\sum_{i=1}^N\Omega_{ii}\sum_{j,k\in N(i)}W_{ijk}(\dif A^\flat)_{ijk}(\dif B^\flat)_{ijk}
		&=\frac{1}{2}\sum_{i=1}^N\sum_{j,k\in N(i)}\kappa_{ijk}\,\omega_A(e_{ijk})\,\omega_B(e_{ijk})\\
		&=\frac{1}{2}\sum_{\text{nodes }e}\omega_A(e)\,\omega_B(e)\,|{\star e}|\\
		&=\llangle[\big]\Lambda(A^\flat),B\rrangle[\big],
	\end{align*}
	where $\kappa_{ijk}$ denotes the (unsigned) area $|\star e_{ijk}\cap C_i|$ of the kite determined by the triplet of adjacent cells $C_i$, $C_j$ and $C_k$.
\end{proposition}

Now, remembering from \eqref{eq:Newton-Fourier} that for the Navier-Stokes-Fourier system we chose $\bm{\sigma}^\text{fr}=2\mu\Def\mathbf{u}+\left(\zeta-\frac{2}{3}\right)(\dive\mathbf{u})I_2$ with $\mu$ the shear viscosity and $\zeta$ the bulk viscosity, so we finally have setting $\tilde{\mu}=\zeta+\frac{4}{3}\mu$:
\begin{equation*}
	-\bm{\sigma}^\text{fr}:\nabla\mathbf{v}
	=-\bm{\sigma}^\text{fr}:\Def\mathbf{v}
	=-\tilde{\mu}\dive\mathbf{u}\dive\mathbf{v}-\mu\curl\mathbf{u}\cdot\curl\mathbf{v}-2\mu\dive\big[\nabla_\mathbf{v}\mathbf{u} - (\dive\mathbf{u})\mathbf{v}\big].
\end{equation*}
Therefore, concerning the variational constraint \eqref{eq:spatial-variational-constraint} we have everything that we need for its discretization, since, as we will see in a while, \emph{any} discretization of the term $\nabla_\mathbf{v}\mathbf{u} - (\dive\mathbf{u})\mathbf{v}$ will vanish in the derivation of the associated Euler-Lagrange equations thanks to lemma \ref{lem:integral-discrete-divergence}. 

This is not the case for the phenomenological constraint. In this case, we have the following according to \ref{cor:phenomenological-constraint-recasting}:
\begin{align*}
	-\bm{\sigma}^\text{fr}:\nabla\mathbf{u}&=-\bm{\sigma}^\text{fr}:\Def\mathbf{u}\\
	&=-\tilde{\mu}(\dive\mathbf{u})^2-\mu\|\curl\mathbf{u}\|^2-2\mu\dive\left[\left(\lie_\mathbf{u}\mathbf{u}^\flat-\frac{1}{2}\dif\|\mathbf{u}\|^2\right)^\sharp - (\dive\mathbf{u})\mathbf{u}\right],
\end{align*}
So the only term that needs attention is $(\dive\mathbf{u})\mathbf{u}$. But as this term sits in a divergence, we can view the $\dive\big((\dive\mathbf{u})\mathbf{u}\big)$ as the infinitesimal action of the vector field $\mathbf{u}$ on the density associated to $\dive\mathbf{u}$. Therefore the appropriate discretization will be $\dive(A)\bullet A$, and the same type of discretization will be used for the variational constraint\footnote{Although this modification is technically not needed, it is motivated by the fact that we should be able to deduce the variational constraint from the phenomenological one, as in the smooth case.}. Below we will denote by $\nabla_A A\in\LieAlgebraFont{d}(\mathbb{M})$ for any $A\in\LieAlgebraFont{d}(\mathbb{M})$ the discrete vector field defined by $(\nabla_A A)^\flat=\lie_A A^\flat-\frac{1}{2}\dif\llangle A^\flat,A\rrangle$. Consequently, as the multiplication of functions is discretized into the cell-by-cell product of discrete functions, the semi-discrete analogue of the quantity $-\bm{\sigma}^\text{fr}:\nabla\mathbf{v}$ is, for any cell $i\in\{1,\dots,N\}$:
\begin{equation}
	-\tilde{\mu}(\dive A)_i(\dive B)_i-\mu\big(\dif A^\flat\wedge\star\dif B^\flat)_i-2\mu(\dive X)_i+2\mu\big[(\dive A)\bullet B\big]_i,\label{eq:sigma-nabla-v-discretization}
\end{equation}
where $A\in\LieAlgebraFont{d}(\mathbb{M})$ approximates $\mathbf{u}$, $B\in\LieAlgebraFont{d}(\mathbb{M})$ approximates $\mathbf{v}$, and $X$ is any discretization of $\nabla_\mathbf{v}\mathbf{u}$, since anyway it will be discarded when we will derive the Euler-Lagrange equations associated to the spatial variational principle as mentioned before. As for the term $-\bm{\sigma}^\text{fr}:\nabla\mathbf{u}$ appearing in the spatial phenomenological constraint \eqref{eq:spatial-phenomenological-constraint}, its discretization reads:
\begin{equation}
	-\tilde{\mu}(\dive A)_i^2-\mu\big(\dif A^\flat\wedge\star\dif A^\flat\big)_i-2\mu(\dive\nabla_A A)_i+2\mu\left[(\dive A)\bullet A\right]_i.\label{eq:sigma-nabla-u-discretization}
\end{equation}

\begin{lemma}\label{lem:covariant-derivative}
	Let $A\in\mathcal{V}$. Then $\nabla_A A=\left[\lie_A A^\flat-\frac{1}{2}\dif\llangle A^\flat,A\rrangle\right]^\sharp\in\mathcal{V}$.
\end{lemma}

\begin{proof}
	$\nabla_A A\in\mathcal{V}$ is equivalent to $(\nabla_A A)_{ji}\Omega_{jj}+\Omega_{ii}(\nabla_A A)_{ij}=0$, for any $i\neq j$. This is straightforward from the definition \ref{def:sharp} of the sharp operator and the skew-symmetry of discrete 1-forms.
\end{proof}

%---------------------------------------------------------------------------------------------------

\subsection{Semi-discretization of the spatial variational principle}\label{ssec:semi-discrete-spatial}

We now discretize in space the spatial formalism presented in \ref{ssec:smooth-spatial}. The configuration space of the fluid, which was previously in the spatial formalism $\LieAlgebraFont{X}_0(M)\times\Den(M)\times\Den(M)\times\smooth{M}\times\Den(M)$ is now $\LieAlgebraFont{d}_0(\mathbb{M})\times\Den_d(\mathbb{M})\times\Den_d(\mathbb{M})\times\Omega_d^0(\mathbb{M})\times\Den_d(\mathbb{M})$, where $\mathbb{M}$ is a mesh approximating the fluid domain $M$. This concretely means that, after spatial discretization following the framework presented in \ref{sec:DEC}, the velocity vector field $\mathbf{u}$ becomes a discrete vector field $A$, the entropy density $s$ becomes a discrete entropy density $S$, the thermal displacement $\gamma$ becomes a discrete function $\Gamma$ and the variable $\sigma$, which is related to total entropy production, becomes a discrete density $\Sigma$. Both sparsity and compressibility constraints will be imposed as nonholonomic constraints in the variational principle. Note also that the mass density $\rho$, which evolves according to the advection equation, becomes a discrete density $D$.

The semi-discrete variational principle in the spatial formalism is obtained from the semi-discrete variational principle in the material formalism established in the previous section through a reduction akin to the Euler-Poincaré reduction recalled in \ref{ssec:smooth-Euler-Poincaré}. Denote by $q\in\mathsf{D}_0(\mathbb{M})$, $D_\text{ref}\in\Den_d(\mathbb{M})$, $S_\text{mat}\in\Den_d(\mathbb{M})$, $\Gamma_\text{mat}\in\Omega^0_d(\mathbb{M})$ and $\Sigma_\text{mat}\in\Den_d(\mathbb{M})$ for the discrete material variables. Then the spatial variables are $A=\dot{q}q^{-1}\in\LieAlgebraFont{d}_0(\mathbb{M})$, $D=D_\text{ref}\bullet q^{-1}\in\Den_d(\mathbb{M})$, $S=S_\text{mat}\bullet q^{-1}$, $\Gamma=\Gamma_\text{mat}\cdot q^{-1}$ and $\Sigma=\Sigma_\text{mat}\bullet q^{-1}$. Parts of the reduction of the phenomenological and variational constraints were given in section \ref{ssec:constraints-discretization}. We also have a new variable $B=\delta qq^{-1}\in\LieAlgebraFont{d}_0(\mathbb{M})$ associated to the variations $\delta q\in T\mathsf{D}(\mathbb{M})$ of the curve $q$. Note that we have $\dot{\Gamma}_\text{mat}\cdot q^{-1}=\dot{\Gamma}+\Gamma\cdot A$, $\delta\Gamma_\text{mat}\cdot q^{-1}=\delta\Gamma+\Gamma\cdot B$, $\dot{\Sigma}_\text{mat}\bullet q^{-1}=\dot{\Sigma}+\Sigma\bullet A$, $\delta\Sigma_\text{mat}\bullet q^{-1}=\delta\Sigma+\Sigma\bullet B$. We suppose the semi-discrete material Lagrangian $L_{D_\text{ref}}$ to be $\mathsf{D}_0(\mathbb{M})$-invariant, so that we can define a reduced Lagrangian $\lag(A,D,S):=\Lag_{D_\text{ref}}(q,\dot{q},S_\text{mat})$\footnote{The Lagrangian $\lag(A,D,S)$ is the spatial discretization of the smooth Lagrangian $\lag(\mathbf{u},\rho,s)$ given in \eqref{eq:spatial-Lagrangian} but we will explicit a possible discretization further in the text.}. With these definitions and remarks in mind, we can now state the spatial version of the semi-discrete variational principle for the Navier-Stokes-Fourier system.

\begin{mdframed}[style=box,frametitle={Semi-discrete variational principle for viscous heat conducting fluids, spatial version:}]
	Let $\mathbb{M}$ denote a mesh approximating the fluid domain $M$. The motion $t\mapsto \big(A(t),D(t),S(t),\Gamma(t),\Sigma(t)\big)\newline\in\LieAlgebraFont{d}_0(\mathbb{M})\times\Den_d(\mathbb{M})\times\Den_d(\mathbb{M})\times\Omega_d^0(\mathbb{M})\times\Den_d(\mathbb{M})$ is critical for the \emph{variational condition}:
	\begin{equation}\label{eq:semi-discrete-spatial-variational-condition}
		\delta\int_0^T\Big[\lag(A,D,S)+\underbrace{\big\langle S-\Sigma,\dot{\Gamma}+\Gamma\cdot A\big\rangle_0}_{\text{thermal coupling term}}\Big]\mathrm{d}t=0,
	\end{equation}
	subject to the \emph{phenomenological} and \emph{variational} constraints for any cell $i\in\{1,\dots,N\}$:
	\begin{align}
		\frac{\delta\lag}{\delta S_i}(\dot{\Sigma}+\Sigma\bullet A)_i&=-\tilde{\mu}(\dive A)_i^2-\mu\big(\dif A^\flat\wedge\star\dif A^\flat\big)_i-2\mu(\dive\nabla_A A)_i\notag\\
		&\quad\quad\quad+2\mu\left[(\dive A)\bullet A\right]_i+\big((\dot{\Gamma}+\Gamma\cdot A)\cdot J_S\big)_i-D_iR_i,\label{eq:semi-discrete-spatial-phenomenological-constraint}\\
		\frac{\delta\ell}{\delta S_i}(\delta\Sigma+\Sigma\bullet B)_i
		&=-\tilde{\mu}(\dive A)_i(\dive B)_i-\mu\big(\dif A^\flat\wedge\star\dif B^\flat\big)_i-2\mu(\dive X)_i\notag\\
		&\quad\quad\quad+2\mu\big[(\dive A)\bullet B\big]_i+\big((\delta\Gamma+\Gamma\cdot B)\cdot J_S\big)_i,\label{eq:semi-discrete-spatial-variational-constraint}
	\end{align}
	where $X$ is any discretization of $\nabla_\mathbf{v}\mathbf{u}$, $J_S\in\LieAlgebraFont{d}(\overline{\mathbb{M}})$ is a discrete entropy flux density and $D_iR_i$ represents the discrete external heat power supplied at the cell $C_i$; as well as both the \emph{sparsity} and \emph{compressibility} constraints:
	\[
		A\in\mathcal{S}\cap\mathcal{V},\quad B\in\mathcal{S}\cap\mathcal{V}.
	\]
	In \eqref{eq:semi-discrete-spatial-variational-condition}, the variation of the action functional has to be taken with respect to variations $$\delta A=\dot{B}-[A,B],\quad\delta D=-D\bullet B,$$ $\delta S$, $\delta\Sigma$ and $\delta\Gamma$, such that $B\in\LieAlgebraFont{d}(\mathbb{M})$ and $\delta\Gamma$ both vanish at the endpoints $t=0,T$, and such that $\frac{1}{2}(\delta\Gamma_i+\delta\Gamma_\partial)=0$ for all $C_i\in\mathbb{M}{\setminus}\mathbb{M}^\circ$.
\end{mdframed}

Now we have all the necessary ingredients to derive the semi-discrete Euler-Lagrange equations associated to this variational principle; it will yield a spatially discretized version of the Navier-Stokes-Fourier system of PDE. The derivation will tightly follow the one we did in section \ref{ssec:smooth-spatial}. First we compute the variations of the action functional in \eqref{eq:semi-discrete-spatial-variational-condition}:
\begin{align*}
	\delta\int_0^T\Big[\ell(A,D,S)&+\big\langle S-\Sigma,\dot{\Gamma}+\Gamma\cdot A\big\rangle_0\Big]\mathrm{d}t\\
	&=\int_0^T\Bigg[\underbrace{\llangle[\bigg]\frac{\delta\ell}{\delta A},\delta A\rrangle[\bigg]+\left\langle\frac{\delta\ell}{\delta D},\delta D\right\rangle_0}_{\circled{1}}+\underbrace{\left\langle\frac{\delta\ell}{\delta S},\delta S\right\rangle_0+\big\langle\delta S,\dot{\Gamma}+\Gamma\cdot A\big\rangle_0}_{\circled{2}}\\
	&\quad\quad\quad\quad\quad\quad\underbrace{-\big\langle\delta\Sigma,\dot{\Gamma}+\Gamma\cdot A\big\rangle_0}_{\circled{3}}+\underbrace{\big\langle S-\Sigma,\delta\dot{\Gamma}\big\rangle_0}_{\circled{5}}+\underbrace{\langle S-\Sigma,\Gamma\cdot\delta A\rangle_0}_{\circled{6}}+\underbrace{{\langle S-\Sigma,\delta\Gamma\cdot A\rangle_0}}_{\circled{5}}\Bigg]\mathrm{d}t.
\end{align*}
This expression is equal to zero for variations $\delta A$, $\delta D$, $\delta\Sigma$ and $\delta\Gamma$ satisfying all the mentioned variational constraints. Thanks to an integration by parts, section \ref{ssec:discrete-Lie-derivatives} and lemma \ref{lem:pairing-change}, the terms in \circled{1} are arranged as follows:
\begin{equation*}
	\int_0^T\Bigg[\llangle[\bigg]\frac{\delta\ell}{\delta A},\delta A\rrangle[\bigg]+\left\langle\frac{\delta\ell}{\delta D},\delta D\right\rangle_0\Bigg]\mathrm{d}t
	=-\int_0^T\llangle[\bigg]\frac{\mathrm{d}}{\mathrm{d}t}\frac{\delta\ell}{\delta A}+\Omega^{-1}\left[A^\mathsf{T},\Omega\frac{\delta\ell}{\delta A}\right]+D\frac{\delta\ell}{\delta D}^\mathsf{T},B\rrangle[\bigg]\,\mathrm{d}t.
\end{equation*}
The terms in \circled{2} yield:
\[
	\frac{\delta\ell}{\delta S}=-(\dot{\Gamma}+\Gamma\cdot A)=-\Theta,
\]
where $\Theta\in\Omega_d^0(\mathbb{M})$ is now the (semi-discrete) \emph{spatial temperature} (the reduction of the material temperature $\Theta_\text{mat}$). Then, we use the variational constraint \eqref{eq:semi-discrete-spatial-variational-constraint} and lemma \ref{lem:integral-action-density} for handling \circled{3}:
\begin{align*}
	-\big\langle\delta\Sigma,\dot{\Gamma}+\Gamma\bullet A\big\rangle_0
	&=\sum_i\Omega_{ii}\frac{\delta\ell}{\delta S_i}\delta\Sigma_i\\
	&=\sum_i\Omega_{ii}\Big[\Theta_i(\Sigma\bullet B)_i\underbrace{-\tilde{\mu}(\dive A)_i(\dive B)_i}_{\circled{1}}\underbrace{-\mu\big(\dif A^\flat\wedge\star\dif B^\flat\big)_i}_{\circled{1}}+\underbrace{\big((\delta\Gamma+\Gamma\cdot B)\cdot J_S\big)_i}_{\circled{4}}\Big]
\end{align*}
In the sum above, we obtain from lemma \ref{lem:div-adjoint} and proposition \ref{pro:discrete-curl-curl} that the terms in \circled{1} transform into 
\[
	\llangle\tilde{\mu}\,\dif\dive A,B\rrangle-\llangle\mu\Lambda(A^\flat),B\rrangle=\llangle[\big](\tilde{\mu}-\mu)\,\dif\dive A,B\rrangle[\big]+\llangle[\big]\mu\Delta A^\flat,B\rrangle[\big]
\]
according to definition \ref{def:discrete-Hodge-De-Rham-Laplacian}; they have to be added to the others terms we found before that couple with $B$. The term in \circled{4} is also treated with the help of lemma \ref{lem:div-adjoint}, taking into account the boundary condition on $\delta\Gamma$\footnote{In the case this boundary condition is not imposed, then the variational formalism yields the condition $J_S^\partial=0$, which is the discrete analogue of $\mathbf{j}_S\cdot\mathbf{n}=0$.}. We obtain:
\[
	-\big\langle\delta\Sigma,\dot{\Gamma}+\Gamma\bullet A\big\rangle_0
	=\underbrace{\langle\Theta,\Sigma\bullet B\rangle_0}_{\circled{8}}+\langle\tilde{\mu}\,\dif\dive A,B\rangle-\langle\mu\Lambda(A^\flat),B\rangle\underbrace{-\langle\dive J_S,\delta\Gamma\rangle_0}_{\circled{5}}\underbrace{-\langle\dive J_S,\Gamma\cdot B\rangle_0}_{\circled{7}}.
\]
After an integration by parts in time, the terms in \circled{5} combine to yield the following useful identity:
\begin{equation}\label{eq:semi-discrete-delta-Gamma}
	\big(\dot{S}-\dot{\Sigma}\big)+(S-\Sigma)\bullet A=-\dive J_S.
\end{equation}
We now handle the term \circled{6}. For this term, we compute:
\begin{align*}
	\int_0^T\langle S-\Sigma,\Gamma\cdot\delta A\rangle_0\,\mathrm{d}t
	&=\int_0^T\left[\big\langle S-\Sigma,\Gamma\cdot\dot{B}\big\rangle_0-\big\langle S-\Sigma,\Gamma\cdot[A,B]\big\rangle_0\right]\mathrm{d}t\\
	&=\int_0^T\left[-\big\langle\dot{S}-\dot{\Sigma},\Gamma\cdot B\big\rangle_0-\big\langle S-\Sigma,\dot{\Gamma}\cdot B\big\rangle_0-\big\langle S-\Sigma,\Gamma\cdot[A,B]\big\rangle_0\right]\mathrm{d}t\\
	&=\int_0^T\Big[\underbrace{-\big\langle(\dot{S}-\dot{\Sigma})+(S-\Sigma)\bullet A,\Gamma\cdot B\big\rangle_0}_{\circled{7}}\underbrace{-\langle S-\Sigma,\Theta\cdot B\rangle_0}_{\circled{8}}\Big]\mathrm{d}t.
\end{align*}
The terms in \circled{7} combine together and vanish thanks to \eqref{eq:semi-discrete-delta-Gamma}, and the terms in \circled{8} combine to give:
\[
	-\langle\Theta\cdot B,\Sigma\rangle_0-\langle S-\Sigma,\Theta\cdot B\rangle_0
	=-\langle\Theta\cdot B,S\rangle_0
	=\big\langle S\Theta^\mathsf{T},B\big\rangle,
\]
where at the end we used lemma \ref{lem:pairing-change}. Finally, using \eqref{eq:semi-discrete-delta-Gamma} to write $\dot{\Sigma}+\Sigma\bullet A=\dot{S}+S\bullet A+\dive J_S$, we obtain that the semi-discrete evolution equations for the Navier-Stokes-Fourier system are for any $i\in\{1,\dots,N\}$ and $j\in N(i)$:
\begin{equation}\label{eq:semi-discrete-abstract-Navier-Stokes-Fourier}
	\begin{cases}
		\displaystyle\mathbf{P}\left(-\frac{\mathrm{d}}{\mathrm{d}t}\frac{\delta\ell}{\delta A}-\Omega^{-1}\left[A^\mathsf{T},\Omega\frac{\delta\ell}{\delta A}\right]-D\frac{\delta\ell}{\delta D}^\mathsf{T}-S\frac{\delta\ell}{\delta S}^\mathsf{T}+(\tilde{\mu}-\mu)\,\dif\dive A+\mu{\Delta}A^\flat\right)_{ij}=0\\[7pt]
		\displaystyle\dot{D}+\Omega^{-1}A^\mathsf{T}\Omega D=0\\[7pt]
		\displaystyle\Theta_i\big(\dot{S}+S\bullet A+\dive J_S\big)_i\\[7pt]
		\displaystyle\quad=\tilde{\mu}(\dive A)^2_i+\mu\big(\dif A^\flat\wedge\star\dif A^\flat\big)_i+2\mu(\dive\nabla_A A)_i-2\mu\big[(\dive A)\bullet A\big]_i-(\Theta\cdot J_S)_i+D_iR_i
	\end{cases}
\end{equation}

The first equation corresponds to the semi-discrete momentum equation, the second to the semi-discrete continuity equation (which is obtained by reduction from $D=D_\text{ref}\cdot q^{-1}$) and the third to the semi-discrete entropy production equation. We can go further and rearrange the terms in order to make the Laplacian of the discrete temperature appear. We compute that: 
\[
	\Theta_i(\dive J_S)_i+(\Theta\cdot J_S)_i
	=2\Theta_i(J_S)_{ii}-\sum_{j=1}^{N+1}(J_S)_{ij}\Theta_j
	=-\sum_{\substack{j=1\\ j\neq i}}^{N+1}(\Theta_i+\Theta_j)(J_S)_{ij}.
\]
Then, we define, for any $j\in N(i)$, $j\neq i$:
\begin{equation}\label{eq:discrete-entropy-flux}
	(J_S)_{ij}:=\lambda\frac{\Theta_i-\Theta_j}{\Theta_i+\Theta_j}\frac{|h_{ij}|}{\Omega_{ii}|{\star h}_{ij}|},
\end{equation}
as a discretization of Fourier's law, so this way we obtain $\Theta_i(\dive J_S)_i+(\Theta\cdot J_S)_i=-\lambda(\Delta\Theta)_{i}$, $i\in\{1,\dots,N\}$, with the discrete Laplace-Beltrami operator defined in \ref{def:discrete-Laplace-Beltrami-Laplacian}, and where $\lambda$ is the thermal conductivity of the fluid. Then $(J_S)_{ii}:=-\sum_{j=1,\ j\neq i}^{N+1}(J_S)_{ij}$, so that $J_S\in\overline{\mathcal{S}}$, and note that we also have $J_S\in\overline{\mathcal{V}}$. In \eqref{eq:discrete-entropy-flux} we define the dual edge $\star h_{i\partial}$ (which appears when $j=\partial$) as the segment connecting the circumcenter of the cell $C_i$ to the midpoint of the shared edge with the environment $C_\partial$.
 
Therefore, the discrete entropy production equation (phenomenological constraint) takes the form:
\[
		\Theta_i\big(\dot{S}+S\bullet A\big)_i
		=\tilde{\mu}(\dive A)^2_i+\mu\big(\dif A^\flat\wedge\star\dif A^\flat\big)_i+2\mu(\dive\nabla_A A)_i-2\mu\big[(\dive A)\bullet A\big]_i+\lambda(\Delta\Theta)_i+D_iR_i.
\]

We now rewrite the previous semi-discrete equations for a particular Lagrangian. Remember from section \ref{ssec:smooth-spatial} that in the smooth setting, the spatial Lagrangian is given by \eqref{eq:spatial-Lagrangian}:
\[
	\ell(\mathbf{u},\rho,s)=\int_M\left[\frac{1}{2}\rho\langle\mathbf{u}^\flat,\mathbf{u}\rangle-\varepsilon(\rho,s)\right]\mathrm{d}x.
\]
Consequently one possible semi-discrete Lagrangian is given by:
\begin{equation}\label{eq:semi-discrete-spatial-Lagrangian}
	\ell(A,D,S)=\sum_{i=1}^N\Omega_{ii}\left[\frac{1}{2}D_i\sum_{j\in N(i)}A_{ij}^\flat A_{ij}-\widetilde{\varepsilon}(D_i,S_i)\right],
\end{equation}
where $\widetilde{\varepsilon}(D_i,S_i)$ is the real number approximating the average of the smooth internal energy $\varepsilon(\rho,s)$ on the cell $C_i$; we choose to make this value depend on the average cell value of $D$ as well as $S$. Recall also from section \ref{sec:DEC} that the operator $\mathbf{P}:\LieAlgebraFont{gl}(N)\to\Omega_d^1(\mathbb{M})$ is given for any $L\in\gl(N)$ by
\[
	\mathbf{P}(L)_{ij}=\frac{1}{2}(L_{ij}-L_{ji}-L_{ii}+L_{jj}),\quad
	\forall(i,j)\in\{1,\dots,N\}^2.
\]
We compute with the Lagrangian \eqref{eq:semi-discrete-spatial-Lagrangian}:
\[
	\frac{\delta\ell}{\delta A_{ij}}=D_iA_{ij}^\flat,\quad
	\frac{\delta\ell}{\delta D_{i}}=\sum_{j\in N(i)}\frac{1}{2}A_{ij}^\flat A_{ij}-\frac{\partial\widetilde{\varepsilon}}{\partial D_i},\quad
	\frac{\delta\ell}{\delta S_{i}}=-\frac{\partial\widetilde{\varepsilon}}{\partial S_i},
\]
and using the notation $\overline{D}_{ij}=\frac{1}{2}(D_i+D_j)$, we compute that the projected terms are:
\begin{align*}
	&\mathbf{P}\left(-\frac{\mathrm{d}}{\mathrm{d}t}\frac{\delta\ell}{\delta A}\right)_{ij}
	=-\frac{\mathrm{d}}{\mathrm{d}t}\mathbf{P}\left(\frac{\delta\ell}{\delta A}\right)_{ij}
	=-\frac{\mathrm{d}}{\mathrm{d}t}\left(\overline{D}_{ij}A_{ij}^\flat\right)
	=\overline{D\bullet A}_{ij}A_{ij}^\flat-\overline{D}_{ij}\dot{A}_{ij}^\flat,\\
	&\mathbf{P}\left(D\frac{\delta\ell}{\delta D}^\mathsf{T}\right)_{ij}
	=\overline{D}_{ij}\frac{\delta\ell}{\delta D}_{j}-\overline{D}_{ij}\frac{\delta\ell}{\delta D}_{i}
	=\overline{D}_{ij}\left(\sum_{k\in N(j)}\frac{1}{2}A_{jk}^\flat A_{jk}-\frac{\partial\widetilde{\varepsilon}}{\partial D_j}\right)-\overline{D}_{ij}\left(\sum_{k\in N(i)}\frac{1}{2}A_{ik}^\flat A_{ik}-\frac{\partial\widetilde{\varepsilon}}{\partial D_i}\right),\\
	&\mathbf{P}\left(S\frac{\delta\ell}{\delta S}^\mathsf{T}\right)_{ij}
	=\overline{S}_{ij}\frac{\delta\ell}{\delta S}_{j}-\overline{S}_{ij}\frac{\delta\ell}{\delta S}_{i}
	=-\overline{S}_{ij}\frac{\partial\widetilde{\varepsilon}}{\partial S_j}+\overline{S}_{ij}\frac{\partial\widetilde{\varepsilon}}{\partial S_i}.
\end{align*}
Note that the term with the Lie derivative has been already computed separately in proposition \ref{pro:Lie-derivative-1-form-density}. Therefore, the spatially discretized version of the system of PDE \eqref{eq:smooth-classical-Navier-Stokes-Fourier} is:
\begin{equation}\label{eq:semi-discrete-classical-Navier-Stokes-Fourier}
	\begin{cases}
		\displaystyle\frac{\mathrm{d}}{\mathrm{d}t}\big(\overline{D}_{ij}A_{ij}^\flat\big)+\lie_A(DA^\flat)_{ij}+\frac{1}{2}\overline{D}_{ij}\left(\sum_{k\in N(i)}\frac{1}{2}A_{ik}^\flat A_{ik}-\sum_{k\in N(j)}\frac{1}{2}A_{jk}^\flat A_{jk}\right)\\[0.5cm]
		\displaystyle\quad\quad=\overline{D}_{ij}\left(\frac{\partial\widetilde{\varepsilon}}{\partial D_i}-\frac{\partial\widetilde{\varepsilon}}{\partial D_j}\right)+\overline{S}_{ij}\left(\frac{\partial\widetilde{\varepsilon}}{\partial S_i}-\frac{\partial\widetilde{\varepsilon}}{\partial S_j}\right)+2(\tilde{\mu}-\mu)(A_{jj}-A_{ii})+\mu\big({\Delta}A^\flat\big)_{ij}\\[0.5cm]
		\displaystyle\dot{D}+\Omega^{-1}A^\mathsf{T}\Omega D=0\\[0.5cm]
		\displaystyle\Theta_i\big(\dot{S}+S\bullet A\big)_i\\[7pt]
		\displaystyle\quad\quad=\tilde{\mu}(\dive A)^2_i+\mu\big(\dif A^\flat\wedge\star\dif A^\flat\big)_i+2\mu(\dive\nabla_A A)_i-2\mu\big[(\dive A)\bullet A\big]_i+\lambda(\Delta\Theta)_i+D_iR_i
	\end{cases}
\end{equation}

\begin{remark}[Discrete pressure]
	In the first equation of \eqref{eq:semi-discrete-classical-Navier-Stokes-Fourier}, there is no discrete gradient of pressure, that is, an expression of the form $\dif P$ for some $P\in\Omega^0_d(\mathbb{M})$, although it could have been expected comparing \eqref{eq:semi-discrete-classical-Navier-Stokes-Fourier} with \eqref{eq:smooth-classical-Navier-Stokes-Fourier}. In \cite{PaMuToKaMaDe2011} a discrete pressure gradient was obtained at the semi-discrete level for the incompressible Euler equation because the nature of the discrete momenta was different and the authors resorted to the projector $\mathbf{
	R}$ explained in remark \ref{rq:projector-incompressible-Euler}, which naturally makes a discrete pressure gradient appears in the final discrete momentum equation. Here, the discrete momenta are different, and the equation is coupled to any $B\in\mathcal{V}\subset\LieAlgebraFont{d}(\mathbb{M})$ (see the derivation of the associated Euler-Lagrange equations above). Therefore we need to resort to the projector $\mathbf{P}$ and use the lemma \ref{cor:discrete-momenta-projection}. Nevertheless, observe that the expression
	\begin{equation}\label{eq:semi-discrete-pressure-gradient}
		\overline{D}_{ij}\left(\frac{\partial\widetilde{\varepsilon}}{\partial D_i}-\frac{\partial\widetilde{\varepsilon}}{\partial D_j}\right)+\overline{S}_{ij}\left(\frac{\partial\widetilde{\varepsilon}}{\partial S_i}-\frac{\partial\widetilde{\varepsilon}}{\partial S_j}\right)
	\end{equation}
	is an approximation of $\dif p=\rho\dif\frac{\partial\varepsilon}{\varepsilon\rho}+s\dif\frac{\partial\varepsilon}{\partial\rho}$, where the pressure $p$ is given by $p=\frac{\partial\varepsilon}{\partial\rho}\rho+\frac{\partial\varepsilon}{\partial s}s-\varepsilon$. Indeed, the integral $\int_{[\star C_i,\star C_j]}\dif p$ is approximated by $$\frac{D_i+D_j}{2}\int_{[\star C_i,\star C_j]}\dif\frac{\partial\varepsilon}{\partial\rho}+\frac{S_i+S_j}{2}\int_{[\star C_i,\star C_j]}\dif\frac{\partial\varepsilon}{\partial s},$$ which in turn is approximated by \eqref{eq:semi-discrete-pressure-gradient}. Note that the discrete pressure field is given by any discretization of the defining relation $p=\frac{\partial\varepsilon}{\partial\rho}\rho+\frac{\partial\varepsilon}{\partial s}s-\varepsilon$.
\end{remark}

We now focus on the discrete analogues of the balance of entropy and the balance of energy. On a cell $C_i$, we can formulate a local balance of entropy as follows. From the phenomenological constraint above we have:
\begin{align*}
	&\big(\dot{S}+S\bullet A+\dive J_S\big)_i\\
	&\quad=\frac{1}{\Theta_i}\left[\tilde{\mu}(\dive A)^2_i+\mu\big(\dif A^\flat\wedge\star\dif A^\flat\big)_i+2\mu(\dive\nabla_A A)_i-2\mu\big[(\dive A)\bullet A\big]_i-(\Theta\cdot J_S)_i\right]+\frac{D_i R_i}{\Theta_i}.
\end{align*}
Then, the global entropy balance on the mesh $\mathbb{M}$ is obtained by integration. According to lemma \ref{lem:integral-action-density} and lemma \ref{lem:integral-discrete-divergence-general} we obtain:
\begin{align*}
	&\sum_{i=1}^N\Omega_{ii}\dot{S}_i+\sum_{\substack{i\in N(\partial)\\ i\neq\partial}}\Omega_{ii}(J_S)^\partial_i\\
	&\quad=\sum_{i=1}^N\frac{\Omega_{ii}}{\Theta_i}\left[\tilde{\mu}(\dive A)^2_i+\mu\big(\dif A^\flat\wedge\star\dif A^\flat\big)_i+2\mu(\dive\nabla_A A)_i-2\mu\big[(\dive A)\bullet A\big]_i-(\Theta\cdot J_S)_i\right]\\
	&\quad\quad\quad+\sum_{i=1}^N\Omega_{ii}\frac{D_i R_i}{\Theta_i},
\end{align*}
which is the discrete analogue of corollary \ref{cor:smooth-spatial-entropy-balance}.

Now for the energy balance. We define the total energy of the semi-discrete system by:
\[
	\mathsf{E}(A,D,S)=\llangle[\bigg]\frac{\delta\ell}{\delta A},A\rrangle[\bigg]-\ell(A,D,S).
\]
We compute using the advection equation $\dot{D}+D\bullet A=0$ and lemma \ref{lem:pairing-change}:
\begin{align*}
	\frac{\mathrm{d}\mathsf{E}}{\mathrm{d}t}
	&=\llangle[\bigg]\frac{\mathrm{d}}{\mathrm{d}t}\frac{\delta\ell}{\delta A},A\rrangle[\bigg]-\left\langle\frac{\delta\ell}{\delta D},\dot{D}\right\rangle_0-\left\langle\frac{\delta\ell}{\delta S},\dot{S}\right\rangle_0\\
	&=\llangle[\bigg]\frac{\mathrm{d}}{\mathrm{d}t}\frac{\delta\ell}{\delta A}+D\frac{\delta\ell}{\delta D}^\mathsf{T},A\rrangle[\bigg]+\big\langle\Theta,\dot{S}\big\rangle_0.
\end{align*} 
Then using lemma \ref{lem:projection} and the phenomenological constraint we have:
\begin{align*}
	\frac{\mathrm{d}\mathsf{E}}{\mathrm{d}t}
	&=\llangle[\Bigg]\mathbf{P}\left(\frac{\mathrm{d}}{\mathrm{d}t}\frac{\delta\ell}{\delta A}+D\frac{\delta\ell}{\delta D}^\mathsf{T}\right),A\rrangle[\Bigg]\\
	&\quad\quad\quad+\sum_{i=1}^N\Omega_{ii}\Big[-\Theta_i(S\bullet A)_i+\tilde{\mu}(\dive A)^2_i+\mu\big(\dif A^\flat\wedge\star\dif A^\flat\big)_i+2\mu(\dive\nabla_A A)_i\\
	&\quad\quad\quad\quad\quad\quad-2\mu\big[(\dive A)\bullet A\big]_i+\lambda(\Delta\Theta)_i+D_iR_i\Big].
\end{align*}
Now observe that the term $2\mu(\dive\nabla_A A)_i$ as well as the term $2\mu\big[(\dive A)\bullet A\big]_i$ vanish thanks to lemmas \ref{lem:covariant-derivative}, \ref{lem:integral-discrete-divergence} and \ref{lem:integral-action-density}. 

Using the semi-discrete momentum equation in \eqref{eq:semi-discrete-abstract-Navier-Stokes-Fourier} then lemma \ref{lem:projection} we further obtain:
\begin{align*}
	\frac{\mathrm{d}\mathsf{E}}{\mathrm{d}t}
	&=\llangle[\bigg]-\Omega^{-1}\left[A^\mathsf{T},\Omega\frac{\delta\ell}{\delta A}\right]\underbrace{-S\frac{\delta\ell}{\delta S}^\mathsf{T}}_{\circled{1}}+\underbrace{(\tilde{\mu}-\mu)\,\dif\dive A+\mu{\Delta}A^\flat}_{\circled{2}},A\rrangle[\bigg]\\
	&\quad+\sum_{i=1}^N\Omega_{ii}\Big[\underbrace{-\Theta_i(S\bullet A)_i}_{\circled{1}}+\underbrace{\tilde{\mu}(\dive A)^2_i+\mu\big(\dif A^\flat\wedge\star\dif A^\flat\big)_i}_{\circled{2}}+\lambda(\Delta\Theta)_i+D_iR_i\Big].
\end{align*}
Now observe that terms in \circled{1} cancel each other thanks to lemma \ref{lem:pairing-change}, whereas terms in \circled{2} cancel each other thanks to lemma \ref{lem:div-adjoint} and proposition \ref{pro:discrete-curl-curl}. The remaining term is 
\[
	\llangle[\bigg]\ad_A^*\frac{\mathrm{d}}{\mathrm{d}t}\frac{\delta\ell}{\delta A},A\rrangle[\bigg]
\]
which clearly equals 0. Using lemma \ref{lem:div-adjoint} we obtain:
\[
	\frac{\mathrm{d}\mathsf{E}}{\mathrm{d}t}=
	-\sum_{\substack{i\in N(\partial)\\ i\neq \partial}}\Omega_{ii}\frac{\Theta_i+\Theta_\partial}{2}(J_S)_i^\partial+\sum_{i=1}^N\Omega_{ii}D_iR_i,
\]
which has to be compared to the corollary \ref{cor:smooth-spatial-energy-balance}.

%---------------------------------------------------------------------------------------------------

\section{A variational integrator for the Navier-Stokes-Fourier system}\label{sec:variational-discretization-Navier-Stokes-Fourier}

%---------------------------------------------------------------------------------------------------

\subsection{Variational integrators}

In computational geometric mechanics, variational integrators are geometric integrators based on the discretization of the variational principle from which comes the ordinary differential equation of interest, rather than the discretization of the equation itself. This subject is rich in geometric structures, which are the discrete counterpart of the geometric structures that one can find in traditional geometric mechanics, and brings desired and important qualitative features into numerical simulation in general.

\subsubsection{Variational integrators for finite-dimensional systems}

Variational integrators are numerical schemes that arise from a discrete version of Hamilton's principle, or Lagrange-d'Alembert's principle in the case external forces act on the system. These geometric integrators are thoroughly reviewed in \cite{MaWe2001}, we simply recall the broad idea hereafter. Let $Q$ be a configuration manifold and $\Lag:TQ\to\R$ be a Lagrangian. Given a time step $h$, $[0,T]$ is discretized into the sequence $t_k=kh$, $k\in\{0,\dots,n\}$. A curve $q$ in $Q$ is discretized into a sequence $q_d=(q_k)_{0\leq k\leq n}$, and a variation $\delta q$ of $q$ is discretized into a sequence $\delta q_d=(\delta q_k)_{0\leq k\leq n}$, such that $\delta q_k\in T_{q_k}Q$, for any $k\in\{0,\dots,n\}$. The Lagrangian $\Lag$ is discretized into a \emph{discrete Lagrangian} $\Lag_d:Q\times Q\to\R$ such that the following approximation holds:
\[
	\Lag_d(q_k,q_{k+1})\approx\int_{t_k}^{t_{k+1}}\Lag(q(t),\dot{q}(t))\,\mathrm{d}t,
\]
where the curve $t\mapsto q(t)$ is the solution of the Euler-Lagrange equations with endpoints $q_k$ and $q_{k+1}$. Usually this approximation is related to some numerical quadrature rule of the integral above. Then the discrete analogue of Hamilton's principle for the discrete action defined by
\[
	S_d(q_d)=\sum_{k=0}^{n-1}\Lag_d(q_k,q_{k+1})
\]
is $\delta S_d(q_d)\cdot\delta q_d=0$ for all variations $\delta q_d$ of $q_d$ with vanishing endpoints. After taking variations and applying a discrete integration by parts formula (change of indices), we obtain the \emph{discrete Euler-Lagrange} equations:
\begin{equation}\label{eq:discrete-Euler-Lagrange-equations}
\mathbf{D}_2\Lag_d(q_{k-1},q_k)+\mathbf{D}_1\Lag_d(q_k,q_{k+1})=0,\quad\forall k\in\{1,\dots,n-1\}.
\end{equation}
These equations define, under appropriate conditions, an algorithm which solves for $q_{k+1}$ knowing the two previous configuration variables $q_k$ and $q_{k-1}$. These integrators are symplectic, meaning that at each time step, the discrete symplectic form is preserved; for more details see \citet[section 1.3]{MaWe2001}.

\subsubsection{Variational integrators for finite-dimensional Euler-Poincaré systems}\label{ssec:discrete-Euler-Poincaré}

Let's now describe how variational integrators can be constructed in the case of Euler-Poincaré systems whose configuration space is a finite-dimensional Lie group $G$. The tangent space $TG$ is discretized into $G\times G$ as usual, and the analogue of the projection map $\bar{\pi}:TG\to TG/G\cong\LieAlgebraFont{g}$ is given by one of the two maps $\bar{\pi}^\pm:G\times G\to(G\times G)/G\cong G$ defined by $\bar{\pi}^+(g_0,g_1)=g_1g_0^{-1}$ and $\bar{\pi}^-(g_0,g_1)=g_0g_1^{-1}$. In what follows we choose to only work with $\bar{\pi}^+$ and will write $\Xi_k$ for $g_{k+1}g_k^{-1}$; this is the discrete analogue of $TG/G$ being identified with $\LieAlgebraFont{g}$ with the help of the right Maurer-Cartan form. Note however that this is just a matter of choice. Also note that in contrast to Euler-Poincar\'e reduction in the smooth setting, the discrete reduced tangent space is represented by the manifold $G$ rather than the vector space $\LieAlgebraFont{g}$.

Suppose that in addition to acting on the right of itself, $G$ acts on a parameter manifold $P$. Suppose that we have a discrete Lagrangian $L_{d,a_\text{ref}}:G\times G\to\R$, with $a_\text{ref}\in P$, which is $G_{a_\text{ref}}$-invariant; this means that for any $h_0$, $h_1$ and $g\in G$ we have:
\[
	\Lag_{d,a_\text{ref}}(h_0g,h_1g)=\Lag_{d,a_\text{ref}}(h_0,h_1).
\]
Then we define the \emph{reduced discrete Lagrangian} $\mathsf{L}_d:G\times G\times\Orb(a_\text{ref})\to\R$ for any $g_0$, $g_1\in G$ by:
\[
	\mathsf{\Lag}_d(g_1g_0^{-1},a_\text{ref}g_0^{-1})=\Lag_{d,a_\text{ref}}(g_0,g_1).
\]

Below is a summary of Euler-Poincaré reduction for discrete Euler-Poincaré systems as we defined above. It can be found in \citet{MaPeSh1999} without advected parameters and has been extended in \cite{CoGB2018} to include thermodynamics.

\begin{theorem}[\textbf{Discrete Euler-Poincar\'e reduction with advected parameters}]\label{thm:discrete-Euler-Poincaré-reduction}
	Let $G$ be a finite dimensional Lie group and denote by $\LieAlgebraFont{g}$ its Lie algebra. Suppose that $G$ acts on the right of a manifold $P$. For a fixed parameter $a_\text{ref}\in P$, let $\Lag_{d,a_\text{ref}}:G\times G\to\R$ be a discrete $G_{a_\text{ref}}$-invariant Lagrangian, and $\mathsf{\Lag_d}:G\times\Orb(a_\text{ref})\to\R$  be the associated reduced map. Then the following assertions are equivalent:
	\begin{enumerate}[label=(\arabic*), ref=\thetheorem.(\arabic*)]
		
		\item The discrete curve $g_d$ is critical for the discrete variational principle
		\begin{equation} 
			\delta\sum_{k=0}^{n-1}\Lag_{d,a_\text{ref}}(g_k,g_{k+1})=0.
		\end{equation}
	
		\item The discrete curve $g_d$ satisfies the discrete Euler-Lagrange equations \eqref{eq:discrete-Euler-Lagrange-equations}.
		
		\item\label{eq:reduced-discrete-variational-principle} The discrete curve $(\Xi_d,a_d)$ defined by $\Xi_k=g_{k+1}g_k^{-1}\in G$ and $a_k=a_\text{ref}g_k^{-1}\in\Orb(a_{\rm ref})$, is critical for the \emph{reduced discrete variational principle}
		\begin{equation*}			
			\delta\sum_{k=0}^{n-1}\mathsf{\Lag}_d(\Xi_k,a_k)=0,
		\end{equation*}
		subject to the discrete Euler-Poincar\'e constraints
		\[
			\delta\Xi_k=T_eR_{\Xi_{k+1}}\eta_{k+1}-T_eL_{\Xi_k}\eta_{k},\quad\delta a_k=-(\eta_k)_P(a_k),
		\]
		where $\eta_d$ is any discrete curve in $\LieAlgebraFont{g}$ with vanishing endpoints.
		
		\item\label{eq:discrete-Euler-Poincaré-equations} The discrete curve $(\Xi_d,a_d)$ is the solution of the \emph{discrete Euler-Poincar\'e equations}	
		\begin{empheq}[left=\empheqlbrace]{align}
			& \mathbf{D}_1\mathsf{\Lag}_d(\Xi_k,a_k)T_eL_{\Xi_k}=\mathbf{D}_1\mathsf{\Lag_d}(\Xi_{k-1},a_{k-1})T_eR_{\Xi_{k-1}}
			-\mathbf{J}\big(\mathbf{D}_2\mathsf{\Lag}_d(\Xi_k,a_k)\big)\label{eq:discrete-Euler-Poincaré-1}\\ 
			& a_{k+1}=a_k\Xi_k^{-1} \label{eq:discrete-advection-equation-1}
		\end{empheq}
	\end{enumerate}
\end{theorem}

\begin{proof}
	The equivalence between (1) and (2) comes directly from discrete variational mechanics, as in \ref{eq:discrete-Euler-Lagrange-equations}. Showing that (3) and (4) are equivalent is done as usual; note that the advection equation comes from the definition of $a_k$:
	\begin{equation*}
		a_{k+1}=a_\text{ref}g_{k+1}^{-1}=a_\text{ref}g_k^{-1}g_kg_{k+1}^{-1}=a_k\Xi_k^{-1}.
	\end{equation*}
	
	It remains to show that (1) and (3) are equivalent. Firstly, given a variation $\delta g_k$ of $g_k$,  we have by setting $\eta_k=\delta g_kg_k^{-1}=T_eR_{g_k^{-1}}\delta g_k$ that
	\begin{align*}
		\delta\Xi_k&=\delta g_{k+1}g_{k+1}^{-1}g_{k+1}g_k^{-1}-g_{k+1}g_k^{-1}\delta g_k g_k^{-1}\\
		&=T_eR_{\Xi_k}\eta_{k+1}-T_eL_{\Xi_k}\eta_{k}.
	\end{align*}
	Secondly we find that $\delta a_k=-(\eta_k)_M(a_k)$ using a computation similar to the one we did in the proof of the theorem \ref{thm:Euler-Poincaré-reduction}.
	
	Conversely, suppose that we are given the curve $\Xi_d$ and a variation $\delta\Xi_d$. We want to find a discrete curve $g_d$ and a discrete variation $\delta g_d$ of $g_d$ starting from curves $\Xi_d$ and $\eta_d$ as above. This is achieved by computing successively $g_{k+1}=\Xi_kg_k$, and by setting $\delta g_k=\eta_kg_k$. Since $\eta_d$ is arbitrary and zero at endpoints, $\delta g_d$ is arbitrary and zero at endpoints. We conclude using the right $G_{a_\text{ref}}$-invariance to obtain the variational principle (1) from the variational principle (3).
\end{proof}

As can be seen from Theorem \ref{thm:discrete-Euler-Poincaré-reduction}, the reduced discrete evolution equations now take place on $G\times \operatorname{Orb}(a_\text{ref})$. However, numerically speaking, solving differential equations on manifolds is more difficult that solving differential equations on vector spaces, as it is difficult to design a numerical scheme which ensures that the discrete evolution actually takes place in the manifold. We will write down discrete evolution equations on the Lie algebra $\LieAlgebraFont{g}$ instead. We will do so by using a \emph{group difference map} as introduced in \citet[section 4]{BRMa2008}. Simply put, these maps are approximations of the exponential map $\exp:\LieAlgebraFont{g}\to G$, but that still share its main algebraic properties.

\begin{definition}[\textbf{Group difference map}]
	Let $G$ be a Lie group and denote by $\LieAlgebraFont{g}$ its Lie algebra. A \emph{group difference map} is a local diffeomorphism $\tau:\LieAlgebraFont{g}\to G$ mapping a neighborhood $\mathcal{N}_0$ of $0\in\LieAlgebraFont{g}$ to a neighborhood of $e\in G$, and such that $\tau(0)=e$ and $\tau(\xi)^{-1}=\tau(-\xi)$, for any $\xi\in\mathcal{N}_0$.
\end{definition}

For the purpose of this paper, we will need the \emph{left} trivialized tangent of a group difference map, which is defined below.

\begin{definition}[\textbf{Left trivialized tangent of a group difference map}]
	Let $G$ be a Lie group, $\LieAlgebraFont{g}$ its Lie algebra and $\tau:\LieAlgebraFont{g}\to G$ a group difference map. The \emph{left trivialized tangent} of $\tau$ is the map $\mathrm{d}\tau:\LieAlgebraFont{g}\times\LieAlgebraFont{g}\to\LieAlgebraFont{g}$ defined by
	\[
		\mathbf{D}\tau(\xi)(\delta)=T_eL_{\tau(\xi)}\mathrm{d}\tau_\xi(\delta),
	\]
	for all $\xi$, $\delta\in\LieAlgebraFont{g}$. The \emph{inverse left trivialized tangent} of $\tau$ is the map $\mathrm{d}\tau^{-1}:\LieAlgebraFont{g}\times\LieAlgebraFont{g}\to\LieAlgebraFont{g}$ defined by
	\[
		\mathbf{D}\tau^{-1}(\tau(\xi))(\delta)=\mathrm{d}\tau^{-1}_\xi(T_eL_{\tau(-\xi)}\delta),
	\]
	for all $\xi$, $\delta\in\LieAlgebraFont{g}$. Thus $\mathrm{d}\tau_\xi(\mathrm{d}\tau^{-1}_\xi(\delta))=\delta$, for all $\xi$, $\delta\in\LieAlgebraFont{g}$. Note that $\mathrm{d}\tau$ and $\mathrm{d}\tau^{-1}$ are always linear in their second argument, but not necessarily in the first.
\end{definition}

\begin{remark}\label{rq:trivialized-tangents}
	Denote by $\mathrm{d}\tau^-$ the left trivialized tangent of the group difference map $\tau$. The right trivialized tangent $\mathrm{d}\tau^+$ of a $\tau$ is defined in a similar way:
	\[
		\mathbf{D}\tau(\xi)(\delta)=T_eR_{\tau(\xi)}\mathrm{d}\tau^+_\xi(\delta).
	\]
	Therefore the two trivialized tangents of $\tau$ are related by $\mathrm{d}\tau_\xi^+(\delta)=\Ad_{\tau(\xi)}\mathrm{d}\tau_\xi^-(\delta)$.
\end{remark}

\begin{proposition}[{\citet[section 4]{BRMa2008}}]\label{pro:group-difference-map}
	Let $G$ be a Lie group, $\LieAlgebraFont{g}$ its Lie algebra and $\tau:\LieAlgebraFont{g}\to G$ a group difference map. The left trivialized tangent $\mathrm{d}\tau$ of $\tau$ satisfies the following properties:
	\begin{enumerate}[label=(\arabic*), ref=\thetheorem.(\arabic*)]
		\item\label{eq:RTT1} $\mathrm{d}\tau_{-\xi}(\delta)=\Ad_{\tau(\xi)}\mathrm{d}\tau_{\xi}(\delta)$,
		\item\label{eq:RTT2} $\mathrm{d}\tau_{-\xi}^{-1}(\Ad_{\tau(\xi)}\delta)=\mathrm{d}\tau_\xi^{-1}(\delta)$.
	\end{enumerate}
\end{proposition}
\begin{proof}
	For the first property, for any $\xi\in\LieAlgebraFont{g}$ we have $\mu(\tau(\xi),\tau(-\xi))=e$, where $\mu$ denotes the multiplication law in $G$. Differentiating this relation we obtain for any $\delta\in\LieAlgebraFont{g}$:
	\[
		T_{\tau(\xi)}R_{\tau(-\xi)}\mathbf{D}\tau(\xi)(\delta)-T_{\tau(-\xi)}L_{\tau(\xi)}\mathbf{D}\tau(-\xi)(\delta)=0,
	\]
	and finally using the definition of the left trivialized tangent of $\tau$ we obtain
	\[
		\mathrm{d}\tau_{-\xi}(\delta)=T_{\tau(\xi)}R_{\tau(-\xi)}T_eL_{\tau(\xi)}\mathrm{d}\tau_{\xi}(\delta).
	\]
	The second property results from an application of the first one, with $\delta$ replaced by $\mathrm{d}\tau^{-1}_\xi(\delta)$.
\end{proof}

Approximations of the exponential map are available in terms of rational fractions, these are the well-known Pad\'e approximants of the exponential. The $(1,1)$ Pad\'e approximant of the exponential is also known as the \emph{Cayley map}, and is widely used in computational geometric mechanics when $G$ is a quadratic Lie group. For more details, see \citet[Section III.4.1 and IV.8.3]{HLW2006} and \citet[Section 4.6]{BRMa2008}. Also remark that a left trivialized tangent of a group difference map can be defined and used in a similar way, this is not related to which side $G$ acts on itself.

Let $\tau:\LieAlgebraFont{g}\to G$ be a group difference map and $h>0$ be a time step. We are now going to transport the reduced discrete variational formulation \ref{eq:reduced-discrete-variational-principle} and the associated discrete equations \ref{eq:discrete-Euler-Poincaré-equations} to the Lie algebra $\LieAlgebraFont{g}$ using this group difference map $\tau$. We define a new discrete Lagrangian $\lag_d:\LieAlgebraFont{g}\times\Orb(a_\text{ref})\to\R$ by
\[
	\lag_d(\xi_k,a_k)=\mathsf{\Lag}_d(\Xi_k,a_k),\text{ with }\xi_k:=\frac{1}{h}\tau^{-1}(\Xi_k).
\]
This can simply be considered as a change of variable in the Lagrangian, and the relation $g_{k+1}=\tau(h\xi_k)g_k$ amounts to discretize the reconstruction equation $\dot{g}(t)=\xi(t)g(t)=\frac{\mathrm{d}}{\mathrm{d}\varepsilon}\exp\big[\varepsilon\xi(t)\big]g(t)\big|_{\varepsilon=0}$, since the group difference map $\tau:\LieAlgebraFont{g}\to G$ is an approximation of the exponential map $\exp:\LieAlgebraFont{g}\to G$. This definition naturally extends to other quantities. Therefore, using the proposition \ref{pro:group-difference-map}, the discrete variational formulation \ref{eq:reduced-discrete-variational-principle} can be reformulated for a curve $(\xi_d,a_d)$ as:
\begin{equation*}
	\delta\sum_{k=0}^{n-1}\lag_d(\xi_k,a_k)=0,
\end{equation*} 
subject to the discrete Euler-Poincar\'e constraints
\begin{equation*}
	\delta\xi_k=\frac{1}{h}\delta\tau^{-1}(\Xi_k)=\frac{1}{h}\mathrm{d}\tau^{-1}_{h\xi_k}(\eta_{k+1})-\frac{1}{h}\mathrm{d}\tau^{-1}_{-h\xi_k}(\eta_k),\quad\delta a_k=-(\eta_k)_M(a_k),
\end{equation*} 
where $\eta_d$ is any discrete curve in $\LieAlgebraFont{g}$ with vanishing endpoints. From the discrete variational formulation above we get the following reformulation of the equations \ref{eq:discrete-Euler-Poincaré-equations}:
\begin{empheq}[left=\empheqlbrace]{align}
	& (\mathrm{d}\tau^{-1}_{-h\xi_k})^*\mathbf{D}_1\lag_d(\xi_k,a_k)=(\mathrm{d}\tau^{-1}_{h\xi_{k-1}})^*\mathbf{D}_1\lag_d(\xi_{k-1},a_{k-1})
	-h\mathbf{J}\big(\mathbf{D}_2\lag_d(\xi_k,a_k,S_k,S_{k+1})\big)\label{eq:discrete-Euler-Poincaré-2}\\
	& a_{k+1}=a_k\tau(-h\xi_k)\label{eq:discrete-advection-equation-2}
\end{empheq}

\begin{remark}
	For the kind of variational integrators that we have presented, if we want to define energy properly at the discrete level, then the time step $h$ has to be promoted to a full dynamic variable as well, meaning that the {discrete} Lagrangian of the system {is interpreted as being time-dependent} and thus the {discrete} Euler-Lagrange equations decompose into the usual dynamical part and an equation enforcing energy conservation in addition. See \cite{KMO1999} and \cite{LeDi2002} for such an approach in the case of mechanical systems without thermal effects.
\end{remark}

\begin{remark}
	A Kelvin-Noether theorem is available for this reformulation of the discrete Euler-Poincaré equations, see \citet[section 2.1.4]{DeGaGBZe2014} and \citet[section 2.3 and formula 2.14]{CoGB2018}.
\end{remark}

The variational discretization of Euler-Poincaré systems has been applied to the rigid body in \cite{BoRa2007}, the heavy top, incompressible fluids, incompressible magnetohydrodynamics, nematic liquid crystals and microstretch fluids in \cite{PaMuToKaMaDe2011}, rotating stratified fluids in \cite{DeGaGBZe2014}, anelastic and pseudo-incompressible flows in \cite{BaGB2017}. The spatial discretization of the compressible Euler equations has been done in \cite{BaGB2018}, but not the variational (time) discretization which can be obtained as a particular case of this paper.

\subsubsection{Variational integrators for simple systems with thermodynamics}

We first review the main points of the variational formalism for nonequilibrium thermodynamics of simple systems that has been introduced in \cite{GBYo2017a}. Recall that by a \emph{simple system} we mean a thermodynamical system for which we only need one entropy variable $S$ and a finite set of mechanical variables $(q,\dot{q})$ in order to describe entirely the state of the system. Given such a simple system, let $Q$ be the finite dimensional configuration manifold associated to the mechanical variables of the system. The Lagrangian of such a system is a map $\Lag:TQ\times\R\to\R$, $(q, \dot q, S) \mapsto \Lag(q, \dot q, S)$. We also have forces that act on the system: external forces $F^\text{ext}:TQ\times\R\to T^*Q$ that do not derive from a potential and friction forces $F^\text{fr}:TQ\times\R\to T^*Q$ which ultimately encode all the irreversible processes in the simple system and are responsible for internal entropy production. These maps are assumed to be fiber preserving, that is, $F^\text{ext}(q,\dot{q},S)$, $F^\text{fr}(q,\dot{q},S)\in T_q^*Q$ for any $q\in Q$, $\dot{q}\in T_q Q$ and $S\in\R$. It is also possible that the system exchanges heat with its environment, and we will denote by $P^\text{ext}_H:TQ\times\R\to\R$ the power due to heat transfer with the exterior of the system. A curve $(q,S)$ satisfies the \emph{variational formulation for nonequilibrium thermodynamics} if and only if it satisfies the \emph{variational condition}
\begin{equation}\label{eq:simple-variational-condition}
	\delta\int_0^T\Lag(q,\dot{q},S)\,\mathrm{d}t+\int_0^T\big\langle F^\text{ext}(q,\dot{q},S),\delta q\big\rangle\,\mathrm{d}t=0,
\end{equation}
for all variations $\delta q$ and $\delta S$ satisfying the \emph{variational constraint}
\begin{equation}\label{eq:simple-variational-constraint}
	\frac{\partial\Lag}{\partial S}(q,\dot{q},S)\delta S=\left\langle F^\text{fr}(q,\dot{q},S),\delta q\right\rangle,
\end{equation}
with $ \delta q(0)= \delta q(T)=0$ and if it also satisfies the \emph{phenomenological constraint}
\begin{equation}\label{eq:simple-phenomenological-constraint}
	\frac{\partial\Lag}{\partial S}(q,\dot{q},S)\dot{S}=\left\langle F^\text{fr}(q,\dot{q},S),\dot{q}\right\rangle-P^\text{ext}_H(q,\dot{q},S).
\end{equation}

The associated Euler-Lagrange equations are obtained by taking variations in the variational condition \eqref{eq:simple-variational-condition} and using the constraints \eqref{eq:simple-variational-constraint} and \eqref{eq:simple-phenomenological-constraint}, we obtain the following system of differential equations:
\begin{empheq}[left=\empheqlbrace]{align*}
	&\frac{\mathrm{d}}{\mathrm{d}t}\frac{\partial\Lag}{\partial\dot{q}}(q,\dot{q},S)-\frac{\partial\Lag}{\partial q}(q,\dot{q},S)=F^\text{ext}(q,\dot{q},S)+F^\text{fr}(q,\dot{q},S),\\[2mm]
	&\frac{\partial\Lag}{\partial S}(q,\dot{q},S)\dot{S}=\left\langle F^\text{fr}(q,\dot{q},S),\dot{q}\right\rangle-P_H^\text{ext}(q,\dot{q},S).
\end{empheq}
In a nutshell, the variational formalism reviewed above yields the time evolution equations for the thermomechanical system considered, in accordance with the axiomatic formulation of thermodynamics of Stueckelberg (\citet[chapter 1]{St1974}); see \citet[cection 3]{GBYo2017a} for more details.

We now review the variational discretization of thermodynamical simple systems as introduced in \citet[cection 3.1]{GBYo2017c} as it will give insights of what is done in the next section. The entropy curve $t\mapsto S(t)\in\R$ is discretized into a sequence $S_d=(S_k)_{0\leq k\leq N}$. The discrete Lagrangian is now a map $\Lag_d:Q\times Q\times\R\times\R\to\R$ such that:
\[
	\Lag_d(q_k,q_{k+1},S_k,S_{k+1})\approx\int_{t_k}^{t_{k+1}}\Lag(q(t),\dot{q}(t),S(t))\,\mathrm{d}t.
\]
The discrete counterparts of these forces are given by four maps $F^{\text{ext}+}_d$, $F^{\text{ext}-}_d$, $F^{\text{fr}+}_d$, $F^{\text{fr}-}_d:Q\times Q\times\R\times\R\to T^*Q$ such that the following approximation holds:
\begin{align*}
	\int_{t_k}^{t_{k+1}}\left\langle F^{\text{ext}}(q(t),\dot{q}(t),S(t)),\delta q(t)\right\rangle\mathrm{d}t
	\approx\big\langle & F_{d}^{\text{ext}-}(q_k,q_{k+1},S_k,S_{k+1}),\delta q_k\big\rangle\\
	&\quad+\big\langle F_{d}^{\text{ext}+}(q_k,q_{k+1},S_k,S_{k+1}),\delta q_{k+1}\big\rangle
\end{align*}
and similarly for $F^\text{fr}_d$. These discrete forces are required to be fiber-preserving in the sense that $\pi_{T^*Q}\circ F^{\text{ext}\pm}_d=\pi_Q^\pm$ and similarly for $F^\text{fr}_d$, with $\pi_{T^*Q}:T^*Q\to Q$ being the canonical projection and the maps $\pi_Q^\pm:Q\times Q\times\R\times\R\to Q$ being defined by $\pi_Q^-(q_0,q_1,S_0,S_1)=q_0$ and $\pi_Q^+(q_0,q_1,S_0,S_1)=q_1$. Concretely this means, for instance, that $F^{\text{ext}+}_d(q_0,q_1,S_0,S_1)\in T_{q_1}^*Q$. See also \citet[section 3.2.1]{MaWe2001} for a description of forces in the discrete setting.

We now need to discretize the phenomenological constraint. As stated in \citet[section 3]{GBYo2017c}, this is done with the help of a \emph{finite difference map} $\chi:Q\times Q\times\R\times\R\to TQ\times T\R$, a notion which was introduced in \citet[section 4]{McPe2006} for the development of variational integrators for systems with nonholonomic constraints. Essentially, such maps are directly responsible for the discretization of $q(t)$, $\dot{q}(t)$, $S(t)$ and $\dot{S}(t)$ in terms of $q_k$, $q_{k+1}$, $S_k$ and $S_{k+1}$, therefore their use is not limited to the discretization of constraints and they are also used for the discretization of the Lagrangian too. The phenomenological constraint can be seen as the zero-level set $\mathcal{C}$ of the map $P:T(Q\times\R)\to\R$ defined by:
\[
	P(q,\dot{q},S,\dot{S})=\frac{\partial\Lag}{\partial S}(q,\dot{q},S)\dot{S}-\big\langle F^\text{fr}(q,\dot{q},S),\dot{q}\big\rangle+P_H^\text{ext}(q,\dot{q},S),
\]
for any $(q,\dot{q})\in TQ$ and $(S,\dot{S})\in T\R$. The discrete counterpart of $\mathcal{C}$ is $\mathcal{C}_d$, defined as
\begin{equation}\label{eq:discrete-simple-phenomenological-constraint}
	\mathcal{C}_d=\chi^{-1}(\mathcal{C})\subset Q\times Q\times\R\times\R.
\end{equation}
Therefore, $\mathcal{C}_d$ can be seen as the zero-level set of the map $P_d=P\circ\chi:Q\times Q\times\R\times\R\to\R$\footnote{Note that the way in which the entropy is discretized in the discrete Lagrangian $\Lag_d$ is not necessarily the one used for the discrete phenomenological constraint $P_d$.}.

We can now state the discrete version of this variational formulation. A discrete curve $(q_d,S_d)$ satisfies the \textit{discrete variational formulation for nonequilibrium thermodynamics} if first it satisfies the \emph{discrete variational condition}
\begin{align*}
	\delta\sum_{k=0}^{N-1}\Lag_d(q_k&,q_{k+1},S_k,S_{k+1})\notag\\
	&+\sum_{k=0}^{N-1}\big\langle F^{\text{ext}-}_d(q_k,q_{k+1},S_k,S_{k+1}),\delta q_k\big\rangle
	+\big\langle F^{\text{ext}+}_d(q_k,q_{k+1},S_k,S_{k+1}),\delta q_{k+1}\big\rangle=0,
\end{align*}
for all variations $\delta q_d$ and $\delta S_d$ satisfying the \emph{discrete variational constraint}
\begin{align*}
	D_3\Lag_d(q_k,q_{k+1},S_k,S_{k+1})&\delta S_k+D_4\Lag_d(q_k,q_{k+1},S_k,S_{k+1})\delta S_{k+1}\notag\\
	&=\big\langle F^{\text{fr}-}_d(q_k,q_{k+1},S_k,S_{k+1}),\delta q_k\big\rangle
	+\big\langle F^{\text{fr}+}_d(q_k,q_{k+1},S_k,S_{k+1}),\delta q_{k+1}\big\rangle,
\end{align*}
for all $k\in\{0,\dots,N-1\}$, where $ \delta q_d$ vanishes at the endpoints, and if it also satisfies the \emph{discrete phenomenological constraint} 
\begin{equation*}
	P_d(q_k,q_{k+1},S_k,S_{k+1})=0,
\end{equation*}
for all $k\in\{0,\dots,N-1\}$. Taking variations and applying a discrete integration by parts (change of indices) yield the \emph{discrete equations for the thermodynamic of simple closed systems}:
\begin{empheq}[left=\empheqlbrace,right=,]{align*}
	&D_1\Lag_d(q_k,q_{k+1},S_k,S_{k+1})+D_2\Lag_d(q_{k-1},q_k,S_{k-1},S_k)\notag\\[2mm]
	&\quad\quad\quad+(F^{\text{ext}-}_d+F^{\text{fr}-}_d)(q_k,q_{k+1},S_k,S_{k+1})+(F^{\text{ext}+}_d+F^{\text{fr}+}_d)(q_{k-1},q_k,S_{k-1},S_k)=0\\[2mm]
	&P_d(q_k,q_{k+1},S_k,S_{k+1})=0
\end{empheq}
for all $k\in\{0,\dots,N-1\}$, see \cite{GBYo2017c}.

The most important point to recall is that these variational integrators qualitatively preserve the balance of energy of the system. Note that, although the phenomenological and variational constraints are tightly related in the smooth setting, there are not so much anymore in the discrete setting. 

In \cite{CoGB2018} variational integrators for Euler-Poincaré simple systems with thermodynamics were developed. For an Euler-Poincaré system on a finite-dimensional Lie group $G$, a finite difference map $\chi:G\times G\times\R\times\R\to T(G\times\R)$ is needed, where $\R$ is the space in which the entropy variable lives. As mentioned above, this finite difference map $\chi$ allows us to discretize the velocities $\dot{g}(t)$ and $\dot{S}(t)$ of a curve $(g,S)$ satisfying the Euler-Poincaré equations with thermodynamics. More precisely, we can write $\chi(g_k,g_{k+1},S_k,S_{k+1})=(g_k,\dot{g}_k,S_k,\dot{S}_k)$ with $\dot{g}_k\approx\dot{g}(t_k)$ and $\dot{S}_k\approx\dot{S}(t_k)$. Then, since in the smooth setting we have the relation $\dot{g}(t)=\xi(t)g(t)$, it is natural to define $\dot{g}_k:=\xi_k g_k$. Note that using this finite difference map is not problematic because usually we are ultimately much more interested in the evolution equations \eqref{eq:discrete-Euler-Poincaré-2} and \eqref{eq:discrete-advection-equation-2} after reduction and conversion to $\LieAlgebraFont{g}$ than the evolution equations \eqref{eq:discrete-Euler-Poincaré-1} and \eqref{eq:discrete-advection-equation-1} on $G$; see \citet[section 3.2]{CoGB2018}. Also note that the present article generalizes this approach to systems that are not simple, that is, the entropy of the system is not described by a single real number, and the Lie group acts on the entropy variable.

%---------------------------------------------------------------------------------------------------

\subsection{Variational integrator for the spatial variational principle}\label{ssec:discrete-spatial}

We now proceed to the reduction of the variational integrator presented in the previous section by extending what has been done in \citet[section 3.2]{CoGB2018}; this will yield a variational integrator in the spatial formalism for the Navier-Stokes-Fourier system. As in the variational discretization of any Euler-Poincaré system, the reduction is done in three stages that we recalled in \ref{ssec:discrete-Euler-Poincaré}. The first is writing down a variational principle on $\mathsf{D}_0(\mathbb{M})^2$ involving $q_k$ and $q_{k+1}\in\mathsf{D}_0(\mathbb{M})$; this was explained in the previous section. The second stage is to perform a reduction, then the reduced variational principle is expressed on $\mathsf{D}_0(\mathbb{M})$ in terms of $\Xi^k=q^{k+1}(q^k)^{-1}\in\mathsf{D}_0(\mathbb{M})$ (and the other reduced quantities). The last stage is akin to make a change of variable by defining $A^k:=\frac{1}{h}\tau^{-1}\big(q^{k+1}(q^k)^{-1}\big)$ where $\tau$ is a group difference map for $\mathsf{D}_0(\mathbb{M})$; it allows us to write a variational principle on $\LieAlgebraFont{d}(\mathbb{M})$ only. We will do the two last stages in one in this section. The reference mass density $D_\text{ref}\in\Den_d(\mathbb{M})$ becomes a dynamic mass density $D\in\Den_d(\mathbb{M})^{n+1}$ defined by $D^k=D_\text{ref}\bullet(q^k)^{-1}$. If we denote by $S_\text{mat}\in\Den_d(\mathbb{M})^{n+1}$, $\Gamma_\text{mat}\in\Omega_d^0(\mathbb{M})^{n+1}$ and $\Sigma_\text{mat}\in\Den_d(\mathbb{M})^{n+1}$ the other discrete curves that represent discrete material variables, then we have new discrete curves $S\in\Den_d(\mathbb{M})^{n+1}$, $\Gamma\in\Omega_d^0(\mathbb{M})^{n+1}$ and $\Sigma\in\Den_d(\mathbb{M})^{n+1}$ defined by $S^k=S_\text{mat}^k\bullet(q^k)^{-1}$, $\Gamma^k=\Gamma_\text{mat}^k\cdot(q^k)^{-1}$ and $\Sigma^k=\Sigma_\text{mat}^k\bullet(q^k)^{-1}$. The entropy flux density $J_{S,\text{mat}}$ and the external heat supply $R_\text{mat}$ are transformed into $J_S^k=J_{S,\text{mat}}^k(q^k)^{-1}$ and $R^k=R^k_\text{mat}\bullet(q^k)^{-1}$. We also have a new variable $B\in\LieAlgebraFont{d}(\mathbb{M})^{n+1}$ defined by $B^k=\delta q^k(q^k)^{-1}$, where $\delta q_d$ is a discrete variation of the discrete curve $q_d$. Note that $\Sigma_\text{mat}^{k+1}-\Sigma_\text{mat}^k$ is reduced to $\Sigma^{k+1}\bullet q^{k+1}(q^k)^{-1}-\Sigma^k=\Sigma^{k+1}\bullet\tau(h A^k)-\Sigma^k$, $\Gamma_\text{mat}^{k+1}-\Gamma_\text{mat}^k$ to $\Gamma^{k+1}\cdot\tau(hA^k)-\Gamma^k$, $\delta\Sigma_\text{mat}^k$ to $\delta\Sigma^k+\Sigma^k\bullet B^k$ and $\delta\Gamma_\text{mat}^k$ to $\delta\Gamma^k+\Gamma^k\cdot B^k$. We will also suppose that $L_{d,D_\text{ref}}$ is $\mathsf{D}_0(\mathbb{M})$-invariant so that we can define a reduced discrete Lagrangian by $\lag_d(A^k,D^k,S^k,S^{k+1})=\Lag_{d,\text{ref}}(q^k,q^{k+1},S_\text{mat}^k,S_\text{mat}^{k+1})$.

\begin{lemma}\label{lem:commutation-actions}
	Let $G$ be a Lie group, denote by $\LieAlgebraFont{g}$ its Lie algebra, and let $V$ be a manifold on which $G$ acts on the right. Then we have $v\cdot T_eR_g(\xi)=(v\cdot\xi)\cdot g$ for any $g\in G$, $\xi\in\LieAlgebraFont{g}$ and $v\in V$.
\end{lemma}

\begin{proof}
	We compute:
	\[
		(v\cdot\xi)\cdot g =\frac{\mathrm{d}}{\mathrm{d}\varepsilon}\bigg|_{\varepsilon=0}\left[v\cdot\exp(\varepsilon\xi)\right]\cdot g
		=\frac{\mathrm{d}}{\mathrm{d}\varepsilon}\bigg|_{\varepsilon=0}v\cdot\left[\exp(\varepsilon\xi)\cdot g\right]
		=v\cdot T_eR_g(\xi).
	\]
\end{proof}

Below we will denote by $\big[\mathbf{D}_Z\lag_d\big]^k$ the quantity $\mathbf{D}_Z\lag_d(A^k,D^k,S^k,S^{k+1})$, where $Z$ is one of the four variables of $\lag_d(A,D,S_0,S_1)$. In this article, we will limit ourselves to the simplest discrete Lagrangian, $\lag_d(A,D,S_0,S_1)=h\lag(A,D,S_0)$. However, other discrete Lagrangians can be considered as well as other discretizations of the thermal coupling term.

\begin{mdframed}[style=box,frametitle={Variational integrator for viscous heat conducting fluids, spatial version:}]
	Let $\mathbb{M}$ denote a mesh approximating the fluid domain $M$ and a sequence $(t_k)_{0\leq k\leq n}$ discretizing the time interval. The discrete motion $\big(A_d,D_d,S_d,\Gamma_d,\Sigma_d\big)\in\LieAlgebraFont{d}_0(\mathbb{M})^{n+1}\times\Den_d(\mathbb{M})^{n+1}\times\Den_d(\mathbb{M})^{n+1}\times\Omega_d^0(\mathbb{M})^{n+1}\times\Den_d(\mathbb{M})^{n+1}$ is critical for the \emph{discrete variational condition}:
	\begin{equation}\label{eq:discrete-spatial-variational-condition}
		\delta\sum_{k=0}^{n-1}\Big[\lag_d(A^k,D^k,S^k,S^{k+1})+\big\langle S^k-\Sigma^k,\Gamma^{k+1}\cdot\tau(hA^k)-\Gamma^k\big\rangle_0\Big]=0,
	\end{equation}
	subject to the \emph{discrete phenomenological} and \emph{variational} constraints for any cell $i\in\{1,\dots,N\}$ and time step $k\in\{0,\dots,n-1\}$:
	\begin{align}
		&\Theta^k_i\big(\Sigma^{k+1}\bullet\tau(hA^k)-\Sigma^k\big)_i\notag\\
		&\quad=h\tilde{\mu}(\dive A^k)_i^2+h\mu\big(\dif A^{k\flat}\wedge\star\dif A^{k\flat}\big)_i+2h\mu(\dive\nabla_{A^k} A^k)_i-2h\mu\left[(\dive A^k)\bullet A^k\right]_i\notag\\
		&\hspace{7cm}-\Big[\big(\Gamma^{k+1}\cdot\tau(hA^k)-\Gamma^k\big)\cdot J_S^k\Big]_i+hD_i^kR_i^k,\label{eq:discrete-spatial-phenomenological-constraint}\\
		&\big[\mathbf{D}_{S_{0}}\lag_d\big]^k_i\big(\delta{\Sigma}^k+\Sigma^k\bullet B^k\big)_i+\big[\mathbf{D}_{S_{1}}\lag_d\big]_i^{k}\big(\delta{\Sigma}^{k+1}+\Sigma^{k+1}\bullet B^{k+1}\big)_i\notag\\
		&\quad=-h\tilde{\mu}(\dive A^k)_i(\dive B^k)_i-h\mu\big(\dif A^{k\flat}\wedge\star\dif B^{k\flat}\big)_i-2h\mu(\dive X)_i+2h\mu\left[(\dive A^k)\bullet B^k\right]_i\notag\\
		&\hspace{9.7cm}+h\big((\delta{\Gamma}^k+\Gamma^k\cdot B^k)\cdot J_S^k\big)_i,\label{eq:discrete-spatial-variational-constraint}
	\end{align}
	where $X$ is any discretization of $\nabla_\mathbf{v}\mathbf{u}$, $J_S\in\LieAlgebraFont{d}(\overline{\mathbb{M}})$ is a discrete entropy flux density and $D_iR_i$ represents the discrete external heat power supplied at the cell $C_i$; as well as both the \emph{sparsity} and \emph{compressibility} constraints:
	\[
		A^k\in\mathcal{S}\cap\mathcal{V},\quad B^k\in\mathcal{S}\cap\mathcal{V}.
	\]
	In \eqref{eq:discrete-spatial-variational-condition}, the variation of the action functional has to be taken with respect to variations 
	\[
		\delta A^k=\frac{1}{h}\mathrm{d}\tau^{-1}_{hA^k}\big(B^{k+1}\big)-\frac{1}{h}\mathrm{d}\tau^{-1}_{-hA^k}\big(B^k\big),\quad\delta D^k=-D^k\bullet B^k,
	\] 
	$\delta S_d$, $\delta\Sigma_d$ and $\delta\Gamma_d$, such that $B_d\in\LieAlgebraFont{d}(\mathbb{M})^{n+1}$ and $\delta\Gamma_d\in\Omega_d^0(\mathbb{M})^{n+1}$ both vanish at the endpoints $t=0,T$, and such that $\frac{1}{2}(\delta\Gamma^k_i+\delta\Gamma^k_\partial)=0$ for all $C_i\in\mathbb{M}{\setminus}\mathbb{M}^\circ$ and time step $k$. Here $\mathrm{d}\tau^{-1}$ denotes the inverse \emph{left} trivialized tangent of the group difference map $\tau$.
\end{mdframed}

We now derive the associated discrete Euler-Lagrange equations. Taking variations of the discrete action functional in \eqref{eq:discrete-spatial-variational-condition} we obtain:
\begin{align*}
	&\sum_{k=0}^{n-1}\Bigg[\underbrace{\llangle[\Big]\big[\mathbf{D}_A\lag_d]^k,\delta A^k\rrangle[\Big]+\Big\langle\big[\mathbf{D}_D\lag_d]^k,\delta D^k\Big\rangle_0}_{\circled{1}}+\underbrace{\Big\langle\big[\mathbf{D}_{S_0}\lag_d\big]^k,\delta S^k\Big\rangle_0+\Big\langle\big[\mathbf{D}_{S_1}\lag_d\big]^{k},\delta S^{k+1}\Big\rangle_0}_{\circled{2}}\\
	&\quad+\underbrace{\big\langle\delta S^k,\Gamma^{k+1}\cdot\tau(hA^k)-\Gamma^k\big\rangle_0}_{\circled{2}}\underbrace{-\big\langle\delta\Sigma^k,\Gamma^{k+1}\cdot\tau(hA^k)-\Gamma^k\big\rangle_0}_{\circled{3}}+\underbrace{\big\langle S^k-\Sigma^k,\delta\Gamma^{k+1}\cdot\tau(hA^k)\big\rangle_0}_{\circled{4}}\\
	&\quad+\underbrace{\big\langle S^k-\Sigma^k,\Gamma^{k+1}\cdot\mathbf{D}\tau(hA^k)(h\delta A^k)\big\rangle_0}_{\circled{5}}\underbrace{-\big\langle S^k-\Sigma^k,\delta\Gamma^k\big\rangle_0}_{\circled{4}}\Bigg]=0.
\end{align*}
Using an integration by parts (change of indices) and lemma \ref{lem:pairing-change}, terms in \circled{1} are equal to:
\begin{align*}
	\llangle[\Big]\frac{1}{h}(\mathrm{d}\tau^{-1}_{-hA^{k-1}})^*\big[\mathbf{D}_A\lag_d\big]^{k-1}-\frac{1}{h}(\mathrm{d}\tau^{-1}_{hA^k})^*\big[\mathbf{D}_A\lag_d\big]^k-D^k[\mathbf{D}_D\lag_d\big]^{k\mathsf{T}},B^k\rrangle[\Big].
\end{align*}

Using another change of indices, we obtain from terms in \circled{2}:
\[
	\big[\mathbf{D}_{S_0}\lag_d\big]^k+\big[\mathbf{D}_{S_1}\lag_d\big]^{k-1}=-\big(\Gamma^{k+1}\cdot\tau(hA^k)-\Gamma^k\big)=:-h\Theta^k.
\]

Using this definition and a change of indices as well as lemmas \ref{lem:integral-action-density} and \ref{lem:integral-discrete-divergence-general}, we rearrange the terms in \circled{3} with the help of the variational constraint as follows:
\begin{align*}
	\sum_{k=0}^{n-1}-\big\langle\delta\Sigma^k,h\Theta^k\big\rangle_0
	=\sum_{k=0}^{n-1}\sum_{i=1}^N&\Omega_{ii}\bigg[\underbrace{-\big[\mathbf{D}_{S_{0}}\lag_d\big]^k_i\big(\Sigma^k\bullet B^k\big)_i-\big[\mathbf{D}_{S_{1}}\lag_d\big]_i^{k}\big(\Sigma^{k+1}\bullet B^{k+1}\big)_i}_{\circled{5}}\\
	&\underbrace{-h\tilde{\mu}(\dive A^k)_i(\dive B^k)_i-h\mu\big(\dif A^{k\flat}\wedge\star\dif B^{k\flat}\big)_i}_{\circled{1}}+\underbrace{h(\delta\Gamma^k\cdot J_S^k)_i}_{\circled{4}}+\underbrace{h\big((\Gamma^k\cdot B^k)\cdot J_S^k\big)_i}_{\circled{5}}\bigg].
\end{align*}

The terms $-h\tilde{\mu}(\dive A^k)_i(\dive B^k)_i-h\mu\big(\dif A^{k\flat}\wedge\star\dif B^{k\flat}\big)_i$ will have to be added to other terms in \circled{1} found previously, through the use of lemma \ref{lem:div-adjoint} and proposition \ref{pro:discrete-curl-curl}. Then collecting terms in \circled{4} we obtain thanks to a change of indices the following relation, similar to \eqref{eq:semi-discrete-delta-Gamma}:
\begin{equation}\label{eq:discrete-delta-Gamma}
	(S^{k-1}-\Sigma^{k-1})\bullet\tau(-hA^{k-1})-(S^k-\Sigma^k)=h\dive J_S^k.
\end{equation}

Now, using the properties of the inverse left trivialized tangent of the group difference map $\tau$ we compute that:
\begin{align*}
	\big\langle S^k-\Sigma^k,\Gamma^{k+1}\cdot&\mathbf{D}\tau(hA^k)(h\delta A^k)\big\rangle_0\\
	&=\Big\langle S^k-\Sigma^k,\Gamma^{k+1}\cdot\Big[T_eL_{\tau(hA^k)}\mathrm{d}\tau_{hA^k}\big(\mathrm{d}\tau^{-1}_{-hA^k}(B^{k+1})-\mathrm{d}\tau^{-1}_{hA^k}(B^k)\big)\Big]\Big\rangle_0\\
	&=\big\langle S^k-\Sigma^k,\Gamma^{k+1}\cdot\big[T_eL_{\tau(hA_k)}\Ad_{\tau(-hA^k)}(B^{k+1})-T_eL_{\tau(hA^k)}(B^k)\big]\big\rangle_0\\
	&=\big\langle S^k-\Sigma^k,\Gamma^{k+1}\cdot\big(B^{k+1}\tau(hA^k)-\tau(hA^k)B^k\big)\big\rangle_0,	
\end{align*}
which combined with other terms in \circled{5}, yield, using lemma \ref{lem:commutation-actions}, the definition of temperature, as well as:
\begin{align*}
	\sum_{k=0}^{n-1}\Big[\big\langle(S^{k-1}-\Sigma^{k-1})\bullet\tau(-hA^{k-1}),\Gamma^k\cdot B^k\big\rangle_0&-\big\langle S^k-\Sigma^k,(h\Theta^k+\Gamma^k)\cdot B^k\big\rangle_0\\
	&+\langle\Sigma^k\bullet B^k,h\Theta^k\rangle_0-h\langle\dive J_S^k,\Gamma^k\cdot B^k\rangle_0\Big].
\end{align*}
Then using the relation \eqref{eq:discrete-delta-Gamma} and lemma \ref{lem:two-actions}, we finally obtain that the above sum is equal to 
\[
	\sum_{k=0}^{n-1}\llangle[\big]hS^k\Theta^{k\mathsf{T}},B^k\rrangle[\big].
\] 
The fully discrete Navier-Stokes-Fourier equations are therefore (see section \ref{ssec:discrete-momenta} for the definition and use of the projection $\mathbf{P}$):
\begin{equation}\label{eq:discrete-abstract-Navier-Stokes-Fourier}
\begin{cases}
	\displaystyle\mathbf{P}\bigg(\frac{1}{h}(\mathrm{d}\tau^{-1}_{hA^k})^*D^kA^{k\flat}-\frac{1}{h}(\mathrm{d}\tau^{-1}_{-hA^{k-1}})^*D^{k-1}A^{(k-1)\flat}\bigg)_{ij}\\
	\displaystyle\hspace{5cm}+\frac{1}{2}\overline{D}_{ij}^k\left(\sum_{p\in N(i)}\frac{1}{2}A_{ip}^{k\flat} A_{ip}^k-\sum_{p\in N(j)}\frac{1}{2}A_{jp}^{k\flat} A_{jp}^k\right)\\[0.5cm]
	\displaystyle\quad\quad=\overline{D}_{ij}^k\left(\frac{\partial\widetilde{\varepsilon}}{\partial D_i^k}-\frac{\partial\widetilde{\varepsilon}}{\partial D_j^k}\right)+\overline{S}_{ij}^k\left(\frac{\partial\widetilde{\varepsilon}}{\partial S_i^k}-\frac{\partial\widetilde{\varepsilon}}{\partial S_j^k}\right)+2(\tilde{\mu}-\mu)(A_{jj}^k-A_{ii}^k)+\mu\big({\Delta}A^{k\flat}\big)_{ij}\\[7pt]
	\displaystyle D^{k+1}=D^k\bullet\tau(-hA^k)=\Omega^{-1}\tau(-hA^k)^\mathsf{T}\Omega D^k\\[7pt]
	\displaystyle\Theta^k_i\big(S^{k+1}\bullet\tau(hA^k)-S^k+h\dive J_S^{k+1}\bullet\tau(hA^k)\big)_i\\[7pt]
	\displaystyle\quad\quad=h\tilde{\mu}(\dive A^k)_i^2+h\mu\big(\dif A^{k\flat}\wedge\star\dif A^{k\flat}\big)_i+2h\mu(\dive\nabla_{A^k} A^k)_i-2h\mu\left[(\dive A^k)\bullet A^k\right]_i\\[7pt]
	\displaystyle\hspace{3cm}-\big[\big(\Gamma^{k+1}\cdot\tau(hA^k)-\Gamma^k\big)\cdot J_S^k\big]_i+hD_i^kR_i^k
\end{cases}
\end{equation}
to which we need to add the Dirichlet boundary conditions on the temperature $\Theta_i^k=\Theta_\partial$ for any cell $C_i$ in $\overline{\mathbb{M}}\setminus\mathbb{M}^\circ$, and also the no-slip boundary conditions for the discrete vector field $A^k$, at each time step.

Concerning the choice of group difference map $\tau$, contrary to the incompressible case, there is no reason to use the Cayley map as $\LieAlgebraFont{d}_0(\mathbb{M})$ is not a quadratic Lie algebra. Moreover, the discrete velocity field $A$ has to satisfy the compressibility constraint $\mathcal{V}$, which means in particular that $A=L+D$ with $L$ a $\Omega$-skew symmetric matrix and $D$ a diagonal matrix, and in this case $I\pm A$ is possibly not invertible, preventing us to use the Cayley map in this particular case.

%--------------------------------------------------------------------------------------------------

\section*{Conclusion and perspectives}
\addcontentsline{toc}{section}{Conclusion and perspectives}

This article opens new perspectives concerning the geometric integration of fluids; we would like to give a few of them. On the \emph{geometrical side}, the reduction of the structure-preserving property established in \cite{GBYo2018} is a topic of high importance, as it replaces the preservation of the symplectic form, in both smooth and discrete settings. The geometric structure, in terms of a variational principle, for interactions between a fluid and a rigid body, or between two fluids, is an active subject of research, see \cite{JaVa2014} for instance. It would be interesting to see how the framework developed in this paper can be enlarged to handle the interactions with a rigid body, for instance with a more sophisticated model yet similar to the one developed in \citet[section 3]{CoGB2018}. Another research direction is the incorporation of stochastics into geophysical fluid flows, see for instance \cite{GBHo2018}. Indeed, measurement campaigns performed by satellites have revealed the stochastic behavior of drifting devices in the ocean. It would be interesting to investigate discrete variational integrators based on the stochastic models developed in \cite{GBHo2018}. It could be then also for the development of a Kalman filter that would use this stochastic variational integrator for making predictions (see \cite{ScFlMu2017} for an account of the use of variational integrators as predictors in Kalman filters).

More on the \emph{modeling side}, it would be interesting to understand the discretization of other irreversible processes, for instance in magnetohydrodynamics, and also to see if the variational integrator developed in this paper can be used to optimally control the equilibrium state of a plasma, which is of great interest because of the special properties of geometric integrators. Also, multiple irreversible processes are included in the atmospheric model developed in \cite{GB2018}; it would be important to discretize this model and implement it as a circulation model in mind and see the benefits of geometric integration on long runs of the integrator. Discretization of complex fluids is also a possible direction of research. 

At last, a lot remains to be done on the \emph{numerical side}. Other discretizations of the thermal coupling term have to be considered, as well as others methods for vector field discretization and reconstruction. The numerical analysis (stability, convergence) of this variational integrator, starting from the one in the original work of \cite{Pa2009} for ideal fluids, is a subject of uttermost importance. This subject will bring many standard tools of numerical analysis into computational geometric mechanics and forge interesting partnerships between these fields. Finally, it would be interesting to benchmark this integrator against known numerical schemes for the Navier-Stokes-Fourier system qualitatively, but also in terms of (empirical) convergence, and computational and memory cost.

\bigbreak\textbf{Acknowledgments:} The authors were financed by the ANR project GEOMFLUID (ANR-14-CE23-0002). The first author thanks S. Shamekh for helpful discussions concerning Rayleigh-Bénard convection.

\end{document}